\documentclass[english,x11names]{amsart}

\PassOptionsToPackage{hyphens}{url}
\usepackage[top=1.4in, bottom=1.3in, left=1.4in, right=1.4in]{geometry}
\usepackage{amsxtra}
\usepackage{mathtools}
\usepackage{amsfonts}
\usepackage[toc,page]{appendix}
\usepackage{amsthm}
\usepackage[most]{tcolorbox}
\usepackage{graphicx}
\usepackage{amssymb}
\usepackage[protrusion=true, tracking=true, expansion=true, final]{microtype}
\usepackage{tikz-cd}
\usepackage{textcomp}
\usepackage{fontawesome}
\usepackage{hyphenat}
\usepackage{mathabx}
\usepackage[bbgreekl]{mathbbol}
\usepackage{relsize}
\usepackage[T1]{fontenc}
\usepackage{charter}
\usepackage{seqsplit}
\usepackage{breqn}
\usepackage{upgreek}
\usepackage{mathrsfs}
\usepackage[colorlinks=true, citecolor=DodgerBlue2, linkcolor=HotPink2, urlcolor=Aquamarine3]{hyperref}

\theoremstyle{plain}
\newtheorem{thm}{Theorem}[section]
\newtheorem{cor}[thm]{Corollary}
\newtheorem{pro}[thm]{Proposition}
\newtheorem{lem}[thm]{Lemma}
\theoremstyle{definition}
\newtheorem{dfn}[thm]{Definition}
\newtheorem{example}[thm]{Example}
\newtheorem{rmk}[thm]{Remark}
\newtheorem{conjecture}[thm]{Conjecture}
\newtheorem{question}[thm]{Question}

\DeclareSymbolFontAlphabet{\mathbb}{AMSb}
\DeclareSymbolFontAlphabet{\mathbbl}{bbold}

\let\onto\twoheadrightarrow
\def\into{\hookrightarrow}

\def\Hom{\mathop{\textup{Hom}}\nolimits}
\def\Isom{\mathop{\textup{Isom}}\nolimits}

\def\Perf{\mathop{\textbf{\textup{Perf}}}\nolimits}
\def\Perfd{\mathop{\textbf{\textup{Perfd}}}\nolimits}
\def\perf{\mathop{\textup{perf}}\nolimits}

\def\QSyn{\mathop{\textbf{\textup{QSyn}}}\nolimits}
\def\Spec{\mathop{\textbf{\textup{Spec}}}\nolimits}
\def\Spa{\mathop{\textbf{\textup{Spa}}}\nolimits}

\def\Spd{\mathop{\textbf{\textup{Spd}}}\nolimits}
\def\Sht{\mathop{\textbf{\textup{Sht}}}\nolimits}
\def\Bun{\mathop{\textbf{\textup{Bun}}}\nolimits}
\def\deg{\mathop{\textup{deg}}\nolimits}

\def\mod{\mathop{\textup{mod}}\nolimits}
\def\kernel{\mathop{\textup{Ker}}\nolimits}

\def\Lie{\mathop{\textup{Lie}}\nolimits}

\def\GL{\mathop{\textup{GL}}\nolimits}
\def\Gr{\mathop{\textup{Gr}}\nolimits}

\def\Ad{{\textup{Ad}}}

\def\id{{\textup{id}}}

\def\Def{{\textup{\emph{Def}}}}

\def\Set{\mathop{\textbf{Set}}\nolimits}

\newcommand{\prism}{{\mathlarger{\mathbbl{\Delta}}}}

\newcommand{\BF}{{\mathbb{F}}}
\newcommand{\BG}{{\mathbb{G}}}

\newcommand{\BI}{{\mathbb{I}}}

\newcommand{\BM}{{\mathbb{M}}}
\newcommand{\BN}{{\mathbb{N}}}

\newcommand{\BZ}{{\mathbb{Z}}}

\newcommand{\Fa}{{\mathfrak{a}}}

\newcommand{\Fh}{{\mathfrak{h}}}

\newcommand{\Fm}{{\mathfrak{m}}}

\newcommand{\FF}{{\mathfrak{F}}}

\newcommand{\FH}{{\mathfrak{H}}}

\newcommand{\CB}{{\mathcal{B}}}
\newcommand{\CC}{{\mathcal{C}}}
\newcommand{\CD}{{\mathcal{D}}}

\newcommand{\CF}{{\mathcal{F}}}
\newcommand{\CG}{{\mathcal{G}}}

\newcommand{\CN}{{\mathcal{N}}}
\newcommand{\CO}{{\mathcal{O}}}

\newcommand{\CT}{{\mathcal{T}}}

\newcommand{\CY}{{\mathcal{Y}}}

\newcommand{\SA}{{\mathscr A}}
\newcommand{\SB}{{\mathscr B}}

\newcommand{\SF}{{\mathscr F}}
\newcommand{\SG}{{\mathscr G}}
\newcommand{\SH}{{\mathscr H}}

\newcommand{\SL}{{\mathscr L}}
\newcommand{\SM}{{\mathscr M}}
\newcommand{\SN}{{\mathscr N}}

\newcommand{\SP}{{\mathscr P}}
\newcommand{\SQ}{{\mathscr Q}}

\newcommand{\SW}{{\mathscr W}}

\DeclareMathOperator*{\doublewedge}{\wedge\mkern-15mu\wedge}

\title{Deformations of Prismatic Higher $(G,\mu)$-Displays over Quasi-Syntomic Rings}

\author{S. Mohammad Hadi Hedayatzadeh}
\address{Institute for Research in Fundamental Sciences (IPM),
Niavaran Sq., Tehran, Iran}
\email{hadi@ipm.ir}

\author{Ali Partofard}
\address{Institute for Research in Fundamental Sciences (IPM),
Niavaran Sq., Tehran, Iran}
\email{alipartofard@ipm.ir}

\begin{document}

\begin{abstract}
    We prove a deformation theorem for prismatic higher $(G,\mu)$-displays over quasi-syntomic rings. As an application, we extend the classification of $p$-divisible groups via prismatic Dieudonn\'e modules to a class of rings, properly containing quasi-syntomic rings. Finally, we relate the stack of prismatic higher $(G,\mu)$-displays to integral local Shimura varieties.\vspace{7pt}
\end{abstract}

\maketitle

\begin{center}
    \noindent\rule{0.8\linewidth}{0.4pt}\vspace{5pt}


    \noindent\rule{0.8\linewidth}{0.4pt}
\end{center}
\vspace{0.5cm}

\tableofcontents


\setlength{\parskip}{1em}    

\setcounter{section}{0}
\phantomsection
\section*{Introduction}


The quest for linear algebraic invariants of motives with $G$-structure has a long history. The first examples are the classification of $p$-divisible groups through Dieudonn\'e modules over perfect fields and the classification of complex Abelian varieties via Hodge structures.

Zink introduced displays as a generalization of Dieudonn\'e modules when the base is not a perfect field. He, together with Lau, showed that, over a large class of rings, there is an equivalence of categories between the category of $p$-divisible groups and that of displays \cite{zink2002display,Lau:2008aa}.

Later, Langer and Zink, in the spirit of the Tannakian formalism, generalized the notion of displays to allow Hodge polygons with more than two slopes. The concept of \emph{higher displays} is closely related to the de Rham-Witt cohomology. Langer and Zink proved that under some conditions, the de Rham-Witt cohomology of a proper and smooth scheme of characteristic $p$ has a natural structure of higher display \cite{langer2007rham}.

In another direction, in order to avoid pathologies of the ring of Witt vectors of certain rings, Zink introduced the notion of \emph{windows} over \emph{frames} and showed that for ``nice'' rings, the category of windows over frames attached to these rings is equivalent to the category of displays over them \cite{Zink_2001,zbMATH05827284}.

When $G$ is a reductive (or affine smooth) group scheme and $\mu$ is a cocharacter of $G$, B\"ultel \cite{bultel2008pel}, B\"ultel-Hedayatzadeh \cite{bueltel2020g}, and Lau \cite{lau2021higher} introduced and studied $(G,\mu)$-displays and $(G,\mu)$-windows over arbitrary frames.

B\"ultel used the stack of $(G,\mu)$-displays to construct an integral model for certain Shimura varieties that are not of Abelian type \cite{bultel2008pel}. B\"ultel and Pappas used the stack of $(G,\mu)$-displays to construct an integral model for Rapoport-Zink spaces of Hodge type \cite{bultel2020displays}, and recently Bartling \cite{bartling2022mathcal} proved that this integral model agrees with the model constructed in \cite{scholze2020berkeley}.

Prismatic cohomology, developed by Bhatt and Scholze \cite{10.4007/annals.2022.196.3.5}, specializes to other known cohomology theories, including the de Rham-Witt cohomology. It is natural to seek a stack of linear algebraic objects such that the prismatic cohomology of ``good motives'' can be seen as objects in this stack. It is expected that this stack has a functor to the stack of displays, which is compatible with the specialization map from the prismatic cohomology of a good motive to its de Rham-Witt cohomology.

One step in this direction was the work of Ansch\"utz-Le Bras, which classifies $p$-divisible groups by a certain class of prismatic $F$-crystals (which they call \emph{admissible prismatic Dieudonn\'e crystals}) over quasi-syntomic rings \cite{anschutz2023prismatic}. In the introduction to their paper, they raised some questions and problems. Among them was extending the classification beyond the quasi-syntomic case and studying the deformation theory of these objects.

One goal of this paper is to partially answer these questions. We also introduce a new version of admissible crystals with $G$-structure called \emph{prismatic higher $(G,\mu)$-displays}, which are more adapted for quasi-syntomic rings, and we relate them to the points of integral local Shimura varieties over quasi-syntomic rings.

We prove a theorem in the spirit of Grothendieck-Messing for deformations of prismatic higher $(G,\mu)$-displays for a large class of nilpotent thickenings of quasi-syntomic rings. In fact, we allow for nilpotent thickenings by rings that are not necessarily quasi-syntomic rings. This enables us to obtain a classification of $p$-divisible groups in terms of prismatic $F$-crystals for a larger class of rings containing quasi-syntomic rings. We believe that the same methods should yield similar results in a more general setting.

We define a functor from the stack of prismatic higher $(G,
\mu)$-displays over a quasi-syntomic ring to the points of the integral local Shimura variety (attached to the datum $(G,\mu)$) over that ring. In the case where the integral local Shimura variety is representable and satisfies some conditions, we prove that this functor is an equivalence. The proof relies on some facts about the deformation theory of schemes over Henselian thickenings, and it might be possible to modify the proof to remove the condition of representability.

While we were working on this project, some independent works appeared that address similar questions for complete regular local rings \cite{imai2023prismatic, ito2023deformation} (and to some extent for smooth $\CO_K$-algebras, where $\CO_K$ is the ring of integers of some local field $K$). Although there are some similarities between our construction and these works (at least in their philosophy or techniques), it is worth noting that there exist several significant examples of quasi-syntomic rings that are not smooth $\CO_K$-algebras. Furthermore, some of our techniques could be useful in understanding the prismatic Dieudonn\'e crystals beyond the quasi-syntomic rings.

Let us now give a more detailed introduction to the content of the paper. Let $G$ be a reductive group over $\BZ_p$ and $\mu$ a minuscule cocharacter of $G$ (for the sake of simplicity of the introduction, we assume that $\mu$ is defined over $\BZ_p$). To this \emph{display datum} and a (bounded) orientable prism $(A,I)$ we attach the \emph{display group} $G^\mu(A,I)$, which is the subgroup of $G(A)$ consisting of those morphisms $\CO_G\to A$ that map the filtration on $\CO_G$ given by $\mu$ to the Nygaard filtration on $A$. The display group has a natural action on $G(A)$ via $\mu$-conjugation (see \S \ref{SectionBanalObjects}). We define the groupoid of \emph{banal prismatic higher $(G,\mu)$-displays} as the quotient groupoid $\CB^{\textup{ban}}_{(G,\mu)}(A,I):=\left[G(A):G^\mu(A,I)\right]$.

The assignment $(A,I)\mapsto \CB^{\textup{ban}}_{(G,\mu)}(A,I)$ enables us to define presheaves and, by sheafification, sheaves on various sites (see \S \ref{StackprismaticGdisplays}); on the relative prismatic site $(A,I)_{\prism}$, we obtain the stack $\mathcal{B}^{\prism_{/(A,I)}}_{(G,\mu)}$. Noting that \emph{quasi-regular semiperfectoid rings} form a basis for the quasi-syntomic site, and using the fact that the prismatic site of a quasi-regular semiperfectoid ring $R$ has a final object $(\prism_R,\BI_R)$, we obtain the stack $\mathcal{B}^{\textup{qsyn}}_{(G,\mu)}$. Finally, on $\Perf/\Spd \BZ_p$ (the site of perfectoid spaces of characteristic $p$ over $\Spd\BZ_p$, with the \'etale topology), and using $A_{\textup{inf}}$ prisms, we have the stack $\mathcal{B}^{\textup{perf}}_{(G,\mu)}$. One can define objects of these stacks, called \emph{prismatic higher $(G,\mu)$-displays}, as $G$-torsors (over the relevant spaces) with additional data. When $G$ is $\GL_n$, and over quasi-syntomic rings, the stack $\mathcal{B}^{\textup{qsyn}}_{(G,\mu)}$ is nothing but the stack of admissible prismatic Dieudonn\'e crystals, which is equivalent to the stack of $p$-divisible groups \cite{anschutz2023prismatic}.

One can attach a \emph{Hodge filtration} to a prismatic higher $(G,\mu)$-display. When $G$ is $\GL_n$, this is the same Hodge filtration attached to the $p$-divisible group corresponding to the admissible $F$-crystal.

We can now state one of the main theorems of this paper:

\begin{thm}[Theorem \ref{grothendique-messinggeneral}]
    Let $f: S'\onto S$ be a first-order thickening of quasi-syntomic rings. There is a natural bijection between the isomorphism classes of deformations of an adjoint-nilpotent\footnote{Please see Definition \ref{adjointnilpotentdefinition} for the definition of adjoint-nilpotence.} prismatic higher $(G,\mu)$-display on $S$ along $f$ and the isomorphism classes of deformations of the associated Hodge filtration from $S$ to $S'$.
\end{thm}

In fact, we prove a little more. We prove the theorem for \emph{admissible thickenings} (Theorem \ref{ThmDefAdmThink}). Let us explain. Let $R$ be an integral perfectoid $\widehat{\mathbb{Z}_p[\zeta_{p^\infty}]}$-algebra. We say that a quasi-regular semiperfectoid ring $S$ is 
\emph{finitely generated perfectly presented} if it can be written as $S=R/(f_1,\dots,f_r)$ where the sequence $(f_1,\dots,f_r)$ is $p$-completely regular relative to $\widehat{\mathbb{Z}_p[\zeta_{p^\infty}]}$ and each $f_i$ has a compatible system of $p$-power roots in $R$. The initial prism $\prism_S$ attached to such rings behaves nicely, and general statements (e.g., for the Nygaard filtration) about quasi-syntomic rings are reduced to these special rings. Now assume that $R$ is flat over $\mathbb{Z}_p[\zeta_{p^\infty}]$. We say that a thickening $\tilde{S}\onto S$ is \emph{admissible} if we can write $\tilde{S}$ as a quotient $R/I^2$ such that $S$ is $R/I$. Note that $\tilde{S}$ is quasi-syntomic only when $I$ is a principal ideal.

In the special case of $G=\GL_n$, this partially answers a question of Ansch\"utz and Le Bras in \cite{anschutz2023prismatic}. It seems plausible that this theorem extends to more general thickenings of rings, but we may have to modify the definition of our stack.

The general strategy of the proof is similar to that in previous work on the deformation theory of displays. However, we need new ingredients to make the strategy work. Let us now give an outline of the proof. We proceed by \emph{d\'evissage}. We first show that deformations of banal prismatic higher $(G,\mu)$-displays along nilpotent thickenings are banal. Then, by \'etale descent and the invariance of the \'etale site under nilpotent thickenings, we reduce the problem to the banal case. In the banal situation, we use Andr\'e's Lemma (\cite[Theorem 7.14]{10.4007/annals.2022.196.3.5}) to find a suitable cover by finitely generated perfectly presented algebras. Then, we reduce to the case of a ``universal'' first-order thickening
    \begin{multline*}
        S'(r):=\widehat{\mathbb{Z}_p[\zeta_{p^\infty}]}\langle X_1^{1/p^\infty},\cdots,X_n^{1/p^\infty}\rangle/(X_1^2,\cdots,X_r^2)\onto\\
        S(r):=\widehat{\mathbb{Z}_p[\zeta_{p^\infty}]}\langle X_1^{1/p^\infty},\cdots,X_n^{1/p^\infty}\rangle/(X_1,\cdots,X_r)
    \end{multline*} 
    
Using K\"unneth Formula (Proposition \ref{KunnethFormulae}), we can further reduce to the case of the thickening $S'(1)\onto S(1)$. For this thickening, we translate the deformation problem to that of deformations of windows on the Nygaard frames $\SA:=\SA_{\textup{Nyg}}(\prism_{S(1)},\BI_{S(1)})$ and $\SA':=\SA_{\textup{Nyg}}(\prism_{S'(1)},\BI_{S'(1)})$ naturally attached to $S(1)$ and $S'(1)$ (see Example \ref{ExNygFramePrism}).  The general practice in such scenarios is to modify the frame structure on the ring of $\SA'$, so that the groupoid of banal windows over the modified frame becomes equivalent to that over $\SA$. One can then reduce to comparing windows (and their deformations) over frames sharing the same underlying ring. A useful tool in the business of proving equivalences of categories of windows over different frames is the Unique Lifting Lemma (see Proposition \ref{uniqueliftinglemma} for the version we use). However, we cannot apply it from the start because the map $\prism_{S'(1)}\to \prism_{S(1)}$ is not surjective (unlike the case of Witt frames), and so the standard method does not work. Instead, we introduce some auxiliary frames and use some facts that we prove about the Nygaard filtration to apply the Unique Lifting Lemma to the truncation of those frames. We then lift the equivalence of windows to get a compatible equivalence of windows over all truncations of the frames. Finally, we use a limit argument to prove the desired equivalence of windows over the original frames.

As an application of Theorem \ref{grothendique-messing}, we give a classification of $p$-divisible groups over a class of rings (containing quasi-syntomic rings) in terms of prismatic higher displays over a prism canonically attached to these rings.

\begin{thm}[Theorem \ref{ClasspDivGrps}]
    Let $R$ be an integral perfectoid flat $\mathbb{Z}_p[\zeta_{p^\infty}]$-algebra and $\tilde{S}$ be an $R$-algebra. Assume that there exists an admissible thickening $\tilde{S}\onto S$ for some quasi-syntomic $R$-algebra $S$. There is an equivalence of categories between the category of formal $p$-divisible groups over $\tilde{S}$ and the category of adjoint-nilpotent prismatic higher displays over $\prism_{\tilde{S}}$.
\end{thm}

Over quasi-syntomic rings, we expect that prismatic higher $(G,\mu)$-displays are related to the points of local Shimura varieties. Assume that $(G,\mu,b)$ is an integral Shimura datum. One can associate a prismatic higher $(G,\mu)$-display $\CD_b$ to this datum (see \cite{bueltel2020g}). Over the site of $p$-torsion-free quasi-regular semiperfectoid rings (with quasi-syntomic topology), we define a stack $\mathcal{B}^{\textup{qsyn}}_{(G,\mu,b)}$, by sending $R$ to the pairs $(\CD,\eta)$ where $\CD\in \mathcal{B}^{\textup{qsyn}}_{(G,\mu)}(R)$ and $\eta$ is a quasi-isogeny $\eta:\bar{\CD}\dashrightarrow \bar{\CD}_b$ between the reductions mod $p$. We show that there is a natural functor $\mathfrak{Fib}_G:\mathcal{B}^{\textup{qsyn}}_{(G,\mu,b)}\to{\mathcal{M}^{\textup{int}}(G,\mu,b)}$, and by a group-theoretic adaptation of the arguments in \cite{anschutz2023prismatic}, we prove the following.

\begin{thm}[Theorem \ref{equivalencetolocalshimura}]
     If the integral local Shimura variety $\mathcal{M}^{\textup{int}}(G,\mu,b)$ is representable by the diamond of a formal scheme, then the natural transformation $\mathfrak{Fib}_G$ is an equivalence.
\end{thm}

We need the representability assumption because our proof uses the deformation theory of integral local Shimura varieties. Of course, this representability is expected to hold in general; however, we believe that one should be able to prove the statement without this assumption.

\begin{conjecture}
    The natural transformation $\mathfrak{Fib}_G:\mathcal{B}^{\textup{qsyn}}_{(G,\mu,b)}\to{\mathcal{M}^{\textup{int}}(G,\mu,b)}$ is an isomorphism.
\end{conjecture}

When $R$ is a general non-quasi-syntomic ring (outside the class considered in this paper), one expects that the appropriate framework is that of derived prisms, in the sense of Bhatt–Lurie’s prismatic F-gauges \cite{bhatt2022prismatic}, together with the geometrization via derived prismatization/syntomification. We expect to address this perspective in future work.

Following the completion of the present paper, Gardner and Madapusi \cite{gardner2024algebraicity} established a far-reaching generalization of Theorem \ref{grothendique-messing} by means of derived prismatization and derived algebraic geometry. In particular, their work yields a completely general classification of $p$-divisible groups over $p$-complete rings in terms of vector bundles on the syntomification, extending earlier results in the quasi-syntomic setting. In the non–quasi-syntomic case, the associated syntomifications are genuinely derived stacks. In contrast, the results obtained here are established by more concrete arguments and yield deformation results for prismatic Dieudonné modules and, more generally, for higher $(G,\mu)$-displays. In the range of rings considered in this paper, the resulting structures remain within a classical (underived) geometric framework. These deformation results, which are of independent interest (for instance, in the study of local Shimura varieties), may thus be viewed as complementary to the derived-geometric classification.

\subsection*{Structure of the Paper}

In Section \ref{review}, we review the theory of prisms and quasi-syntomic rings. In Section \ref{SecNygFil}, we prove some facts we need about the structure of the Nygaard filtration of a quasi-syntomic ring and the thickenings between such rings. This section is used only in Section \ref{deformationtheory}. In Section \ref{Secprisdisp}, we define the notion of prismatic higher displays over a prism and the stack of prismatic higher displays. This notion generalizes the notion of admissible prismatic Dieudonn\'e modules of \cite{anschutz2023prismatic}. In Subsection \ref{SectionBanalObjects}, we define the display group and the groupoid of banal prismatic higher $(G,\mu)$-displays. In Subsection \ref{StackprismaticGdisplays}, we define the stack of prismatic higher $(G,\mu)$-displays over various sites, and prove some descent theorems about them. In Section \ref{SecWindows}, we review the definition and basic properties of windows and frames and use them to compare prismatic higher $(G,\mu)$-displays with Witt displays. We also introduce the notion of \emph{adjoint-nilpotence} in this section. Section \ref{deformationtheory} is the heart of the first part of the paper. We prove our main theorem about the deformation theory of adjoint-nilpotent prismatic higher $(G,\mu)$-displays over finitely generated perfectly presented rings. Then, we use descent to prove it for thickenings of more general quasi-syntomic rings. Next, we define a canonical prism for rings that admit an admissible thickening by a quasi-syntomic ring and use this prism to define prismatic higher $(G,\mu)$-displays for such rings. At the end of this section, we generalize the classification of formal $p$-divisible groups for such rings. In Section \ref{removingcondition}, we use simplicial methods to remove condition $\leftmoon$. Section \ref{Secrelationshtukas} states the relation between prismatic higher $(G,\mu)$-displays and shtukas over perfectoid rings. The result of this section appeared partially in \cite{bartling2022mathcal}. In Section \ref{Secrelationlocalshimuravariety}, we relate the prismatic higher $(G,\mu)$-displays to the points of integral local Shimura varieties. In Appendix \ref{proofofstructure} we prove the structure theorem of the display group. Appendix \ref{shtukas} is a review of the theory of shtukas and integral local Shimura varieties. In Appendix \ref{SecHenselianpair} we review some results about Henselian pairs that we need in the paper.  

\subsection*{Notations and Conventions}\label{Not&Conv}

Here is a list of notations and conventions that we will use throughout the paper. Some of these notations will make sense later, but we have collected them here for the sake of completeness and convenience of the reader.

\small

\begin{itemize}
    \renewcommand\labelitemi{--}
   
    \item $p$ is a prime number.
    \item $q$ is a power of $p$. It will also be a variable (see Section \ref{SecNygFil}). We believe that the risk of confusion is very small. 
    \item Unless specified to the contrary, all rings are commutative with $1$.
    \item We denote the group of units of a ring $R$ by $R^\times$.
    \item The Jacobson radical of a ring $R$ is denoted by $\textup{Rad}(R)$.
    \item The Frobenius morphism $x\mapsto x^p$ of a ring of characteristic $p$ is denoted by $\textup{Frob}_p$.
    \item All prisms are assumed to be bounded.
    \item By a \emph{perfectoid} ring, we mean an \emph{integral perfectoid} ring (cf. Definition \ref{DefIntPerfRing}).
    \item For a ring $R$, we denote by $W(R)$ the ring of ($p$-typical) Witt vectors. We denote by $[\cdot]:R\into W(R)$ the Teichmüller lift.
    \item $\BZ_q:=W(\BF_q)$.
    \item The $p$-adic completion of a ring $R$ is denoted by $R^\wedge$.
    \item In a $\BZ_p[q]$-algebra (here $q$ is a variable), and for a natural number $n$, we set $[n]_q:=q^{n-1}+q^{n-2}+\dots+q+1$.
    \item The $(p,[p]_q)$-adic completion of a $\BZ_p[q]$-algebra $R$ is denoted by $R^\doublewedge$ (cf. Section \ref{SecNygFil}).
    \item $R_0$ is the integral perfectoid ring $\mathbb{Z}_p[\zeta_{p^\infty}]^\wedge$.
    \item $A_0$ is the $(p,q-1)$-adic completion of $\mathbb{Z}_p[q^{1/p^\infty}]$.
    \item There exist elements $w\in R_0^\flat$ and $u\in A_0^\times$ (which we fix once and for all), such that $[p]_q=[w^p]+\sigma(u)p$. Here $\sigma$ is the Frobenius lift, sending $q$ to $q^p$.
    \item For $n\in\BN$, we set $R_n:=R_0\langle X_1^{1/p^\infty},\cdots,X_n^{1/p^\infty}\rangle$.
    \item For $n\in\BN$, we set $S(n):=R_n/(X_1,\cdots,X_n)$.
    \item For $n\in\BN$, we set $S'(n):=R_n/(X_1^2,\cdots,X_n^2)$.
    \item For $n\in\BN$, we set $A_n:=A_{\textup{inf}}(R_n)$.
    \item Let $R$ be a ring and $I\subset R$ a Cartier divisor. For $n\in\BN$, we set $I^{\otimes -n}:=\Hom_R(I^{\otimes n},R)$. We have inclusions $I^{\otimes -n}\subset I^{\otimes -(n+1)}$, and we set $R[\frac{1}{I}]:=\varinjlim_{n\in\BN}I^{\otimes -n}$. If $I$ is generated by an element $d\in R$, then this is canonically isomorphic to the localization $R[\frac{1}{d}]$.  
\end{itemize}

\normalsize

\section{Review of Prisms and Perfectoid Spaces}\label{review}

For ease of reference, we review some definitions and theorems that we need about prisms and perfectoid spaces.

We start with a review of prisms. The main reference is \cite{10.4007/annals.2022.196.3.5}. Let $p$ be a prime number.

\begin{dfn}
    A \emph{$\delta$-ring} is a pair $(A,\delta_A)$, where $A$ is a commutative ring and $\delta_A:A \to A$ is a map of sets that satisfies the following identities for all elements $a,b\in A$:

        \begin{enumerate}
            \item $\delta_A(0)=0=\delta_A(1)$  
            \item $\delta_A(ab)=a^p\delta_A(b)+b^p\delta_A(a)+p\delta_A(a)\delta_A(b)$
            \item $\delta_A(a+b)=\delta_A(a)+\delta_A(b)-\sum_{i=1}^{p-1}\frac{1}{p}\binom{p}{i}a^ib^{p-i}$
        \end{enumerate}

    We define a Frobenious lift on $A$ by $\sigma_A(a)=a^p+p\delta_A(a)$. Morphisms of $\delta$-rings are morphisms of the underlying rings that commute with $\delta$. By abuse of terminology and when confusion is not likely, we say that $A$ is a $\delta$-ring. We may also drop the subscript $A$ from the Frobenius lift $\sigma_A$. If $A\to B$ is a morphism of $\delta$-rings, we say that $B$ is a $\delta$-ring over $A$.
\end{dfn}

\begin{pro}{\cite[Lemmas 2.17,2.18]{10.4007/annals.2022.196.3.5}}\label{stabilityunderetalemap}
    Let $(A,\delta_A)$ be a $\delta$-ring. Assume that $f:A\to B$ is an \'etale morphism. Then, there is a canonical structure of a $\delta$-ring on $B$ such that $f$ is a morphism of $\delta$-rings. If an ideal $J\subset A$ is stable under $\delta_A$, then the $J$-adic completion of $A$ has a canonical structure of a $\delta$-ring such that the completion $A\to \hat{A}$ is a morphism of $\delta$-rings.
\end{pro}

\begin{pro}{\cite[Théorème 2]{zbMATH03955079}}\label{DeltaWittAdj}
    The forgetful functor from the category of $\delta$-rings to the category of rings has a right adjoint, namely the functor of ($p$-typical) Witt vectors. In particular, for any $\delta$-ring $A$, there is a natural morphism of $\delta$-rings $A\to W(A)$ lifting the projection $w_0:W(A)\to A$, called the \emph{Cartier (diagonal) map}.
\end{pro}

\begin{dfn}
    A \emph{prism} is a pair $(A,I)$ where $A$ is a $\delta$-ring and $I$ is a Cartier divisor in $A$, such that $p\in I+\sigma_A(I)$ and $A$ is derived $(p,I)$-adically complete. A map of prisms $(A,I)\to (B,J)$ is a map of $\delta$-rings that sends $I$ to $J$. Let $(A,I)$ be a prism. It is called \emph{orientable} if $I$ is a principal ideal. If a generator of $I$ is chosen, the prism is called \emph{oriented}. We say that $(A,I)$ is \emph{bounded} if the ring $A/I$ has bounded $p^\infty$-torsion, i.e., there is an integer $n\geq1$, such that $(A/I)[p^\infty]=(A/I)[p^n]$. If $A$ is a $\BZ_q$-algebra, we say that $(A,I)$ is a $\BZ_q$-prism. Finally, we say that $(A,I)$ is \emph{perfect} if the Frobenius lift $\sigma_A$ is a bijection. Let $R$ be a ring. A prism over $R$ is a prism $(A,I)$ endowed with a morphism $R\to A/I$.
\end{dfn}

\begin{rmk}
    From now on, all prisms are assumed to be bounded. In that case, $A$ is $(p,I)$-adically complete.
\end{rmk}

\begin{pro}[Rigidity]{\cite[Lemma 3.5]{10.4007/annals.2022.196.3.5}}\label{rigidity}
    If $(A,I)\to (B,J)$ is a map of prisms, then we have $J=IB$. 
\end{pro}

\begin{dfn}\label{DefAbsPrismSite}
    Let $R$ be a ring. The \emph{(absolute) prismatic site} over $R$ is the site whose underlying category is the opposite category of prisms over $R$ and endowed with the topology for which covers are given by $(p,I)$-completely \'etale maps. We denote this site by $(R)_\prism$. If $(A,I)$ is a prism over $R$, we denote the slice of the prismatic site over $(A,I)$ by $(R)_{\prism/(A,I)}$, or by abuse of notations, $(R)_{\prism/A}$. If this category has a final object, we denote it by $(\prism_{R/A}, \BI_{R/A})$.
\end{dfn}

\begin{dfn}
    Let $(A,I)$ be a prism. The \emph{Nygaard filtration} on $A$ is the decreasing filtration on $A$, whose $i$-th term is given by $\CN_A^i:=\sigma_A^{-1}(I^i)$ (for $i\geq 0$).
\end{dfn}

\begin{pro}{\cite[Lemma 4.28]{anschutz2023prismatic}}\label{Henselianpair}
    Let $(A,I)$ be a prism. Then the pair $(A, \CN^1_A)$ is Henselian.
\end{pro}

Now we want to define the prismatic envelope of a regular sequence and its relation to the pd-envelope.

\begin{pro}{\cite[Corollary 2.44]{10.4007/annals.2022.196.3.5}}
    Let $A$ be a $\delta$-ring. Let $x$ be an element of $A$ and $I$ a Cartier divisor in $A$. Then the slice category of $\delta$-rings $B$ over $A$ such that $x\in IB$ has an initial object, denoted by $A\{\frac{x}{I}\}$. If $I$ is generated by a (non-zero divisor) $y\in A$, we denote this initial object by $A\{\frac{x}{y}\}$.
\end{pro}

\begin{pro}{\cite[Corollary 2.39]{10.4007/annals.2022.196.3.5}}\label{pdenvelope}
    Let $A$ be a $p$-torsion-free $\delta$-ring, and $x_1, x_2,\cdots, x_r$ be a regular sequence in $A/p$. Then the $\delta$-ring $A\left\{\frac{x_1}{p},\frac{x_2}{p},\cdots,\frac{x_r}{p}\right\}$ is the pd-envelope of the ideal $(x_1,x_2,\cdots,x_r)$ in $A$.
\end{pro}

\begin{pro}{\cite[Lemma V.5.1]{bhatt2018lectures}}
    Let $(A,I)$ be a prism and $x_1,x_2,\cdots,x_r$ be a regular sequence in $A/I$. Then $(C,I)=\left(A\left\{\frac{x_1}{I},\frac{x_2}{I},\cdots,\frac{x_r}{I}\},IA\{\frac{x_1}{I},\frac{x_2}{I}\right\}\right)$ is the initial object in the category whose objects are prisms $(B,IB)$ together with a map $\left(A,I+\sum_{i=1}^rx_iA\right)\to (B,IB)$. This object is called the \emph{prismatic envelope} of $\left(A,(I,x_1,\cdots,x_r)\right)$.
\end{pro}

\begin{dfn}{\cite{bhatt2018integral}}\label{DefIntPerfRing}
    A ring $R$ is \emph{integral perfectoid} if it satisfies the following conditions:
        \begin{enumerate}
            \item it is $\pi$-complete for some $\pi\in R$ with $p\in \pi^pR$ 
            \item $R/p$ is semiperfect
            \item the kernel of the Fontaine map $\theta_R:A_{\textup{inf}}(R):=W(R^\flat)\to R$ is a principal ideal.
        \end{enumerate}
\end{dfn}

We also need some facts about quasi-syntomic rings.

\begin{dfn}\label{quasi-syntomic}
    A ring $A$ is \emph{quasi-syntomic} if it is $p$-complete, has bounded $p^\infty$-torsion, and the cotangent complex $L_{A/\mathbb{Z}_p}$ has $p$-complete Tor amplitude in the interval $[-1,0]$. A homomorphism $A\to B$ of $p$-complete rings with bounded $p^\infty$-torsion is a \emph{quasi-syntomic cover} if it is $p$-completely faithfully flat and $L_{B/A}$ has $p$-complete Tor amplitude in $[-1,0]$. The category of quasi-syntomic rings with quasi-syntomic covers is a site (cf. \cite[Lemma 4.17]{bhatt2019topological}) denoted by $\QSyn$.
\end{dfn}

\begin{dfn}
    Let $S=\Spa R$ be a formal $p$-complete scheme. We denote the site whose underlying category is the opposite category of $p$-complete bounded $p^\infty$-torsion $R$-algebras with quasi-syntomic covers by $(\tilde{S})_{\textup{QSYN}}$.
\end{dfn}

The site of quasi-syntomic rings has a basis with good properties:

\begin{dfn}
    A ring $R$ is \emph{quasi-regular semiperfectoid} (or \emph{qrsp} for short) if $R$ is quasi-syntomic and there is a surjection $S\onto R$ from a perfectoid ring $S$.
\end{dfn}

\begin{example}
    If $S=R/(f_1,f_2,\cdots,f_n)$ where $R$ is a perfectoid ring and $f_1,f_2,\cdots,f_n$ is a regular sequence, then $S$ is quasi-regular semiperfectoid.
\end{example}

\begin{dfn}
    Let \emph{$\textbf{QRSPerfd}$} be the site of qrsp rings with quasi-syntomic covers (this is actually a site by \cite[Lemma 4.27]{bhatt2019topological}).
\end{dfn}

\begin{pro}{\cite[Lemma 4.28]{bhatt2019topological}}\label{propertiesofqrsp}
    Every quasi-syntomic ring admits a quasi-syntomic cover by a qrsp ring. 
\end{pro}

\begin{pro}{\cite[Proposition 4.31]{bhatt2019topological}}\label{equivalenceofsites}
    For any presentable ($\infty-$)category $\CC$ there is an equivalence of categories between the category of sheaves on the site $\QSyn$ with values in $\CC$ and the category of sheaves on the site $\textbf{QRSPerfd}$ with values in $\CC$.
\end{pro}

\begin{pro}{\cite[Proposition 7.2]{10.4007/annals.2022.196.3.5}}\label{InitObjPrisSite}
    The absolute prismatic site of a qrsp ring has a final object denoted by $(\prism_R, \BI_R)$. The formation of $\prism_R$ is functorial and the ideal $\BI_R$ is principal (because by functionality there is a map from a perfect prism to $\prism_R$.)
\end{pro}

We also need to work over the slice categories $\QSyn/\mathbb{Z}_q$ and $\textbf{QRSPerfd}/\mathbb{Z}_q$. To do this, we need the following lemma.

\begin{lem}\label{initialinZqprisms}
    Let $R$ be a qrsp ring. Every morphism $\mathbb{Z}_q\to R$ has a canonical lift $\mathbb{Z}_q\to \prism_R$.
\end{lem}

\begin{proof}
    Pick $f:\BZ_q\to R$. First, assume that $R$ is perfectoid. The reduction $\bar{f}:\BF_q\to R/p$ lifts canonically to a map $\BF_q\to R^\flat=\varprojlim R/p$, and therefore we have a canonical lift $\BZ_q=W(\BF_q)\to W(R^\flat)=\prism_R$. Now, assume that $R$ is an arbitrary qrsp ring and choose a surjection $S\onto R$ from a perfectoid ring. By \cite[Corollary 2.10]{anschutz2023prismatic}, we can assume that $S$ is Henselian with respect to the kernel of $S\onto R$, and so $f$ lifts up to a map $\BZ_q\to S$. Now, by the first part, there is a canonical lift to $\prism_S$. By functoriality of $\prism$ (Proposition \ref{InitObjPrisSite}), we have a map $\prism_S\to\prism_R$, and so there is a lift of $f$ to $\prism_R$. The fact that the lift is canonical can be proved similarly to the proof of \cite[Lemma 4.8]{10.4007/annals.2022.196.3.5}, noting that since the extension $\BZ_q/\BZ_p$ is \'etale, the cotangent complex $L_{\prism_S/\BZ_q}$ vanishes after derived $p$-completion.
\end{proof}
\begin{pro}\label{qrsptoflat}\cite[Proposition 7.5]{bhatt2022prismatization}
Let $R\to R'$ be a morphism between qrsp rings, then $\prism_R\to \prism_{R'}$ is flat.
\end{pro}
\begin{pro}[K\"unneth Formula]\label{KunnethFormulae}
    Let $R, R'$ and $R''$ be qrsp rings and assume that $R'$ and $R''$ are $R$-algebras. Then we have $$\CN_{\prism_{R'}\otimes_{\prism_{R}}\prism_{R''}}^i=\sum_{r+s=i} \CN_{\prism_{R'}}^r\otimes_{\prism_{R}} \CN_{\prism_{R''}}^s.$$
\end{pro}

\begin{proof}
    The prismatic sites of $R'$ and $R''$ have final objects $\prism_{R'}$ and $\prism_{R'}$ and we have isomorphisms $\prism_{R'/\prism_R}\cong\prism_{R'}$ and $\prism_{R''/\prism_R}\cong\prism_{R''}$ (see Definition \ref{DefAbsPrismSite} for the notations). By \cite[Proposition 3.30]{anschutz2023prismatic} we have an isomorphism $\prism_{R'\otimes_{R} R'}\cong\prism_{R''}\otimes_{\prism_R} \prism_{R''}$. By definition, $\sum_{r+s=i} \CN_{\prism_{R'}}^r\otimes_{\prism_{R}} \CN_{\prism_{R''}}^s$ is a subset of $\CN_{\prism_{R'}\otimes_{\prism_{R}}\prism_{R''}}^i$. To check that we have an equality, note that under the Hodge-Tate comparison, the image of Nygaard filtration is the conjugate filtration \cite[Theorem 12.2]{10.4007/annals.2022.196.3.5}, and we have K\"unneth Formula for the cotangent complex. Therefore the inclusion should be an equality; for the details, look at the proof of \cite[Proposition 12.8]{bhatt2021prismatic}. 
\end{proof}

\begin{pro}\label{etalecoverofqrsp}
    An \'etale cover of a semiperfect (resp. qrsp) ring is semiperfect (resp. qrsp).
\end{pro}

\begin{proof}
    Let $R\to R'$ be an \'etale cover. First, assume that $R$ is semiperfect. Then by \cite[\href{https://stacks.math.columbia.edu/tag/0F6W}{Lemma  0F6W}]{stacks-project} the relative Frobenius $\sigma_{R'/R}:R'\otimes_{\sigma_R} R\to R' $ is an isomorphism. Since $\sigma_{R'/R}$ factors through $\sigma_{R'}$, $\sigma_{R'}$ is surjective, and so $R'$ is also semiperfect.

    Now, assume that $R$ is qrsp, with a map from a perfectoid ring $S$. Then $L_{R'/\mathbb{Z}}=L_{R/\mathbb{Z}}$. By the first part of the proposition, $R'/p$ is semiperfect, and there is the map $S\to R\to R'$ from the perfectoid ring $S$ to $R'$. This implies that $R'$ is qrsp by \cite[Remark 4.22]{bhatt2019topological}.
\end{proof}

\section{Nygaard Filtration of Quasi-Regular Semiperfectoid Rings}\label{SecNygFil}

In this section, we collect the algebraic facts that we need about the final object of the prismatic site of a quasi-regular semiperfectoid ring. Throughout this section, we denote the $p$-adic completion of a ring $A$ by $A^\wedge$. Let $R_0$ be the integral perfectoid ring $\mathbb{Z}_p[\zeta_{p^\infty}]^\wedge$.

Let $(A_0,I_0)$ be the perfect prism associated with $R_0$, that is, $A_0$ is the $(p,q-1)$-adic completion of $\mathbb{Z}_p[q^{1/p^\infty}]$ and $I_0$ is generated by $[p]_q$. Here $q$ is a formal variable and $[p]_q$ is the ``$q$-analogue'' of $p$, where: 
$$[n]_q=q^{n-1}+q^{n-2}+\cdots+q+1$$ and Frobenius is defined by $\sigma(q)=q^p$. Note that the Frobenius of $A_0$ induces an isomorphism $A_0/\sigma^{-1}(I_0)\cong R_0$, and this is the $A_0$-algebra structure, with which we will endow $R_0$. We denote the $(p,[p]_q)$-adic completion of a ring $A$ by $A^\doublewedge$. 

We also define the ``$q$-analogue" of factorial as
$$[n]_q!=\prod_{i\le n}[i]_q$$

The element $[p]_{q^{1/p}}:=\sigma^{-1}([p]_q)$ is a primitive element of $A_0$, so we can write it as $[p]_{q^{1/p}}=[w]+up$ for some $w\in R_0^\flat$ and $u\in A_0^\times$, which we fix. So, we have $$[p]_q=\sigma([p]_{q^{1/p}})=[w^p]+\sigma(u)p.$$

Let us define two rings that will serve as prototypes of the rings that we will consider throughout the paper. For all $n\in\BN$, we set 
\[S(n):=R_0\langle X_1^{1/p^\infty},\cdots,X_n^{1/p^\infty}\rangle/(X_1,\cdots,X_n)\] and 
\[S'(n):=R_0\langle X_1^{1/p^\infty},\cdots,X_n^{1/p^\infty}\rangle/(X_1^2,\cdots,X_n^2)\]

First, we want to consider a key example. Consider the qrsp ring $S(1)$. The final object of the prismatic site of $S(1)$ is the ring
    \begin{equation}\label{FinObjPrismSite}
        \prism_{S(1)}=\left(A_0\langle Y^{1/p^\infty}\rangle\{\frac{Y^p}{[p]_q}\}\right)^{\doublewedge}
    \end{equation}

The structure morphism $\prism_{S(1)}\to {S(1)}$ sends $Y$ to $X$ and on $A_0$ is Frobenius. More explicitly, $\prism_{S(1)}$ is the $(p,[p]_q)$-adic completion of the ring
$$\bigoplus_{i\in\mathbb{N}[\frac{1}{p}]}A_0\frac{Y^i}{[\lfloor i\rfloor]_q!}$$

\begin{lem}\label{Nygaardstructurespecial}
    The image of $\CN^1_{\prism_{S(1)}}$  in $\prism_{S(1)}/[w]$ is the pd-ideal generated by $p$ and the image of $Y$.
\end{lem}

\begin{proof}
    By \cite{10.4007/annals.2022.196.3.5}, the Nygaard filtration of $\prism_{S(1)}$ is given by $$\CN^1_{\prism_{S(1)}}=\bigoplus_{i\in\mathbb{N}[\frac{1}{p}],0\le i<1} [p]_{q^{1/p}}A_0 Y^i\oplus A_0\,Y \oplus\bigoplus_{i\in\mathbb{N}[\frac{1}{p}],i>1} A_0\frac{Y^i}{[\lfloor i\rfloor]_q!}$$ 
    
    Modulo $[w]$ we have $[\lfloor i\rfloor]_q!\equiv i!$ and $[p]_{q^{1/p}}\equiv p$. The result is now clear. 
\end{proof}

Now consider ${S'(1)}=R_0\langle X^{1/p^\infty}\rangle/{X^2}$. Like before, the final object of the prismatic site of ${S'(1)}$ is:
$$\prism_{S'(1)}=\left(A_0\langle Y^{1/p^\infty}\rangle\{\frac{Y^{2p}}{[p]_q}\}\right)^{\doublewedge}$$ or more explicitly, the $(p,[p]_q)$-adic completion of the ring:
$$\bigoplus_{i\in\mathbb{N}[\frac{1}{p}]}A_0.\frac{Y^{2i}}{[\lfloor i\rfloor]_q!}$$

\begin{lem}\label{specialcaseuniversalthickening}
    There is a surjective map $h:\prism_{S(1)}\to S'(1)$ whose precomposition with $\prism_{S'(1)}\to \prism_{S(1)}$ is the canonical map $\prism_{S'(1)}\to S'(1)$. Furthermore, the image of $\kernel h$ in $\prism_{S(1)}/[w]$ is a pd-ideal.
\end{lem}

\begin{proof}
    Consider the ideal $J\subset \prism_S$ generated by $Y^2$ and $\frac{Y^i}{[\lfloor i\rfloor]_q!}$ for $ i\in \BN[\frac{1}{p}], i>2$. It contains the image of $\CN^1_{\prism_{S'(1)}}$ and from the explicit description of $\prism_{S(1)}$ it is clear that $\prism_{S(1)}/J=S'(1)$. Similarly to lemma \ref{Nygaardstructurespecial}, one can check that the image of $J$ is a pd-ideal in $\prism_{S(1)}/[w]$.
\end{proof}

\begin{rmk}\label{JisHensel}
    Note that by construction $J$ is contained in $\CN^1_{\prism_{S(1)}}$, and so the pair $(\prism_{S(1)},J)$ is also Henselian (Proposition \ref{Henselianpair}).
\end{rmk}

Now we want to prove the above lemmas for a larger class of rings that form a basis for the topology of $\QSyn$.

\begin{dfn}\label{perfectlypresented}
    Let $R$ be an integral perfectoid $R_0$-algebra. A qrsp $R$-algebra $S$ is \emph{finitely generated perfectly presented} (fgpp for short) if it can be written as 
        \begin{equation}\label{fgpp}
             S=R/(f_1,\cdots,f_r)
        \end{equation}
    where $f_1,\cdots,f_r$ is a sequence of elements that is $p$-completely regular relative to $R_0$ and each $f_i$ has a compatible system of $p$-power roots in $R$.
\end{dfn}

\begin{example}
    Set $R_r:=R_0\langle X_1^{1/p^\infty},\cdots,X_r^{1/p^\infty}\rangle$. The typical example of an fgpp $R_r$-algebra is $S(r)$. In fact, for a general $R_0$-algebra $R$ and fgpp $R$-algebra $S$ written as $S=R/(f_1,\cdots,f_r)$, there is a canonical isomorphism \begin{equation}\label{EqGenFGPPIso}
         S\cong S(r)\hat{\otimes}_{R_r}R   
    \end{equation}
     and we will see that quite often questions about general cases of fgpp rings can be reduced to this special case. 
\end{example}

\begin{lem}\label{Nygaardofperfectlypresented}
    Let $S$ be an fgpp $R$-algebra (with a presentation as in (\ref{fgpp})). The ideal $\widebar{\CN}^1_{\prism_S}$ (the image of $\CN_S^1$ in $\prism_S/[w]$) is a pd-ideal in $\prism_S/[w]$, generated by elements $f_1^{1/p},\cdots,$ $f_r^{1/p}$ and $p$.
\end{lem}

\begin{proof}
    First, consider the case of $S(r)$. Here, the ring $R$ is $R_r=R_0\langle X_1^{1/p^\infty},\cdots,X_r^{1/p^\infty}\rangle$ and $f_i$ are $X_i$. We have a natural isomorphism \[S(r)\cong\widehat{\bigotimes_{i=1}^r}S(1)\]

    Now, using Lemma \ref{Nygaardstructurespecial} and the K\"{u}nneth formula (see Proposition \ref{KunnethFormulae}), the statement follows. For general $S$, consider the canonical isomorphism (\ref{EqGenFGPPIso}) \[S\cong S(r)\hat{\otimes}_{R_r}R\]
    
    Then we have $$\prism_S=\prism_{S(r)}\overset{\doublewedge}{\otimes}_{A_r}A_{\textup{inf}}(R)$$ and $\CN_S^1=\CN^1_{S(r)}\otimes_{A_r} A_{\textup{inf}}(R)$
    where $A_r:=A_{\textup{inf}}(R_r)$. Therefore, the image of $\CN_S^1$ in $\prism_S/[w]$ is generated by the image of $\CN_{S(r)}^1$ and we are reduced to the previous case (note that by Proposition \ref{qrsptoflat} the map $\prism_{S(r)}\to \prism_S$ is flat and therefore the image of a pd-ideal is a pd-ideal).
\end{proof}

\begin{lem}\label{universalthickening}
    If $S'\onto S$ is a first-order thickening of fgpp $R$-algebras, then there is a surjective map $h:\prism_S\to S'$ such that the composition $\prism_{S'}\to \prism_S\to S'$ is the canonical map. Furthermore, the image of $\kernel h$ in $\prism_{S(1)}/[w]$ is a pd-ideal.
\end{lem}

\begin{proof}
    We have presentations $S=R/I$ and $S'=R/I'$ as in (\ref{fgpp}). As before, one has $S=S(n)\widehat{\otimes}_{R_n} R$ for some $n$. Since $I^2\subset I'$, there is a surjective map from $S''=S'(n)\widehat{\otimes}_{R_0} R$ to $S'$. Therefore, it is enough to prove that there is a surjective map from $h_r:\prism_{S(n)}\to S'(n)$. By K\"{u}nneth formula, one can reduce to the lemma \ref{specialcaseuniversalthickening}. To see that the image of $\kernel h$ is a pd-ideal, it is enough to prove that $\kernel h=(\kernel h_r)\prism_{S}+\CN^1_{S'}\prism_{S}$ because the images of $\kernel h_r$ and $\CN^1_{S'}$ in $\prism_S/[w]$ are pd-ideals. Note that by surjectivity of the composition $\prism_{S''}\to \prism_{S}\to S''$, one can write $x\in \prism_{S}$ as $y+z,y\in (\kernel h_r)\prism_{S}$ and $z\in S''$. We know that $(\kernel h_r)\prism_{S}\subset \kernel h$, therefore, if $x\in \kernel h$ then $z\in \kernel h
    \cap S'=\CN^1_{S'}$ as desired.   
\end{proof}

\begin{rmk}
    If $S$ is a semiperfect ring, the lemma is trivial because $\prism_S=A_{\textup{crys}}(S)$ is the universal pd-thickening of $S$. The fact that the lemma is true in general suggests that there should be a description of $\prism_S$ as a universal thickening.
\end{rmk}

For technical reasons, we first consider the quasi-syntomic rings $S$ that satisfy the following condition: 

\begin{center}
\begin{tcolorbox}[enhanced,
  colframe=white, 
  colback=white,  
  boxrule=0pt,    
  sharp corners,  
  width=\textwidth, 
  halign=center   
]
$(\leftmoon)$\quad There is a quasi-syntomic hypercover $\tilde{S}^\bullet$ of $S$
such that for each $i$, $\tilde{S}^i$ is a qrsp $R_i$-algebra for a
perfectoid ring $R_i$ that is flat over $R_0$.
\end{tcolorbox}
\end{center}

We eventually prove our theorem for general quasi-syntomic rings, but it is a natural step in the proof to first consider the ones that satisfy $\leftmoon$.

The main point of this condition is the following lemma. 

\begin{lem}\label{LemPrismModWptorsionfree}
    Let $R$ a be a flat integral perfectoid $R_0$-algebra, and let $S$ be a fgpp $R$-algebra. Then $\prism_S/[w]$ is $p$-torsion-free and $p,[w]$ is a regular sequence.
\end{lem}

\begin{proof}
    By assumption we have $\prism_S=\prism_{S(m)}\otimes_{A_m} A_{\textup{inf}}(R)$ and this is flat over $A_{\textup{inf}}(R)$ (because by the explicit construction it is clear that $\prism_{S(m)}$ is flat over $A_m$). Being a non-zero divisor does not change after flat base change; therefore, we can assume that $S=R$. Since $[w]$ is a non-zero divisor (it is so in $A_{\textup{inf}}(R)$), it is enough to show that $p$ is a non-zero divisor modulo $[w^p]$. Assume that $px=[w^p]t$ for some $x$ and $t$. We have $[p]_qt=(p\sigma(u)+[w^p])t=p(\sigma(u)t+x)$. This means that the image of $p(\sigma(u)t+x)$ in $R$ is zero. Since in $A_{\textup{inf}}(R)/[p]_q=R$, element $p$ is not a zero-divisor, $[p]_q|\left(\sigma(u)t+x\right)$ and therefore $p|t$, and we have $[w^p]|x$ (since $p$ is not a zero divisor in $A_{\textup{inf}}(R)=W(R^{\flat})$).
\end{proof}

\begin{rmk}\label{pcompleteness}
    Let $S$ be as in the lemma.
    \begin{itemize}
        \item[1)] We have $[p]_q=vp$ in $\prism_S/[w]$ for some unit $v$, and therefore the lemma states that $[w],[p]_q$ is a regular sequence, and therefore for all $n\geq 1$, $[w]^n,[p]_q$ is also a regular sequence.
        
        \item[2)] By definition, $\prism_S$ is $\left(p,[p]_q\right)$-adically complete. Since $[w]^p\in \left(p,[p]_q\right)$, it is also $[w]$-adically complete (cf. \cite[\href{https://stacks.math.columbia.edu/tag/090T}{Lemma  090T}]{stacks-project}). 
        
        \item[3)] $S$ is $[w]$-adically complete. In fact, since $\prism_S$ is $[w]$-adically complete, we need to show that $\CN^1_{\prism_S}$ is $[w]$-adically closed. Let ${x_i}$ be a sequence in $\CN^1_{\prism_S}$ converging $[w]$-adically to an element $x\in \prism_S$. It follows that $\sigma(x_i)$ converges $[w]^p$-adically to $\sigma(x)$. All $\sigma(x_i)$ are divisible by $[p]_q$ and since $[w]^p, [p]_q$ is a regular sequence, $\sigma(x)$ is also divisible by $[p]_q$, and so $x$ belongs to $\CN^1_{\prism_S}$ as desired.
        
        \item[4)] $\prism_S/[w]$ is $p$-adically complete. We have to check that $[w]\prism_S$ is $p$-adically closed. This follows from the fact that $[w],p$ is a regular sequence. 
    \end{itemize}
\end{rmk}

We have that $\delta([w])=0\in [w]\prism_S$, therefore $\prism_S/[w]$ is a $\delta$-ring. By adjunction between $\delta$-rings and Witt vectors (Proposition \ref{DeltaWittAdj}), the natural map $\prism_S/[w]\to S/p$ induces a map $\prism_S/[w]\to W(S/p)$.

\begin{lem}\label{topologicallynilpotent}
    Assume that $S$ satisfies the condition $(\leftmoon)$. Then, the Frobenius $\sigma$ of $\prism_S/[w]$ restricted to the kernel of the map $\prism_S/[w]\to W(S/p)$ is divisible by $p$.
\end{lem}

\begin{proof}
    Note that $\prism_S/[w]$ is a pd-thickening of $S/p$, and so we get a map $$A_{\textup{crys}}(S/p)\to \prism_S/[w]$$ 
    
    One can check that the following diagram is commutative.
        \[
            \begin{tikzcd}
                \prism_S\arrow[rd]\arrow[d]& \\
                A_{\textup{crys}}(S/p)\arrow[r]& \prism_S/[w]\\
            \end{tikzcd}
        \]    
        
    The map $\prism_S\to \prism_S/[w]$ and therefore also $A_{\textup{crys}}(S/p)\to \prism_S/[w]$ are surjective. The Frobenius is divisible by $p$ on the kernel of the map $A_{\textup{crys}}(S/p)\to W(S/p)$. Therefore, it is divisible by $p$ on the kernel of $\prism_S/[w]\to W(S/p)$.
\end{proof}

\section{Prismatic Higher Displays}\label{Secprisdisp}

\begin{dfn}
    Let $(A,(d))$ be an oriented prism with Frobenius $\sigma$ and there exists $\zeta\in A$ such that $\sigma(\zeta)=d$. A \emph{prismatic higher display} on $(A,(d))$ is a datum $\mathcal{D}=\big(P,(P_i)_{i=0}^r,(F_i)_{i=0}^r\big)$, where $P$ is a finite projective $A$-module, \[P=P_0\supset P_1\supset\cdots\supset P_r\] is a filtration of $A$-modules and $F_i:P_i\to P$ are $\sigma$-linear morphisms, subject to the following axioms:
    \begin{enumerate}
        \item $\CN_A^i P\subset P_{i}$
        \item $F_{0|P_i}=d^i F_i$
        \item There exists a decomposition of $A$-modules, called a \emph{normal decomposition}, $P=\bigoplus_{i=0}^{r} L_i$ such that for $0\leq i\leq r$ we have 
        \begin{align}\label{normaldecomposition}
            P_i=\bigoplus_{j=0}^{i}\CN_A^{i-j}L_{j}\oplus\bigoplus_{j=i+1}^rL_j
        \end{align}
        
        and the \emph{divided Frobenius} $${\textup{div}} F:=F_{r|{L_r}}\oplus F_{r-1{|L_{r-1}}}\oplus\cdots \oplus F_{0{|{L_0}}}$$ is a $\sigma$-linear isomorphism. In the special case where $r=1$, the quadruples $(P_0,P_1,F_0,F_1)$ are called \emph{prismatic displays}.
    \end{enumerate}
\end{dfn}

\begin{rmk}
    One can uniquely determine a prismatic higher display by its normal decomposition and divided Frobenius. Indeed, from this data, one can define $P_i$ as in formula (\ref{normaldecomposition}) and $F_i$ by setting for all $0\leq j \leq i$: $$F_{i|_{\CN_A^{i-j}L_j}}=d^{i-j}{\textup{div}} F_{|{\CN_A^{i-j}L_j}}$$
\end{rmk}

\begin{dfn}
    Let $\mu$ be a cocharacter of $\GL_n$, which is identified with a vector $(\mu_1,\dots,\mu_n) \in \mathbb{Z}^n$. We say that a prismatic higher display with a normal decomposition $\oplus L_j$ is of \emph{type $\mu$} if the rank of $L_j$ is equal to the multiplicity of $j$ among $\mu_i$. One can easily see that this is independent of the choice of normal decomposition. 
\end{dfn}

\begin{dfn}
    We call a prismatic higher display \emph{banal} if it has a normal decomposition by free modules.
\end{dfn}

\begin{rmk}\label{changeofbasis1}
    Let $\left(P,(P_i)_{i=0}^r,(F_i)_{i=0}^r\right)$ be a banal prismatic higher display over $(A,I)$. Fix a normal decomposition $P=\bigoplus_{i=0}^{r} L_i$ and a basis $\{v^i_1,\cdots,v^i_{r_i}\}$ for each $L_i$. As we mentioned above, $F_i$ can be decoded in the invertible matrix ${\textup{div}} F$. Now, if we change the decomposition and the basis in such a way that the new decomposition also gives a normal decomposition, and $B=(b_{ij})$ is the change of basis matrix, then for all $0\leq m\leq r$ and all $i>r_0+r_1+\ldots+r_m$, we have $b_{ij}\in \CN^{m}_A$.
\end{rmk}

\begin{rmk}\label{changeofbasis2}
    There is a coordinate-free way to state this. Write $\GL_{n,\mathbb{Z}}=\Spec\, R$ and let $\mu(z)=(z^{i_1},z^{i_2},\ldots,z^{i_n})$ be a cocharacter of $\GL_n$. The action of $\mathbb{G}_m$ on $\GL_n$ by conjugation via $\mu$ induces an action on $R$. This action induces a filtration $\mathcal{F}^i R$ on $R$, where $\mathcal{F}^s R=\mathcal{F}^{s-1} R$ unless we have $s=i_m-i_l$ for some $l,m\in\{1,\dots,n\}$. If we change the coordinates by a matrix $B$, then our condition on $B$ is equivalent to the condition that $B\in \GL_n(A)=\Hom(R,A)$ is given by $k:R\to A$ such that if $x\in \mathcal{F}^i R$, then $k(x)\in \CN^i_A$. Under this change of coordinates, the map ${\textup{div}} F$ changes to $$B^{-1}\cdot {\textup{div}} F\cdot \left(\mu(d) \sigma(B) \mu(d^{-1})\right)$$
\end{rmk}

\begin{dfn}\label{Breuil-Kisin-Fargues-Module}
    A \emph{Breuil-Kisin-Fargues module} over a prism $(A,(d))$ is a pair $(M,F)$, where $M$ is a finitely generated projective $A$-module and $F:\sigma^*M\to M$ is an $A$-linear morphism such that the base change of $F$ to $A[\frac{1}{I}]$ is an isomorphism. We call a Breuil-Kisin-Fargues module \emph{banal} if $M$ is free. 
\end{dfn}
    
\begin{lem}
    Let $(A,(d))$ be an oriented prism and assume that there exists $\zeta\in A$ such that $\sigma(\zeta)=d$. Let $\mathcal{D}=\big(P,(P_i)_{i=0}^r,(F_i)_{i=0}^r\big)$ be a prismatic higher display over $(A,(d))$. The pair $(P,F_0)$ is a Breuil-Kisin-Fargues module over $(A,(d))$, and the functor from the category of prismatic higher displays over $(A,(d))$ to the category of Breuil-Kisin-Fargues modules over $(A,(d))$ is fully faithful.
\end{lem}

\begin{proof}
    To show that $F_0$ is a $\sigma$-linear isomorphism over $P[1/d]$ note that if $x=x_r+x_{r-1}+\cdots+x_0$, we have: $$F_0(x)={\textup{div}}F\left(x_0+\frac{x_1}{\zeta}+\cdots+\frac{x_r}{\zeta^r}\right)$$
    
    Since $d$ is a non-zero-divisor, $F_0$ uniquely determines all $F_i$, and so we need to show that the filtration $(P_i)$ is uniquely determined by the pair $(P_0,F_0)$. Set $$P'_i=\left\{x\in P_0\,|\,F_0x\in d^iP_0\right\}$$ 
    
    We claim that $P_i=P'_i$. The inclusion $P_i\subset P'_i$ follows from condition 2 of the definition of a prismatic higher display. Now, if $x\in P'_i\setminus P_i$, write $x=(x_1,x_2,\cdots,x_r)$ where $x_j\in L_j$ (using the normal decomposition). There should be some $j$ such that $x_j\not\in \CN^{i-j}_AL_j$, therefore $x_j\otimes 1\not\in d^{i-j}\sigma^* L_j$. By the last condition of the definition of a prismatic higher display, ${\textup{div}} F:\sigma^{*}P_0/d^{i-j}\to P_0/{d^{i-j}}$ is an isomorphism, but ${\textup{div}} F(x_j)$ is zero, which is a contradiction.   

    It is clear that this functor is faithful. To see that it is full, note that a morphism of Breuil-Kisin-Fargues modules is a map $P_0\to P'_0$ that is compatible with the Frobenius on both sides. If $(P_0,F_0)$ and $(P'_0,F'_0)$ come from prismatic higher displays, we just define the flirtations and divided Frobenius in terms of $F_0,F'_0$ and so the map is also compatible with these structures and therefore is a morphism of prismatic higher displays.
\end{proof}

\begin{rmk}
    The Breuil-Kisin-Fargues module $(P, F_0)$ is called the \emph{underlying Breuil-Kisin-Fargues module} of the prismatic higher display $\big(P,(P_i)_{i=0}^r,(F_i)_{i=0}^r\big)$. 
\end{rmk}

So, by the lemma, prismatic higher displays are Breuil-Kisin-Fargues modules that satisfy other conditions. We can make these conditions explicit.

\begin{dfn}
Let $(P,F)$ be a Breuil-Kisin-Fargues module over an oriented prism $(A,(d))$. We call $(P, F)$ \emph{admissible} if and only if the following three conditions are satisfied:

    \begin{enumerate}
        \item We have  $dP\subset FP$
        \item The $A/\CN_A^1$-sub module $E$ of $P/d$ generated by the image of the composition $P\xrightarrow{F} P\to P/d$ is projective.
        \item $E\otimes_{A/\CN^1_A,\sigma} A/dA\to P_0/dP_0$ is injective.
    \end{enumerate}

\end{dfn}
\begin{rmk}
    If $R$ is a qsrp ring, an admissible Breuil-Kisin-Fargues module over $\prism_R$ is the same as an admissible Dieudonn\'e module over $R$ in the language of \cite{anschutz2023prismatic}.
\end{rmk}

\begin{pro}\label{char}
    Let $(M,F)$ be a Breuil-Kisin-Fargues module over an oriented prism $(A,(d))$. Assume that $d$ has a $\sigma$-preimage. Then $(M, F)$ is the underlying Breuil-Kisin-Fargues module of a prismatic display if and only if it is admissible. Therefore, admissible Breuil-Kisin-Fargues modules are the same as prismatic displays.
\end{pro}

\begin{proof}
    This is \cite[Proposition 4.29]{anschutz2023prismatic}. Although that proposition is only stated for qrsp rings $R$ and their initial prism $\prism_R$, they don't use anything about these rings anywhere in the proof.
\end{proof}

An alternative method exists for defining prismatic higher displays over a prism. Let $(A,(d))$ be an oriented prism. Assume that $d$ has a $\sigma$-preimage. Define $\underline{A}$ as the graded ring $\oplus_{i\ge 0} N^iA$. Consider the morphisms $\underline{\sigma}, \underline{\tau}:\underline{A}\to A$ where for $x\in N^iA$, $\underline{\sigma}(x)= \frac{\sigma(x)}{d^i}$ and $\underline{\tau}(x)=x$.

\begin{dfn}
    A \emph{higher display} over $\underline{A}$ is a graded projective module $M$ over $\underline{A}$ together with an isomorphism $F:\underline{\sigma}^*M \xrightarrow{\sim}
 \underline{\tau}^* M$.
\end{dfn}

\begin{rmk}\label{alternatedefinition}
    The argument of \cite[\S 3]{lau2021higher} shows that higher displays over $\underline{A}$ are the same as prismatic higher displays over $(A,(d))$.
\end{rmk}

\begin{rmk}\label{prismaticdiudennetheory}
    The main result of \cite{anschutz2023prismatic} states that if $(A,I)$ is the final object of the prismatic site of a quasi-regular semiperfectoid ring $R$, then the category of admissible Breuil-Kisin-Fargues modules on $A$ is equivalent to the category of $p$-divisible groups on $R$. Therefore, a prismatic display can be associated with a $p$-divisible group or an Abelian variety.
\end{rmk}

Let us collect some useful facts about prismatic higher displays.

\begin{pro}
    Let $\kappa:(A,(d))\to (A',(d))$ be a morphism of oriented prisms and $\mathcal{D}$ a prismatic higher display on $(A,(d))$, with $\mathcal{D}=\big(P,(P_i)_{i=0}^r,(F_i)_{i=0}^r\big)$. Then there is a unique prismatic higher display $\mathcal{D}_{A'}$ on $(A',(d))$ with the underlying Breuil-Kisin-Fargues module $(P_{A'},F_{0,A'})$, called the \emph{base change} of $\mathcal{D}$ along $\kappa$. Here $P_{A'}$ is the base change $P\otimes_{A,\kappa}A'$ and $F_{0,A'}$ is $F_0\otimes \sigma_{A'}$.
\end{pro}

\begin{proof}
    Let $P=\oplus L_i$ be a normal decomposition. Consider the  prismatic higher display defined by the normal decomposition $P_{A'}=\oplus (L_{i}\otimes A')$ and divided Frobenius $(\textup{div}F)_{A'}$. We have to show that this is independent of the choice of the normal decomposition for $P$. To see this note that for all $i$,
        \begin{multline*}
            \CN_{A'}^iP_0+\CN_{A'}^{i-1} (P_1\otimes A')+\cdots +\CN^1_{A'} (P_{i-1}\otimes A')+P_i\otimes A'+\cdots+P_r\otimes A' =\\
            \CN_{A'}^{i}(L_0\otimes A')+ \CN_{A'}^{i-1}(L_1\otimes A') +\cdots+ \CN^1_{A'}(L_{i-1}\otimes A')+ L_i\otimes A'+ \cdots+ L_r\otimes A'
        \end{multline*}

    Therefore $P_{A',i}$ is independent of the choice of $L_i$.
\end{proof}

\begin{rmk}\label{monoidalstructure}
    Similarly to the classical theory of higher displays, one can define tensor products, exterior powers, etc. using normal decompositions. In particular, the category of prismatic higher displays over a prism is a tensor category.
\end{rmk}

\begin{dfn}
    The fibered category of prismatic higher displays over the site $\textbf{QRSPerfd}$ is the category whose objects over $R$ are the prismatic higher displays over $(\prism_R,\BI_R)$ and whose morphisms are morphisms between prismatic higher displays.
\end{dfn}

\begin{pro}\label{decentofdisplay}
    The fibered category of prismatic higher displays over the site $\textbf{QRSPerfd}$ is a stack.
\end{pro}

\begin{proof}
    The statement for admissible Breuil-Kisin-Fargues modules is \cite[Proposition 4.9]{anschutz2023prismatic}. By Proposition \ref{char} we get the statement for prismatic displays. The proof for prismatic higher displays is similar: by Remark \ref{alternatedefinition}, it reduces to the descent of projective modules and morphisms between them.
\end{proof}

\begin{pro}\label{banalisation}
    Every prismatic higher display becomes banal after going to a suitable quasi-syntomic cover.
\end{pro}

\begin{proof}
    Choose a normal decomposition $ P_0=\oplus L_i$ for the higher display and consider the base change of this normal decomposition to $R$ to obtain projective $R$-modules $\bar{L}_i$. One can find an \'etale cover $S$ of $R$ such that the base change of these projective $R$ modules becomes free. When we base change the higher display to $S'$, $\oplus (L_i\otimes \prism_{S'})$ is a normal decomposition, and its base change to $S'$ is given by $\bar{L}_i\otimes S'$. The pair $(\prism_{S'},\CN_{S'}^1)$ is Henselian, and the $\bar{L}_i\otimes S'$ are trivial. Therefore, the projective modules $L_i\otimes \prism_{S'}$ are free.
\end{proof}

\begin{dfn}
    The stack of prismatic higher displays over $\QSyn$ is the image of the stack of prismatic higher displays over $\textbf{QRSPerfd}$ under the equivalence of Proposition \ref{equivalenceofsites}.
\end{dfn}

\section{Prismatic Higher \texorpdfstring{$(G,\mu)$}{(G,mu)}-Displays}\label{secGdipslay}

\subsection{Banal Objects}\label{SectionBanalObjects}

We can generalize the construction in the last section to define higher displays with $G$-structure, where $G$ is a reductive group. 

\begin{dfn}\label{DefDispDatum}
    A \emph{display datum} over $\BF_q$ is a pair $(G,\mu)$, where $G$ is a reductive group scheme over $\BZ_p$ and $\mu:\BG_{m,\BZ_q}\to G_{\BZ_q}$ is a cocharacter.
\end{dfn}

 Let $(G,\mu)$ be a display datum over $\BF_q$. The action of $\mathbb{G}_m$ on $G_{\mathbb{Z}_q}$ via $\mu$ induces an action of $\mathbb{G}_m$ on $\mathcal{O}_{G}\otimes_{\mathbb{Z}_p}\mathbb{Z}_q$. Therefore, we get a grading 
    \begin{equation}\label{EqGradDispGrp}
        \mathcal{O}_{G}\otimes_{\mathbb{Z}_p}\mathbb{Z}_q=\oplus_{i=-\infty}^{\infty} \CO_{G,i}
    \end{equation}
    
Consider the decreasing filtration $\mathcal{F}^k=\oplus_{i=k}^\infty \CO_{G,i}$ for $ k\in\mathbb{Z}$.

\begin{dfn}\label{DefDispGrp}
    The \emph{display group} $G^\mu$ is the functor from the category of $\mathbb{Z}_q$-prisms to the category of sets that sends a prism $(A,I)$ to the subset of $G(A)=\Hom(\mathcal{O}_{G}\otimes_{\mathbb{Z}_p}\mathbb{Z}_q,A)$ consisting of elements $k$ satisfying $k(\mathcal{F}^i)\subset \CN_A^i$.
\end{dfn}

\begin{rmk}
    One can also work with a cocharacter $\mu$ defined over a ramified extension $\mathcal{O}_E$ of $\BZ_p$. In that case, we should work with $\mathcal{O}_E$-prisms, as defined in \cite{ito2023prismatic}. Everything we prove in this paper generalizes to this situation, but for the sake of simplicity, we choose to work with the more familiar language of (unramified) prisms.
\end{rmk}

\begin{rmk}
    One can easily check that $G^\mu(A,I)$ is a subgroup of $G(A)$. In fact, in the case of $G=\GL_n$, one can prove it by a direct calculation, and in general, one can embed $G$ in $\GL_n$ and observe that by definition $G^\mu(A)=\GL_n^\mu(A)\cap G(A)$. So, $G^\mu$ is a subgroup functor of $G$. 
\end{rmk}

Assume that $I=(d)$. Consider the map $\Phi_d^\mu:G^\mu(A,(d))\to G(A[\frac{1}{d}])$ defined as
$$\Phi_d^\mu(k)=\mu(d)\sigma(k)\mu(d^{-1})$$
It is not hard to check that this map is a group homomorphism and factors through $G(A)\to G(A[\frac{1}{d}])$. One can show this by a direct calculation for $\GL_n$ and reduce the general case to this case. In fact, we have an identification 
\begin{equation}
    \label{AltDefDispGrp}
    G^\mu(A,I)=\left\{k\in G(A)\,|\,\Phi_d^\mu(k)\in G(A)\right\}
\end{equation}

For more details on the display group, see Appendix \ref{proofofstructure}.

\begin{dfn}
    Let $(A,I)$ be an orientable prism and choose an orientation $d$. Define the map $\Phi_d^\mu:G^{\mu}(A,I)\to G(A)$ as above.
\end{dfn}
    
Now we can define the groupoid of banal $(G,\mu)$-displays. Let $(A,I)$ be an orientable prism. Choose a generator $d$ of $I$, and denote by $G(A)_d$ the group $G(A)$ with the $G^{\mu}(A,I)$-action via the inclusion and $\Phi^{\mu}_d$, that is, for $k\in G^{\mu}(A,I)$ and $X\in G(A)$, we set $k\cdot X:=k^{-1}X\Phi^{\mu}_d(k)$. Now choose another generator $e$ of $I$ and write it as $vd$ with $v\in A^\times$. The map
    \begin{align*}
       \alpha_{e/d}: G(A)_d \to &\, G(A)_{e}\\
       X \mapsto &\, X\mu(v^{-1})
    \end{align*}   
is a $G^{\mu}(A,I)$-equivariant bijection: 
\begin{multline*}
    \alpha_{e/d}(k\cdot X)=\left( k^{-1}X\mu(d)\sigma(k)\mu(d^{-1}) \right).\mu(v^{-1}) = \\ k^{-1}\left(X\mu(v^{-1})\right)\mu(vd)\sigma(k)\mu\left((vd)^{-1}\right)=k\cdot \alpha_{e/d}(X)
\end{multline*}

So, we have a well-defined and canonical action of $G^{\mu}(A,I)$ on $G(A,I):=\varprojlim_d G(A)_d$, where the limit is on the set of generators of $I$ and the transition morphisms are the $\alpha$ defined above.

\begin{dfn}
    Let $(A,I)$ be an orientable prism. The groupoid of \emph{banal prismatic higher $(G,\mu)$-displays} over $(A,I)$ is the quotient groupoid
    $$\mathcal{B}_{(G,\mu)}^{\textup{ban}}(A,I):=\left[G(A,I):G^\mu(A,I)\right]$$
    
    We call morphisms in this category \emph{banal morphisms}. So, a banal prismatic higher $(G,\mu)$-display is represented by an element of $G(A)$, and a banal morphism is an element of $G^\mu(A,I)$. Although $G^\mu(A,I)$ is a subset of $G(A)$, we employ uppercase letters for objects and lowercase letters for morphisms, to improve clarity.
\end{dfn}

\begin{rmk}
     If $\mu$ is minuscule, $R$ is a ring, $A=W(R)$ and $d=p$, banal prismatic higher $(G,\mu)$-displays over $(A,(p))$ are the same as banal $(G,\mu)$-displays over $R$ as defined in \cite{bueltel2020g}.
\end{rmk}

\begin{lem}\label{BanalG-muFunc}
    Let $(A,I)$ be an orientable prism and $(A,I)\to (B,J)$ be a map of prisms. Then the canonical map $G(A)\to G(B)$ induces a functor \[\CB_{(G,\mu)}^{\textup{ban}}(A,I)\to \CB_{(G,\mu)}^{\textup{ban}}(B,J)\]
\end{lem}

\begin{proof}
    Note that by Rigidity (Proposition \ref{rigidity}) the image under $A\to B$ of any generator $d$ of $I$ is a generator of $J$. The following diagram is commutative:
     \[ \begin{tikzcd}
                G^{\mu}(A,I) \arrow{r}{\Phi^{\mu}_d}\arrow{d}& G(A)\arrow{d}\\
                G^{\mu}(B,J) \arrow{r}[swap]{\Phi^{\mu}_d}& G(B)\\
        \end{tikzcd} \]

    Therefore, we have a natural map $G(A,I)\to G(B,J)$ that is equivariant for the induced map $G^{\mu}(A,I)\to G^{\mu}(B,J)$ and the statement of the lemma follows.
\end{proof}

\begin{dfn}
    Let $(A,I)$ be a prism. A \emph{Breuil-Kisin-Fargues module with $G$-structure} over $(A,I)$ is a pair $(\SM,F)$, where $\SM$ is a $G$-torsor over $A$ and $F:\sigma^*\SM\to \SM$ is a morphism over $A$ such that the base change of $F$ to $A[\frac{1}{I}]$ is an isomorphism. We call a Breuil-Kisin-Fargues module with $G$-structure \emph{banal} if $\SM$ is the trivial torsor. A morphism of Breuil-Kisin-Fargues modules is a morphism of torsors that commutes with $F$.
\end{dfn}

\begin{rmk}
    Note that after fixing a trivialization, a banal Breuil-Kisin-Fargues module with $G$-structure is given by an element of $G(A[\frac{1}{I}])$. Therefore, the groupoid of banal Breuil-Kisin-Fargues modules with $G$-structure is the quotient groupoid $$\left[G\left(A[\frac{1}{I}]\right):G(A)\right]$$ where $k\in G(A)$ acts on $G(A[\frac{1}{I}])$ by $k.X=k^{-1}X\sigma(k)$
\end{rmk}

Let $(A,I)$ be an orientable prism and choose an orientation $d$. Then we have $A[\frac{1}{I}]=A[\frac{1}{d}]$. We have a map $G(A,I)\to G\left(A[\frac{1}{d}]\right)$ sending a sequence $(X_e)\in G(A,I)$ to $X_d\mu(d)$. This map is independent of the choice of a generator.

\begin{pro}\label{banalfargueskisinmodule}
    Let $(A,I)$ be an orientable prism. The above map induces a fully faithful functor
    \[\upbeta:\CB_{(G,\mu)}^{\textup{ban}}(A,I)\to \left[G\left(A[\frac{1}{I}]\right):G(A)\right]\] from the groupoid of banal prismatic higher $(G,\mu)$-displays over $(A,I)$ to the groupoid of banal Breuil-Kisin-Fargues modules with $G$-structures over $(A,I)$.
\end{pro}

\begin{proof}
    Choose an orientation $I=(d)$. We need to show that if $k\in G^\mu(A,I)$ is a morphism between displays $X$ and $Y$, then it also defines a morphism between Breuil-Kisin-Fargues modules $X\mu(d)$ and $Y\mu(d)$: $k$ being a morphism between $X$ and $Y$ means that $$Y=k^{-1}X\mu(d)\sigma(k)\mu(d^{-1})$$ and this implies that $$Y\mu(d)=k^{-1}X\mu(d)\sigma(k)$$

    So $k$ is a morphism between $Y\mu(d)$ and $X\mu(d)$. By definition, this functor is faithful. To see that it is full, note that if $$Y\mu(d)=k^{-1}X\mu(d)\sigma(k)$$ then we have
    \[\Phi_d^\mu(k)=\mu(d)\sigma(k)\mu(d)^{-1}=X^{-1}kY\] and so $\Phi_d^\mu(k)\in G(A)$, using (\ref{AltDefDispGrp}), this implies that $k\in G^\mu(A,I)$.
\end{proof}

\subsection{The Stack}\label{StackprismaticGdisplays}

Fix a display datum $(G,\mu)$ and assume that $\mu$ is defined over $\mathbb{Z}_q$. To discuss the stack of prismatic displays, we should first have a site. One option is to consider the prismatic site over a fixed prism. So, fix an orientable $\mathbb{Z}_q$-prism $(A,I)$, and consider the relative prismatic site $(A,I)_{\prism}$ with the $(p,d)$-complete \'etale topology.

\begin{dfn}\label{DefPrimStackPrimDisp}
    Define a presheaf of groupoids by sending a prism $(C,J)$ over $(A,I)$ to $\mathcal{B}_{(G,\mu)}^{\textup{ban}}(C,J)$. This presheaf is, in fact, a prestack (because both objects and morphisms are sheaves).  We define the stack of prismatic higher $(G,\mu)$-displays on $(A,I)_{\prism}$ to be the sheafification of this presheaf and denote it by $\mathcal{B}^{\prism_{/(A,I)}}_{(G,\mu)}$.
\end{dfn}

\begin{rmk}
    One can give an equivalent definition of prismatic higher $(G,\mu)$-displays in the language of Lau's paper, but because the ideal $I$ in the definition of a prism is always a Cartier divisor, we believe that it is simpler to work with this definition.
\end{rmk}

We can also work over a site of geometric objects, where (at least over a basis of the site) there is a canonical prism over the objects. For example, we can take the quasi-syntomic site where quasi-regular semiperfectoid rings form a basis, and their prismatic site has a final object. Similarly, one can consider the site of perfectoid spaces with \'etale topology.

\begin{dfn}\label{DefPrisHighDisp}
    \begin{itemize}
        \item[(i)] Over the site $\textbf{QRSPerfd}/\mathbb{Z}_q$, consider the presheaf of groupoids that sends a qrsp ring $R$ to $\CB_{(G,\mu)}^{\textup{ban}}(\prism_R,\BI_R)$ and sends a map $R\to R'$ to the natural map $\CB_{(G,\mu)}^{\textup{ban}}(\prism_R,\BI_R)\to \CB_{(G,\mu)}^{\textup{ban}}(\prism_{R'},\BI_{R'})$, given by Lemma \ref{BanalG-muFunc} (note that $\BI_R$ is always principal for a qrsp ring).
        We denote the sheafification of this presheaf by $\mathcal{B}^{\textup{qrsp}}_{(G,\mu)}$.
        
        \item[(ii)] We define the stack of prismatic higher $(G,\mu)$-displays on $\QSyn/\BZ_q$, denoted by $\mathcal{B}^{\textup{qsyn}}_{(G,\mu)}$, as the image of $\mathcal{B}^{\textup{qrsp}}_{(G,\mu)}$ under the equivalence of Proposition \ref{equivalenceofsites}.
        
        \item[(iii)] Over the site $\Perf /\Spd \mathbb{Z}_q$, consider a presheaf of groupoids that sends $\Spa(R,R^+)\to \Spd \mathbb{Z}_q$, i.e., an untilt $\Spa(R^\sharp,R^{\sharp+})$, to $\CB_{(G,\mu)}(A_{\textup{inf}}(R^+),\kernel \theta_{R^{\sharp+}})$. We denote the sheafification of this presheaf by $\CB^{\textup{perf}}_{(G,\mu)}$.
    \end{itemize}
\end{dfn}

\begin{pro}\label{fargueskisinmodule}
    Assume that $(R,R^+)$ is an affinoid perfectoid over $\mathbb{Z}_q$. There is a fully-faithful functor from the groupoid $\CB^{\textup{perf}}_{(G,\mu)}(R,R^+)$ to the groupoid of Breuil-Kisin-Fargues modules over $A_{\textup{inf}}(R^+)$.
\end{pro}

\begin{proof}
    Because both sides are stacks, it is enough to check this for banal objects, which is Proposition \ref{banalfargueskisinmodule}.
\end{proof}

\begin{example}
    When $G=\GL_n$ and we are working with the quasi-syntomic site, by Propositions \ref{decentofdisplay} and \ref{banalisation} the $(\GL_n,\mu)$-displays over a qrsp ring $R$ are the same as the admissible prismatic Dieudonn\'e crystal with Hodge polygon $\mu$ over $R$.
\end{example}

\begin{rmk}\label{RemVTop}
    By \cite{ito2023prismatic} the stack $\mathcal{B}^{\perf}_{(G,\mu)}$ is also a stack in the v-topology.
\end{rmk}

We want to describe the stacks $\CB^{\textup{perf}}_{(G,\mu)}$ and $\CB^{\textup{qsyn}}_{(G,\mu)}$ as quotient groupoids.

\begin{dfn}
    Let us define the \emph{positive loop groups} and the associated display groups:
        \begin{itemize}
            \item[(i)] 
                Define presheaves $\mathscr{L}^+G$ and $\mathscr{L}^+G^\mu$ over $\textbf{QRSPerfd}/\mathbb{Z}_q$ by sending a qrsp ring $R$ to $G(\prism_R,\BI_R)$ and $G^{\mu}(\prism_R,\BI_R)$ respectively.
        
            \item[(ii)]
                Define functors $L^+G$ and $L^+G^\mu$ on affinoids in $\Perf/\Spd \mathbb{Z}_q$ by sending $\Spa (R,R^+)$, respectively, to the groups $G\left(A_{\textup{inf}}(R^{\sharp+}),\kernel \theta_{R^+})\right)$ and $G^\mu\left(A_{\textup{inf}}(R^{\sharp+}),\kernel \theta_{R^+}\right)$.
        \end{itemize}
\end{dfn}

\begin{lem}
    Presheaves $\mathscr{L}^+G$ and $\mathscr{L}^+G^\mu$ are sheaves for the quasi-syntomic topology and define sheaves on $\QSyn/\mathbb{Z}_q$. Similarly, $L^+G$ and $L^+G^\mu$ induce v-sheaves on $\Perf/\Spd \mathbb{Z}_q$.
\end{lem}

\begin{proof}
    The functor $R\mapsto G(\prism_R,\BI_R)$ is a sheaf on $\textbf{QRSPerfd}/\mathbb{Z}_q$ because the functor $R\mapsto \prism_R$ is a quasi-syntomic sheaf (\cite[Propsition 7.11]{anschutz2023prismatic}). Now, consider the functor $R\mapsto G^\mu(\prism_R,\BI_R)$. For $G=\GL_n$, this is again a sheaf, since the functors $R\mapsto \CN_{\prism_R}^i$ are sheaves. For an arbitrary $G$, this follows from the fact that if $G\into \GL_n$, then $G^\mu=G\cap \GL_n^\mu$. Using the fact that qrsp rings form a basis for $\QSyn$ (cf. Propositions \ref{propertiesofqrsp} and \ref{equivalenceofsites}), we obtain sheaves on $\QSyn/\mathbb{Z}_q$. The second statement is similar.
\end{proof}
\begin{rmk}
    Note that by the discussion in the previous section, there is an action of $L^+G^{\mu}$ (resp. $\mathscr{L}^+G^\mu$) on $L^+G$ (resp. $\mathscr{L}^+G$).
\end{rmk}

Let us now give an alternative definition of prismatic higher $(G,\mu)$-displays, using torsors:

\begin{pro}\label{PropDispTors}
    \begin{itemize}
        \item[(1)] Let $R$ be a quasi-syntomic $\mathbb{Z}_q$-algebra. Giving a prismatic higher $(G,\mu)$-display over $R$ is the same as giving a pair $(\SQ,\alpha)$, where $\SQ$ is a quasi-syntomic torsor for the sheaf $\mathscr{L}^+G{^\mu}$ over $R$ and $\alpha:\SQ\to \mathscr{L}^+G$ is a $\mathscr{L}^+G^\mu$-equivariant morphism.
        
        \item[(2)] Let $S$ be an object of $\Perf/\Spd \mathbb{Z}_q$. Giving a prismatic higher $(G,\mu)$-display over $S$ is the same as giving a pair $(\SQ,\alpha)$, where $\SQ$ is torsor for the sheaf $L^+G^\mu$ over $S$ and $\alpha:\SQ\to L^+G$ is a $L^+G^\mu$-equivariant morphism.
    \end{itemize}
\end{pro}

\begin{proof}
    We prove the first statement. The proof of the second statement is similar. First, assume that we have a pair $(\SQ,\alpha)$ as in the statement. By definition and Proposition \ref{propertiesofqrsp}, there is a qrsp ring $R'$ that covers $R$ and over which $\SQ$ is trivial. So $\alpha$ is given by an element $X\in \mathscr{L}^+G(R')=G(\prism_{R'},\BI_{R'})$. If we change the trivialization of $\SQ$ over $R'$, $X$ changes to $k\cdot X$, where $k$ is an element of $\mathscr{L}^+G^\mu(R')=G^\mu(\prism_{R'},\BI_{R'})$. This means that we have a well-defined object of the quotient groupoid $\left[G(\prism_{R'},\BI_{R'}):G^\mu(\prism_{R'},\BI_{R'})\right]=\CB_{(G,\mu)}^{\textup{ban}}(R')$. Now, since the pair $(\SQ,\alpha)$ is defined over $R$, we can descend to obtain an object of $\CB_{(G,\mu)}^{\textup{qsyn}}(R)$.

    Conversely, if we have a prismatic higher $(G,\mu)$-display $\CD$ over $R$, there exists a quasi-syntomic cover $R\to R'$ with $R'$ qrsp and such that the base change $\CD_{R'}$ belongs to $[G(\prism_{R'},\BI_{R'}):G^\mu(\prism_{R'},\BI_{R'})]$. This means that we have an element $X$ of $\mathscr{L}^+G(R')$, together with a descent datum $k$ as a section of $\mathscr{L}^+G^\mu$. This implies that we can use $k$ to descend the trivial $\mathscr{L}^+G^\mu$-torsor over $R'$ to obtain a $\mathscr{L}^+G^\mu$-torsor over $R$ and descend $X$ (which is a $\mathscr{L}^+G^\mu$-equivariant morphism from the trivial torsor to $\mathscr{L}^+G$) to obtain a $\mathscr{L}^+G^\mu$-equivariant morphism to $\mathscr{L}^+G$ as desired.
\end{proof}

\begin{cor}\label{assoctiedtorsor}
    Let $R$ be a quasi-syntomic $\mathbb{Z}_q$-algebra. There is a canonical functor from the category of prismatic higher $(G,\mu)$-displays over $R$ to the category of $P_{\mu}$-torsors over $R$. There is a similar statement for prismatic higher $(G,\mu)$-displays on $\Perf/\Spd \mathbb{Z}_q$ and $(A,I)_{\prism}$, where $(A,I)$ is a prism.
\end{cor}

\begin{proof}
    By Proposition \ref{structure}, we have a canonical morphism of quasi-syntomic sheaves $\mathscr{L}^+G^\mu\to P_{\mu}$. The statement now follows from the previous proposition. 
\end{proof}

We need the following definition for the deformation theory of prismatic displays.

\begin{dfn}\label{DefHodgeFilt}
    Let $R$ be a $\BZ_q$-algebra. A \emph{$(G,\mu)$-Hodge filtration} over $R$ is a pair $(\SQ,\SM,\tau)$ where $\SQ$ is a $P_{\mu}$-torsor over $R$, $\SM$ is a $G$-torsor over $R$ and $\tau:\SQ*^{P_{\mu}}G\to \SM$ is a $G$-equivariant isomorphism.
\end{dfn}

\begin{rmk}\label{RemStackofHodgeFilt}
   \begin{itemize}
       \item[1)]  One can easily see that, restricted to $\QSyn/\mathbb{Z}_q$, the category of $(G,\mu)$-Hodge filtrations is a stack over $\QSyn/\mathbb{Z}_q$. We denote this stack by $\SH_{(G,\mu)}$.
       
       \item[2)] Using Corollary \ref{assoctiedtorsor}, we obtain a functor $\FH:\mathcal{B}^{\textup{qsyn}}_{(G,\mu)}\to \SH_{(G,\mu)}$: Let $\CD$ be a prismatic higher $(G,\mu)$-display and let $\SQ$ be the associated $P_{\mu}$-torsor. We send $\CD$ to the triple $(\SQ,\SQ*^{P_{\mu}}G,\id)$. 
   \end{itemize}
\end{rmk}

\begin{rmk}
    The deformation theory of $(G,\mu)$-Hodge filtrations is nice: if $R'\onto R$ is a nilpotent thickening and $(\SQ,\SM,\tau)$ is a $(G,\mu)$-Hodge filtration over $R$, then $\SQ$ and $\SM$ have unique lifts to $R'$. Therefore, we only have to lift $\tau$, which is a point in $\Gr(G,\mu)=G/P_{\mu}$, and so we have the deformation problem of lifting a point in a Grassmannian.
\end{rmk} 

\begin{rmk}
    We also have the Tannakian viewpoint; identically to the proof of the main result of \cite{daniels2021tannakian}, one can prove that there is an equivalence of categories between the category of $(G,\mu)$-displays and the category of fiber functors from the category of representations of $G$ to the category of higher displays of type $\mu$. Similarly, we have a Tannakian viewpoint for the Breuil-Kisin-Fargues modules.
\end{rmk}

\begin{dfn}
    Assume that $R$ is a qrsp $\mathbb{Z}_q$-algebra. Let $\textup{PF}(R)$ be the category of finitely generated projective $\prism_R$-modules $P$ with a filtration \[P=P_0\supset P_1\supset\cdots\] by sub-$\prism_R$-modules, such that $P_i$ come from a normal decomposition in the sense of (\ref{normaldecomposition}). Consider the natural tensor structure on this category (Remark \ref{monoidalstructure}). Consider the fiber functor $\FF_\mu:\textup{Rep}_{\mathbb{Z}_q}(G)\to \textup{PF}(R)$  that sends a representation $(V,\rho)$ to the projective module $V\otimes_{\BZ_q} \prism_R$ with the filtration induced from the normal decomposition $V\otimes_{\BZ_q} \prism_R=\oplus (V_i\otimes_{\BZ_q} \prism_R)$ where $V=\oplus V_i$ is the decomposition of $V$ given by the action of $\BG_m$ via $\mu$. We say that a fiber functor $\FF:\textup{Rep}_{\mathbb{Z}_q}(G)\to \textup{PF}(R)$ is of \emph{type $\mu$} if after passing to some cover of $R$ it becomes isomorphic to $\FF_{\mu}$. Let $\textup{PFF}^\mu(R)$ be the category of fiber functors of type $\mu$.  
\end{dfn}

\begin{pro}\cite[Corollary 3.11]{daniels2021tannakian}\label{tannakianformalism}
    There is an equivalence of categories between the category of $\mathscr{L}^+G^\mu$-torsors over $R$ and $\textup{PFF}^{\mu}(R)$.
\end{pro}

One can give an explicit description of the $\GL_n^\mu$-torsor associated with a prismatic higher display of type $\mu$. Let $P=\oplus L_i$ be a normal decomposition. Consider $T=\oplus F_i$, where $F_i$ is a free module of the same rank as $L_i$. Define $$T_i=\CN^{i}_{\prism_R}F_0\oplus \CN^{i-1}_{\prism_R}F_1 \oplus\cdots\oplus \CN^1_{\prism_R}F_{i-1}\oplus F_i\oplus F_{i+1}\oplus\cdots\oplus F_d$$ 

It is easy to see that $\underline{\Isom}((P_i),(T_i))$, the sheaf of linear isomorphisms from $P$ to $T$ that sends $P_i$ to $T_i$, is a $\GL_n^\mu$-pseudo-torsor. By Proposition \ref{banalisation} this is actually a $\GL_n^\mu$-torsor. It is easy to check that under the equivalence between prismatic higher displays of type $\mu$ and $(\GL_n,\mu)$-displays, this $\GL_n^\mu$-torsor is the $\GL_n^\mu$-torsor associated with $P$.

\begin{pro}\label{banality}
    Let $R$ be a qrsp $\mathbb{Z}_q$-algebra. An element of $\CB^{\textup{qsyn}}_{(G,\mu)}(R)$ is banal if and only if the associated $P_{\mu}$-torsor (over $R$) is trivial.
\end{pro}

\begin{proof}
    By Proposition \ref{tannakianformalism}, we can reduce to the case of $\GL_n$. In that case, having a $P_{\mu}$-torsor is equivalent to having a filtration $\bar{P}=\oplus \bar{L}_i$. Since the pair $(\prism_R,\CN_{\prism_R}^1)$ is Henselian, one can lift this filtration to a normal decomposition $P=\oplus L_i$. If the $P_{\mu}$-torsor is trivial, then all $L_i$ are free and the associated $G^\mu$-torsor $\underline{\Isom}((P_i),(T_i))$ is also trivial.
\end{proof}

\begin{cor} \label{equivofetaleandflat}
    Let $R$ be a qrsp $\mathbb{Z}_q$-algebras. There is a natural equivalence $$\CB^{\prism_{/(\prism_R,\BI_R)}}_{(G,\mu)}\left((\prism_R,\BI_R)\right)\cong \CB^{\textup{qsyn}}_{(G,\mu)}(R)$$
\end{cor}

\begin{proof}
    Let $R\into R'$ be an \'etale cover. Then, $R'$ is again qrsp (cf. Lemma \ref{etalecoverofqrsp}) and $\prism_{R'}$ is \'etale over $\prism_R$. Since the pair $(\prism_R,\CN^1_{\prism_R})$ is Henselian, the \'etale sites of $R$ and $\prism_R$ are equivalent. Pick an object $\CD$ in $\CB^{\prism_{/(\prism_R,\BI_R)}}_{(G,\mu)}\left((\prism_R,\BI_R)\right)$. By definition, and what we said above, there exists an \'etale cover $R\to S$ such that the base change of $\CD$ along $(\prism_R,\BI_R)\to (\prism_{S},\BI_{S})$ belongs to $\CB_{(G,\mu)}^{\textup{ban}}(\prism_{S},\BI_{S})\into \CB_{(G,\mu)}^{\textup{qsyn}}(S)$. Since the quasi-syntomic topology is finer than the \'etale topology, and $\CB_{(G,\mu)}^{\textup{qsyn}}$ is a stack, we can descend to obtain an object of $\CB_{(G,\mu)}^{\textup{qsyn}}(R)$. This construction shows that we have a fully faithful functor $$\CB^{\prism_{/(\prism_R,\BI_R)}}_{(G,\mu)}\left((\prism_R,\BI_R)\right)\to \CB^{\textup{qsyn}}_{(G,\mu)}(R)$$

    Let us now show the essential surjectivity. Pick an object $\CD\in \CB^{\textup{qsyn}}_{(G,\mu)}(R)$. Similar to the above argument, to give an object in $\CB^{\prism_{/(\prism_R,\BI_R)}}_{(G,\mu)}\left((\prism_R,\BI_R)\right)$ we need an \'etale cover, a banal prismatic higher $(G,\mu)$-display on the cover and an \'etale descent datum. Let $\SQ$ be the $P_{\mu}$-torsor associated with $\CD$. Since $P_{\mu}$ is a representable and smooth group scheme, there exists an \'etale cover $R\to S$, such that $\SQ_S$ is trivial. It follows from Lemma \ref{banality} that the base change $\CD_S$ is banal. Since $\CD$ is defined over $R$, we also have the required \'etale descent datum, and we are done.
\end{proof}

\section{Relation between Witt and Prismatic Displays}\label{SecWindows}

In this section, we want to compare different stacks of displays.  Let us start with the general definition of frames and $(G,\mu)$-windows over frames:

\begin{dfn}\label{DefFrame}
    A \emph{frame} is a quadruple $\SF=(A,I,\sigma,\dot{\sigma})$, where $A$ is a $\BZ_p$-algebra, $I\subset A$ is an ideal, $\sigma:A\to A$ is a ring endomorphism, and $\dot{\sigma}:I\to A$ is a $\sigma$-linear morphism, subject to the following conditions:
    \begin{enumerate}
        \item $\sigma$ is a Frobenius lift, i.e., for all $x\in A$, we have $\sigma(x)\equiv x^p\,\mod\, p $
        \item $I+pA\subset {\textup{Rad}}(A)$
        \item $\dot{\sigma}$ is a $\sigma$-linear epimorphism
    \end{enumerate}

    We say that $\SF$ is a \emph{$\BZ_q$-frame} when $A$ is a $\BZ_q$-algebra.
\end{dfn}

\begin{dfn}\label{DefMapFrame}
    Let $\SF=(A,I,\sigma,\dot{\sigma})$ and $\SG=(B,J,\tau,\dot{\tau})$ be frames. A \emph{morphism of frames} $\kappa:\SF \to \SG$ is a ring homomorphism $\kappa:A\to B$ with $\kappa(I)\subseteq J$ such that $\tau\circ\kappa=\kappa\circ \sigma$ and $\dot{\tau}\circ\kappa=v\cdot \kappa\circ\dot{\sigma}$ for some $v\in B^\times$. If we want to specify $v$, we say that $\kappa$ is a \emph{$v$-morphism}. If $v$ is $1$, then $\kappa$ is called \emph{a strict morphism}.
\end{dfn}

\begin{rmk}
    \begin{itemize}
        \item[1)] Let $\SF=(A,I,\sigma,\dot{\sigma})$ be a frame. One can show that there is a unique $\theta\in A$ such that for all $x\in I$, we have $\sigma(x)=\theta\cdot\dot{\sigma}(x)$ (cf. \cite[Lemma 2.2]{zbMATH05827284}.) We call $\theta$ the \emph{frame constant} of $\SF$.
        
        \item[2)] If $\kappa:\SF\to \SG$ is a frame morphism, then there is a unique unit $v\in B$ for which $\kappa$ is a $v$-morphism. This is because $\kappa\left(\dot{\sigma}(I)\right)$ generates $B$ as a $B$-module. If $\theta$ and $\eta$ are the frame constants of $\SF$ and $\SG$ respectively, then we have $\kappa(\theta)=v\cdot\eta$.
    \end{itemize}
\end{rmk}

\begin{example}
    A prime example of a frame is the \emph{Witt frame}, defined as follows. Let $R$ be a ring. Then $\mathscr{W}(R)=(W(R),I(R),F,V^{-1})$ is a frame, where $F$ is the Frobenius of $W(R)$, $I(R)=V\left(W(R)\right)$ and $V^{-1}$ is the inverse of the Verschiebung $V$. The frame constant of $\SW(R)$ is $p$.
\end{example}

\begin{example}\label{ExNygFramePrism}
    Let $(A,(d))$ be an oriented prism. We can attach two frames to $(A,(d))$. The first being $\SA_\prism\left(A,(d)\right):=\left(A,(d),\sigma, \sigma \frac{1}{d}\right)$ (noting that $d$ is a non-zero divisor). The frame constant of $\SA_\prism$ is $d$. For the second, we assume that $d$ has a $\sigma$-preimage, and define the \emph{Nygaard frame} as $\SA_{\textup{Nyg}}\left(A,(d)\right):=\left(A,\CN^1_A,\sigma,\frac{1}{d}\sigma\right)$. The frame constant of $\SA_{\textup{Nyg}}$ is also $d$.
\end{example}

\begin{dfn}\label{DefBanWind}
    Let $(G,\mu)$ be a display datum over $\BF_q$ and assume that $\mu$ is minuscule. Let $\SF$ be a $\BZ_q$-frame and let $\theta$ be its frame constant. Let $G^\mu(\SF)$ be the subgroup $G(A)\times_{G(A/I)}P_\mu(A/I)$ of $G(A)$. It has an action on $G(A)$ (cf. Definition \ref{frameaction}). The groupoid of \emph{banal $(G,\mu)$-windows over $\SF$} is the quotient groupoid $${\textup{Win}}_{(G,\mu)}^{{\textup{ban}}}(\mathscr{F})=\left[G(A):G^\mu(\SF)\right]$$
\end{dfn}

\begin{rmk}\label{RemPrismDisp=Windows}
    \begin{itemize}
        \item[1)] Let $R$ be a ring. Banal $(G,\mu)$-windows over the Witt frame $\SW(R)$ are the same as banal $(G,\mu)$-displays defined in \cite{bueltel2020g}.
        \item[2)] Let $(A,(d))$ be an oriented prism and assume that $d$ has a $\sigma$-preimage. Then we have an equivalence \[{\textup{Win}}_{(G,\mu)}^{{\textup{ban}}}\left(\mathscr{A}_{\textup{Nyg}}(A,(d))\right)\cong\CB^{\textup{ban}}_{(G,\mu)}(A,(d))\]
    \end{itemize}
\end{rmk}

\begin{rmk}\label{RemWittDisplays}
    Using the functoriality of Witt vectors and Witt descent, we can sheafify the groupoid of banal windows over the Witt frame in faithfully flat topology and obtain the stack ${\textup{Win}}_{(G,\mu)}(\mathscr{W}(R))$. For more details, see \cite{bultel2020displays,lau2021higher}.
\end{rmk}

\begin{dfn}\label{adjointnilpotentdefinition}
    Let $A$ be a ring. Fix an element $\theta\in A$. Let $(G,\mu)$ be a display datum over $\BF_q$. An element $X\in G(A)$ is said to be \emph{adjoint-nilpotent} with respect to $\theta$ if the operator 
        \[\left(\Ad_G(X)\otimes\id\right)\circ \left(\id\otimes \textup{Frob}_p\right) \circ \left(\pi\otimes \id\right)\]
         on $\Lie(U_{\mu^{-1}})\otimes \left(A/(\theta,p)\right)$ is nilpotent, where $\pi:\Lie(G)\to \Lie(U_{\mu^{-1}})$ is the projection killing $\Lie(P_{\mu})$, and the groups $U_{\mu^{-1}}$ and $P_{\mu}$ are respectively the unipotent and parabolic subgroups of $G$ given by $\mu^{-1}$ and $\mu$.  
\end{dfn}

\begin{rmk}\label{RemAdjNilpModRadical}
    \begin{itemize}
        \item[1)] An element $X\in G(A)$ is adjoint-nilpotent with respect to $\theta$ if and only if the operator $$\left(\Ad_G(X)\otimes\id\right)\circ \left(\id\otimes \textup{Frob}_p\right) \circ \left(\pi\otimes \id\right)$$ is nilpotent on $\Lie(U_{\mu^{-1}})\otimes \left(A/\sqrt{(\theta,p)}\right)$. This follows from the fact that a matrix $T\in \BM_n(A)$ ($n\ge 1$) is nilpotent, if and only if the reduction $\bar{T}\in \BM_n(A_{\textup{red}})$ is nilpotent (note that the entries of $T^m$, for any $m\ge 1$, are homogeneous polynomials of degree $m$ in entries of $T$).

        \item[2)] It can be shown, as in \cite[Lemma 4.2]{bultel2020displays}, that the subset of adjoint-nilpotent elements of $G(A)$ is invariant under the $\Phi^\mu$-conjugation (cf. Definition \ref{frameaction}).
    \end{itemize}
\end{rmk}

\begin{dfn}
    \begin{itemize}
        \item[(i)]  Let $\SF=(A,I,\sigma,\dot{\sigma})$ be a frame with frame constant $\theta$. A banal window over $\SF$, represented by an element $X\in G(A)$, is said to be \emph{adjoint-nilpotent} if $X$ is adjoint-nilpotent with respect to $\theta$.  We denote the groupoid of banal adjoint-nilpotent windows over $\mathscr{F}$ by ${\textup{Win}}_{(G,\mu)}^{{\textup{ban}, \circ}}(\mathscr{F})$.

        \item[(ii)] Let $(A,(d))$ be an oriented prism. A banal prismatic higher $(G,\mu)$-display over $(A,(d))$, represented by an element $X\in G(A)$, is said to be \emph{adjoint-nilpotent} if $X$ is adjoint-nilpotent with respect to $d$. The groupoid of banal adjoint-nilpotent prismatic higher $(G,\mu)$-displays over $(A,(d))$ will be denoted by ${\CB}_{(G,\mu)}^{{\textup{ban}, \circ}}(A,(d))$.\\
     \end{itemize}

By Remark \ref{RemAdjNilpModRadical} 2), these definitions are independent of the choice of representatives.
\end{dfn}

\begin{rmk}
    Let $\kappa:\SF=(A,I,\sigma,\dot{\sigma})\to\SG=(B,J,\tau,\dot{\tau})$ be a $v$-morphism of frames. The maps 
    \begin{center}
        \begin{tabular}{l l}
            $G(A)\to G(B)$,  & $G^\mu(\SF)\to  G^\mu(\SG)$\\
            $X\mapsto \kappa(X)\mu(v)$ & $k\mapsto  \kappa(k)$
        \end{tabular}
    \end{center}
    induce a \emph{base-change} functor \[{\textup{Win}}_{(G,\mu)}^{{\textup{ban}}}(\mathscr{F})\to {\textup{Win}}_{(G,\mu)}^{{\textup{ban}}}(\mathscr{\SG})\]

    As nilpotent matrices map to nilpotent ones, the above base-change functor induces a functor:
    \[{\textup{Win}}_{(G,\mu)}^{{\textup{ban},\circ}}(\mathscr{F})\to {\textup{Win}}_{(G,\mu)}^{{\textup{ban},\circ}}(\mathscr{\SG})\] between adjoint-nilpotent windows.
\end{rmk}

\begin{rmk}
    If $(A,(d))\to (B,(d))$ is a morphism of oriented prisms, then the base-change functor (cf. Lemma \ref{BanalG-muFunc}) induces a functor 
    \[{\CB}_{(G,\mu)}^{{\textup{ban},\circ}}(A,(d))\to {\CB}_{(G,\mu)}^{{\textup{ban},\circ}}(B,(d))\]
\end{rmk}

\begin{lem}\label{decentofadjointnipotent}
    Let $\SF=(A,I,\sigma,\dot{\sigma})$ and $\SG=(B,J,\tau,\dot{\tau})$ be frames, with frame constants $\theta$ and $\eta$, respectively. Let $\kappa:\SF \to \SG$ a morphism such that the underlying map $A\to B$ is faithfully flat. Then a banal window $\CD$ over $\mathscr{F}$ is adjoint-nilpotent with respect to $\theta$ if and only if the base change $\CD_{\mathscr{G}}$ is adjoint-nilpotent with respect to $\eta$. Similarly, if $(A,(d))\to (B,(d))$ is a covering of oriented prisms, then a banal prismatic higher $(G,\mu)$-display is adjoint-nilpotent over $(A,(d))$, if and only if it is so after base change to $(B,(d))$.
\end{lem}

\begin{proof}
    Let us explain the case of windows. The other case is similar. The base change of a nilpotent operator is again nilpotent; for the inverse, note that $A\to B$ being faithfully flat, the map $\frac{A}{(\theta,p)}\to \frac{B}{(\eta,p)}$ is injective and therefore for any $n\ge 1$, the ring homomorphism $\BM_n\left(\frac{A}{(\theta,p)}\right)\to \BM_n\left(\frac{B}{(\eta,p)}\right)$ is injective.
\end{proof}

\begin{example}
    Let $A$ be a flat $\BZ_p$-algebra, $\sigma$ a Frobenius lift on $A$ and $\Fa$ a pd-ideal. Since for all $x\in \Fa$, $p!\gamma_p(x)=x^p\equiv \sigma(x) \mod p$, we have $\sigma(\Fa)\subset pA$, and so we have a well-defined $\sigma$-linear map $\dot{\sigma}=\frac{1}{p}\sigma:\Fa\to A$. Assume that $pA\subset {\textup{Rad}}(A)$. We obtain a frame $\frame{A}=(A,\Fa,\sigma,\dot{\sigma})$ (note that the assumptions imply that $\Fa+pA\subset {\textup{Rad}}(A)$).
\end{example}

\begin{pro}\label{framemorphismpdtowitt}
    Let $A$ be a flat $\BZ_p$-algebra, $\sigma$ a Frobenius lift on $A$ and $\Fa$ a pd-ideal. Assume that $pA\subset {\textup{Rad}}(A)$. The projection $W(A)\onto W(A/\Fa)$ induces a frame morphism $\frame{A}\to \mathscr{W}(A/\Fa)$.
\end{pro}

\begin{proof}
    Since $A$ is $p$-torsion-free, $W(A)$ is also $p$-torsion-free. We introduce an auxiliary frame $\SW:=(W(A),W(\Fa)+I(A),F, \dot{F})$, where $W(\Fa)$ is the kernel of the projection $W(A)\onto W(A/\Fa)$, and $\dot{F}$ is $\frac{1}{p}F$ as usual. Note that since $\Fa$ is a pd-ideal, $W(\Fa)$ is also a pd-ideal, so $\frac{1}{p}F$ is defined on it. Indeed, denoting by $\gamma_i$ the pd-structure on $\Fa$, we have Zink's logarithmic ghost maps (\cite[\S1.4]{zink2002display})       \begin{align*}
          \textbf{w}'_n:W(\Fa)\to&\, \Fa\\
          (a_0,a_1,\dots)\mapsto& \sum_{i=0}^n(p^i-1)!\gamma_{p^i}(a_{n-i})
        \end{align*}
 
   One can see that these induce a $W(A)$-linear bijection $\log:W(\Fa)\to \Fa^{\BN}$, where for $\underline{x}\in W(A)$ and $\underline{a}=[a_0,a_1,\dots]\in \Fa^{\BN}$, we set $\underline{x}\cdot \underline{a}=[\textbf{w}_0(\underline{x})a_0,\textbf{w}_1(\underline{x})a_1,\dots]$ (the $\textbf{w}_i$ being the classical ghost maps). The subset $\log^{-1}[\Fa,0,0,\dots]\subset W(A)$ is an ideal with a pd-structure inherited from that of $\Fa$. Let us keep the notation $\Fa$ for this ideal. We have $W(\Fa)\subset\Fa\oplus I(A)=\kernel \left(W(A)\onto A\onto A/\Fa\right)$, and since $\Fa\oplus I(A)$ is a pd-ideal, $W(\Fa)$ is also a pd-ideal. Note that $W(\Fa)+I(A)$ is the kernel of $W(A)\onto A/\Fa$.
   
   Since the Cartier diagonal $\Delta:A\to W(A)$ lifts the identity of $A$, it induces a strict frame morphism $\frame{A}\to \SW$. Now, we show that the projection $\pi:W(A)\to W(A/\Fa)$ induces a strict frame morphism $\SW\to \SW(A/\Fa)$. We only need to check the compatibility of $\dot{F}$ and $V_{A/\Fa}^{-1}$ (we use the subscript to emphasize the ring). Since $\dot{F}$ and $V_A^{-1}$ agree on $I(A)$ and the Witt frame is functorial in ring morphisms, we only need to check that for all $w\in W(\Fa)$, we have $\pi(\dot{F}(w))=V_{A/\Fa}^{-1}(\pi(w))=V_{A/\Fa}^{-1}(0)=0$. It is enough to show that $\dot{F}(w)\in W(\Fa)$. Since $W(A)$ has a $\delta$-structure, we have $F(w)=w^p+p\delta(w)=p!\gamma_p(w)+p\delta(w)$ (using the notation $\gamma$ for the pd-structure on $W(\Fa)$), and so $\dot{F}(w)=(p-1)!\gamma_p(w)+\delta(w)$. Since $\pi$ preserves the $\delta$-structures, $\delta(w)\in W(\Fa)$, and the proof is achieved.
\end{proof}

\begin{rmk}\label{Remframemorphismpdtowitt}
    Let $\frame{A}$ be the frame in the example, and $v\in A^\times$. Let us denote by $\frame{A}_v$ the frame $(A,\Fa, \sigma,v\cdot\dot{\sigma})$. The identity map of $A$ then induces a $v^{-1}$-morphism of the frame $\frame{A}_v\to\frame{A}$. Setting $v':=\pi\left(\Delta(v^{-1})\right)\in W(A/\Fa)$, by the proposition we have a $v'$-morphism $\frame{A}_v\to \SW(A/\Fa)$.
\end{rmk}

\begin{rmk}\label{RemPrisDisptoDisp}
    We want to compare prismatic and Witt displays. Let $(A,(d))$ be an oriented prism and set $R:={A}/{\CN^1_A}$. Assume that $d$ has a $\sigma$-preimage. Since $A$ is a $\delta$-ring, the projection $A\onto R$ lifts to a map $A\to W(R)$ (Proposition \ref{DeltaWittAdj}). As in the proof of \cite[proposition 5.21]{anschutz2023prismatic} under this map $\CN^1_A$ goes to $I(R)$, and one can see that $d$ maps to $pu$ for a unit $v$. In order to show that $\mathscr{A}_{\textup{Nyg}}\left(A,(d)\right)\to \mathscr{W}(R)$ is a $v$-morphism, we need to show that $\frac{1}{d}\sigma$ and $V^{-1}$ are compatible. This can be done in some special cases. When $W(R)$ is $p$-torsion-free (for example, when $R$ is a flat $\mathbb{Z}_p$-algebra or when $R$ is a reduced $\mathbb{F}_p$-algebra) we can check the compatibility after multiplying by $d$ which follows from the compatibility of $\sigma$ and $F$. As we explained above when $\mathscr{A}_{\textup{Nyg}}\left(A,(d)\right)\to \mathscr{W}(R)$ is a frame morphism, we get a morphism from the groupoid of banal prismatic higher $(G,\mu)$-displays over $(A,(d))$ to the groupoid of banal $(G,\mu)$-displays over $R$.
\end{rmk}

\begin{pro}\label{ReductionModpPerfVer}
    Let $\Spa(R,R^+)\to\Spd\BZ_q$ be an element of $\Perf/\Spd\BZ_q$, corresponding to an untilt $\Spa(R^\sharp,R^{\sharp+})$ over $\BZ_q$. Assume that $R^{\sharp+}$ is $p$-torsion-free. Then, there is a natural \emph{``reduction mod $p$''} functor
    \begin{align*}
        \mathcal{B}^{\textup{perf}}_{(G,\mu)}(R^\sharp,R^{\sharp+})\to &\, {\textup{Win}}_{(G,\mu)}(\mathscr{W}(R^{\sharp+}/p))\\
        \CD\mapsto &\, \bar{\CD}
    \end{align*}
\end{pro}

\begin{proof}
    As we have seen, there is a frame morphism from $\mathscr{A}_{\textup{Nyg}}(A_{\textup{inf}}(R^{\sharp+}),\kernel \theta^+)$ to $\mathscr{W}(R^{\sharp+})$, and so there is a frame morphism from $\mathscr{A}_{\textup{Nyg}}(A_{\textup{inf}}(R^{\sharp+}),\kernel \theta^+)$ to $\mathscr{W}(R^{\sharp+}/p)$. Under this map, we get a functor between banal displays. The reduction mod $p$ defines a map of sites between the \'etale site on integral perfectoids over $R$ and the \'etale site over schemes over $R^\sharp/p$ (in fact, it is an equivalence of categories by the Almost Purity Theorem \cite[Theorem 7.4.5]{scholze2020berkeley}). So, after sheafification, we get the desired functor. 
\end{proof}

Similarly, we have:

\begin{pro}\label{ReductionModp}
    Assume that $R$ is a $p$-torsion-free quasi-syntomic ring. Then there is a natural \emph{``reduction mod $p$''} functor $$\mathcal{B}^{\textup{qsyn}}_{(G,\mu)}(R)\to {\textup{Win}}_{(G,\mu)}(\mathscr{W}(R/p))$$
    
    We denote the image of $\CD$ under this functor by $\bar{\CD}$.
\end{pro}

\begin{proof}
    When $R$ is qrsp, then we have an explicit description of elements of the groupoid $\mathcal{B}^{\textup{qsyn}}_{(G,\mu)}(R)$, and so, the same arguments as in the proof of the previous proposition give the reduction map. Now, since both $R\mapsto \mathcal{B}^{\textup{qsyn}}_{(G,\mu)}(R)$ and $R\mapsto {\textup{Win}}_{(G,\mu)}(\mathscr{W}(R/p))$ are quasi-syntomic stacks, the reduction map extends to the case of general quasi-syntomic rings.
\end{proof}

In \cite{ito2023prismatic} Ito defines a notion of prismatic displays. His definition is different from ours, as we now explain. We will use the notation and setting of Section \ref{SectionBanalObjects}. Ito works with the filtration $$A\supseteq (d)\supseteq (d^2)\supseteq \cdots $$ instead of the Nygaard filtration. More precisely, let the group $G^{\mu,\min}(A,d)$ be the set of homomorphisms $g: S_
{\BZ_q}\to A$ such that if $s\in \mathcal{F}^i$ then $g(s)\in (d^i)$ and consider the quotient groupoid $$\left[G(A):G^{\mu,\min}(A,d)\right]$$ where $G^{\mu,\min}(A,d)$ acts on $G(A)$ by the formula:
$$k\cdot X=k^{-1}X\sigma\left(\mu(d)k\mu(d)^{-1}\right)$$

\begin{dfn}
    A \emph{minuscule Breuil-Kisin-Fargues module} over the prism $(A,I)$ is a Breuil-Kisin-Fargues module $(M,F)$ over $(A,I)$ such that the linear module generated by $F(\sigma^*M)$ contains $IM$.
\end{dfn}

Ito's prismatic displays are related to minuscule Breuil-Kisin-Fargues modules.  The map $\sigma:G(A)\to G(A)$ induces a map of pairs $$(G(A),G^{\mu}(A,d))\to (G(A),G^{\mu,min}(A,d))$$ which then induces a functor:
$$T:\left[G(A):G^{\mu}(A,d)\right]\to\left[G(A):G^{\mu,min}(A,d)\right]$$

This functor is compatible with the functors to the groupoid of Breuil-Kisin-Fargues modules on both sides, which are fully faithful. Therefore, $T$ is also fully faithful. It is essentially surjective when $\sigma$ is surjective, for example, when $A$ is a perfect prism. 

\section{Deformation Theory of Prismatic Displays}\label{deformationtheory}

In this section, we want to prove a Grothendieck-Messing deformation theorem for prismatic higher $(G,\mu)$-displays. 

First, let us collect some facts about changing frames.

\begin{dfn}
    A morphism of frames is said to be \emph{crystalline} (respectively, \emph{nil-crystalline}) if, for all display data $(G,\mu)$, it induces an equivalence on the associated groupoids of banal $(G,\mu)$-windows (respectively, banal adjoint-nilpotent $(G,\mu)$-windows).
\end{dfn}

The basic tool in the deformation theory of $(G,\mu)$-windows is the Unique Lifting Lemma originally proved by Zink for $G=\GL_n$ and later generalized by others in various settings; for important examples, see \cite{lau2021higher, bultel2020displays, bartling2022mathcal}. We will use the version proved in \cite{bartling2022mathcal}:

\begin{pro}[Unique Lifting Lemma]\label{uniqueliftinglemma}
    Let $\kappa:\mathscr{F}=(A,I,\sigma,\dot{\sigma})\to\mathscr{G}=(B,J,\tau,\dot{\tau})$ be a morphism of frames, and consider the following conditions:
        \begin{enumerate}
          \item $\kappa:A\to B$ is a surjection inducing an isomorphism $A/I\cong B/J$.
          \item $(A,\kernel \kappa)$ is Henselian.
          \item $\kernel \kappa$ is $p$-adically complete.
          \item $\dot{\sigma}$ is pointwise topologically nilpotent on $\ker \kappa$.
        \end{enumerate}

    If $\kappa$ satisfies conditions 1, 2, and 3, then it is \emph{nil-crystalline}. Furthermore, if it also satisfies condition 4, then it is \emph{crystalline}.
\end{pro}

\begin{proof}
    Essential surjectivity follows from \ref{FormalSmoothLift}. Fully-faithfulness is essentially \cite[Proposition 4.11]{bartling2022mathcal}. There, it is assumed that the morphism is strict. But this is not necessary because if $\kappa$ is a $v$-morphism, we can consider the frame $\mathscr{G}_\mu=(B,J,\tau,v^{-1}.\dot{\tau})$. Then $id:B\to B$ induces a (nil-)crystalline morphism $v:\mathscr{G}\to \mathscr{G}_\mu$, and the composition $\mathscr{F}\to \mathscr{G}_\mu$ is a strict morphism that still satisfies the above properties. Therefore, it is (nil-)crystalline, and so $\kappa$ is (nil-)crystalline as desired.
\end{proof}

Let $\left\{\SF_n=(A_n,I_n,\sigma_n,\dot{\sigma}_n)\right\}_n$ be an inverse system of frames, with strict transition morphisms and frame constant $\theta_n$. We obtain a frame $\varprojlim \SF_n$ by taking the limit over all components of the frames. In the other direction, let $\SF=(A,I,\sigma,\dot{\sigma})$ be a frame and $J\subset A$ an ideal, such that $\sigma(J)\subset J$ and for all $n\geq 1$, $\dot{\sigma}(I\cap J^n)\subset J^n$. Assume that $A$ and $A/I$ are $J$-adically complete. Then we have an inverse system of frames $\{\SF_n:=(A/J^n,I/(I\cap J^n),\sigma,\dot{\sigma})\}$ and there is a canonical isomorphism $\SF\cong \varprojlim\SF_n$ in the sense of the following lemma:

\begin{lem}\label{completionofframes}
    Let $(G,\mu)$ be a display datum, and $\SF=(A,I,\sigma,\dot{\sigma})$ a frame, with frame constant $\theta$. Let $J\subset A$ be an ideal such that $\sigma(J)\subset J$ and for all $n\geq 1$, $\dot{\sigma}(I\cap J^n)\subset J^n$. Assume that $A$ and $A/I$ are $J$-adically complete.  The strict morphisms $\SF\to \SF_n$ give rise to an isomorphism 
    \[{\textup{Win}}^{\textup{ban}}_{(G,\mu)}(\SF)\cong 2\text{\small \textendash}\varprojlim{\textup{Win}}^{\textup{ban}}_{(G,\mu)}(\SF_n)\] 
    
    If furthermore, $A/(\theta, p)$ is also $J$-adically complete, the above isomorphism restricts to an isomorphism \[{\textup{Win}}^{\textup{ban},\circ}_{(G,\mu)}(\SF)\cong 2\text{\small \textendash}\varprojlim{\textup{Win}}^{\textup{ban},\circ}_{(G,\mu)}(\SF_n)\]
\end{lem}

\begin{proof}
    The group schemes $G$ and $P_\mu$ are affine and smooth over $\mathbb{Z}_p$. By assumption $A$ and $A/I$ are $J$-adically complete. It follows from Lemma \ref{FormalSmoothLift} that $G(A)\cong\varprojlim G(A/J^n)$ and so $G^\mu(\mathscr{F})\cong\varprojlim G^\mu(\mathscr{F}_n)$. The first isomorphism now follows. The second isomorphism is a direct consequence of the first one.
\end{proof}

\begin{rmk}
    One natural question is that if $R=\varprojlim R_i$, and all $R_i$ are quasi-syntomic rings, under what conditions do we have $$\mathcal{B}^{\textup{qsyn}}_{(G,\mu)}(R)\cong 2\text{\small \textendash}\varprojlim \mathcal{B}^{\textup{qsyn}}_{(G,\mu)}(R_i)$$ 
    
    This question is harder than the case of Witt-displays because if $R'\to R$ is a thickening of qrsp rings, $\prism_{R'}\to \prism_R$ is not surjective, unlike $W(R')\to W(R)$.
\end{rmk}

\begin{dfn}\label{DefNilCat}
    Let $R$ be a quasi-syntomic ring. We denote by $\mathfrak{Nil}_R^1$ the category whose objects are \emph{first-order quasi-syntomic thickenings} of $R$, i.e., pairs $(R',f)$, where $R'$ is a quasi-syntomic ring and $f:R'\onto R$ is a morphism such that $(\kernel f)^2=0$. A morphism from $(R_1',f_1)$ to $(R'_2,f_2)$ is a morphism $g:R_1'\to R_2'$ such that $f_2\circ g=f_1$.
\end{dfn}

\begin{dfn}\label{DefDeformationFunctor}
    Let $S$ be a quasi-syntomic ring and $\CD$ be a prismatic higher $(G,\mu)$-display over $S$. The \emph{deformation} functor attached to $\CD$ is the covariant functor
    \begin{align*}
        \Def_{\CD}:\mathfrak{Nil}^1_S\to& \Set\\
        (S',f)\mapsto& \{(\CD',\iota)\}/\sim
    \end{align*}
    where $\CD'$ is a prismatic higher $(G,\mu)$-display over $S'$ and $\iota:\CD'_{f,S}\cong \CD$ is an isomorphism between the base change of $\CD'$ along $f$ and $\CD$.
\end{dfn}

Similarly, we make the following definition.

\begin{dfn}\label{DefDeformationFunctorHodge}
    Let $R$ be a ring and $\SH=(\SQ,\SM,\tau)$ a $(G,\mu)$-Hodge filtration on $R$. The deformation functor attached to $\SH$ is the functor
    \begin{align*}
        \Def_{\SH}:\mathfrak{Nil}^1_R\to& \Set\\
        (R',f)\mapsto& \{(\SH',\iota)\}/\sim
    \end{align*}
    where $\SH'$ is a $(G,\mu)$-Hodge filtration over $R'$ and $\iota:\SH'_{f,R}\cong \SH$ is an isomorphism between the base change of $\SH'$ along $f$ and $\SH$.
\end{dfn}

\begin{rmk}\label{RemHodgeFiltOnDefs}
    Let $S$ be a quasi-syntomic ring and $\CD\in\mathcal{B}^{\textup{qsyn}}_{(G,\mu)}(S)$. By Remark \ref{RemStackofHodgeFilt} we have a $(G,\mu)$-Hodge filtration $\FH(\CD)\in\SH_{(G,\mu)}(S)$. By functoriality, in fact, we have a natural transformation $\FH:\Def_{\CD}\to\Def_{\FH(\CD)}$. In this subsection, we will prove that, at least for adjoint-nilpotent displays, under some mild conditions on $S$, this transformation is an isomorphism.
\end{rmk}

Let $R'$ be a flat integral perfectoid $R_0$-algebra and $S'\onto S$ a first-order thickening of fgpp R-algebras. To simplify the notations, set $B:=\prism_{S}$ and $B':=\prism_{S'}$.

Set $\SB=\SA_{\textup{Nyg}}\left(B,([p]_q)\right)$ and $\SB'=\SA_{\textup{Nyg}}(B',([p]_q))$. Since $\sigma([p]_{q^{1/p}})=[p]_{q}$, we have identifications \begin{equation}\label{EqPrismDisp=Windows}
    \CB_{(G,\mu)}^{\textup{ban},\circ}(B,([p]_q))={\textup{Win}}^{\textup{ban},\circ}_{(G,\mu)}(\SB)\quad \text{and}\quad \CB_{(G,\mu)}^{\textup{ban},\circ}(B',([p]_q))={\textup{Win}}^{\textup{ban},\circ}_{(G,\mu)}(\SB')
\end{equation}
(see Remark \ref{RemPrismDisp=Windows} 2). Therefore, we are reduced to studying the deformations of $(G,\mu)$-windows over these frames. The strategy for proving Grothendieck-Messing for Witt displays on the thickening $S'\onto S$ was to put a new frame structure on $W\left(S'\right)$ such that the map between this frame and the Witt frame $W(S)$ is crystalline and then compare two frame structures on $W(S')$. In the prismatic setting, the map $B'\to B$ is injective and not surjective.
 
We want to put a new frame structure on $B$ such that the map $B'\to B$ becomes nil-crystalline for the new frame structure. Recall from Lemma \ref{universalthickening} that there is a surjection $B\onto S'$ with kernel $J$, satisfying $J\cap B'=\CN^1_{B'}$, and inducing an isomorphism $B'/\CN^1_{B'}\cong B/J$. Since $[p]_{q^{1/p}}\in \CN^1_{B'}\subset J\subset \CN^1_B$, we have a frame $\mathscr{B}/\mathscr{B'}=(B,J,\sigma,\dot{\sigma}:=\frac{1}{[p]_q}\sigma)$. There is a natural frame morphism $\kappa:\SB'\to \SB/\SB'$. 

\begin{pro}\label{prismaticuniquelifitng}
    The frame morphism $\kappa$ is nil-crystalline.
\end{pro}

\begin{proof}
    For all $n\geq 1$, set $B_n:=B/[w]^n$ and $B'_n:=B'/[w]^n$, and let $J_n$ and $\CN^1_{B',n}$ be the images of $J$ and $\CN^1_{B'}$ in $B_n$ and $B'_n$, respectively. We have $\sigma([w])=[w]^p$, $\dot{\sigma}\left(([w]^n)\cap \CN^1_{B'}\right)\subset [w]^nB'$ and $\dot{\sigma}\left(([w]^n)\cap J\right)\subset [w]^nB$, because $[w]^n,[p]_q$ is a regular sequence (Remark \ref{pcompleteness}). This implies that $\mathscr{B}'_n=:(B'_n,\CN^1_{B',n},\sigma,\dot{\sigma})$ and $\mathscr{B}_n/\mathscr{B'}_n:=(B_n,J_n,\sigma,\dot{\sigma})$ are frames.  Let $\kappa_n:\mathscr{B}'_n\to \mathscr{B}_n/\mathscr{B'}_n$ be the reduction of $\kappa$ mod $[w]^n$.

    By Remark \ref{pcompleteness}, the rings $B$ and $B'$ and the reductions $B/JB$ and $B'/\CN^1_{B'}$ (both isomorphic to $S'$) are $[w]$-adically complete. Since elements $[p]_q$ and $[w]^p$ are equal mod $p$, and the rings $B/p$ and $B'/p$ are $[w]$-adically complete, the quotient rings $B/([p]_q,p)$ and $B'/([p]_q,p)$ are $[w]$-adically complete as well. So, if we show that for all $n\geq 1$, $\kappa_n$ is nil-crystalline, we can apply Lemma \ref{completionofframes}, and conclude that $\kappa$ is nil-crystalline (note that $[p]_q$ is the frame constant for the frames $\SB'$ and $\SB/\SB'$).  We will prove by induction on $n$ that $\kappa_n$ is nil-crystalline.

    We first show that $\kappa_1:\mathscr{B}'_1\to \mathscr{B}_1/\mathscr{B'}_1$ is nil-crystalline. Since $\delta([w])=0$, the rings $B_1$ and $B'_1$ are $\delta$-rings. Consider the projections $\bar{h}:B_1\onto S'/p$ and $\bar{h}':B'_1\onto S'/p$. By adjunction between $\delta$-rings and Witt vectors (Proposition \ref{DeltaWittAdj}) there are canonical lifts $h:B_1\to W(S'/p)$ and $h':B_1'\to W(S'/p)$. By Lemma \ref{Nygaardofperfectlypresented}, $\bar{h}$ and $\bar{h}'$ are pd-thickenings. By Proposition \ref{framemorphismpdtowitt} and Remark \ref{Remframemorphismpdtowitt}, we have two frame morphisms $\Fh:\SB_1/\SB'_1\to \SW(S'/p)$ and $\Fh':\SB'_1\to \SW(S'/p)$. By the universal property of $A_{\textup{crys}}$ and functionality of the Cartier map, the surjective map $A_{\textup{crys}}(S'/p)\onto W(S'/p)$ factors through $h$ and also through $h'$. Therefore, $h$ and $h'$ are surjective. Furthermore, the kernels of $h$ and $h'$ are Henselian, because they are subsets of $\CN^1_BB_1$ and $\CN^1_{B'}B'_1$ respectively (use Proposition \ref{Henselianpair}, Lemmas \ref{LemHenselQuot} and \ref{LemHenselSubId}). Since $B_1$, $B'_1$ and $W(S'/p)$ are $p$-adically complete (Remark \ref{pcompleteness}, and \cite[Proposition 3]{zink2002display}) the kernels are also $p$-adically complete. By Unique Lifting Lemma \ref{uniqueliftinglemma} these maps are nil-crystalline, and so $\kappa_1$ is also nil-crystalline.
    
    Now, in order to show the induction step, we prove that for each $n$, the induced map between the fibers of the vertical maps in the diagram
    \[ \begin{tikzcd}
            {\textup{Win}}_{G,\mu}^{{\textup{ban,$\circ$}}}(\mathscr{B}'_{n+1})\arrow[r]\arrow[d]&{\textup{Win}}_{G,\mu}^{{\textup{ban,$\circ$}}}(\mathscr{B}_{n+1}/\mathscr{B'}_{n+1})\arrow[d]\\
             {\textup{Win}}_{G,\mu}^{{\textup{ban,$\circ$}}}(\mathscr{B}'_n)\arrow[r]&{\textup{Win}}_{G,\mu}^{{\textup{ban,$\circ$}}}(
            \mathscr{B}_n/\mathscr{B'}_n)
        \end{tikzcd} \]
     are equivalences.

     The isomorphism $B'/\CN^1_{B'}\cong B/J$ induces, for all $n\geq 1$, compatible isomorphisms $$B'/\left(\CN^1_{B'}+[w]^nB'\right)\cong B/\left(J+[w]^nB\right)$$
     
     We denote this quotient ring by $R_n$. Let $R_{n+1}\onto R_{n}$ be the projection. We claim that the fibers (up to isomorphism) of the vertical maps in the above diagram are canonically in bijection with the set \[\left(G(R_{n+1})\times_{G(R_n)}P_{\mu}(R_n)\right)/P_{\mu}(R_{n+1})\]

     Again, we show this using the Unique Lifting Lemma. Let $(\Lambda, \Fa)$ be either $(B,J)$ or $(B',\CN^1_{B'})$. By the above discussion, we have $\Lambda/\left(\Fa+[w]^n\Lambda\right)\cong R_n$ and there are canonical isomorphisms 
     \begin{equation}\label{EqLambda}
         R_n\cong \Lambda_{n+1}/(\Fa_{n+1}+[w]^n\Lambda_{n+1})\quad \text{and} \quad R_{n+1}\cong \Lambda_{n+1}/\Fa_{n+1}
     \end{equation}
     
     Consider the frame $$(\Lambda_{n+1},\Fa_{n+1}+[w^n]\Lambda_{n+1},\sigma,\dot{\sigma})$$ where we define $\dot{\sigma}([w^n])=\sigma([w^n])=0$. Then the morphism  $$(\Lambda_{n+1}, \Fa_{n+1}+[w^n]\Lambda_{n+1})\to (\Lambda_{n}, \Fa_{n})$$ satisfies the conditions of the Unique Lifting Lemma and hence is nil-crystalline. Therefore, it remains to determine the fibers of $${\textup{Win}}_{G,\mu}^{{\textup{ban},\circ}}(\Lambda_{n+1},\Fa_{n+1})\to {\textup{Win}}_{G,\mu}^{{\textup{ban},\circ}}(\Lambda_{n+1},\Fa_{n+1}+[w]^n\Lambda_{n+1})$$
     
     These fibers are in bijection with the following set (using identities (\ref{EqLambda}))
    \begin{multline*}
        \Gamma^\mu(\Lambda_{n+1},\Fa_{n+1}+[w^n]\Lambda_{n+1})/\Gamma^\mu(\Lambda_{n+1},\Fa_{n+1}) = \\\left(G(\Lambda_{n+1})\times_{G\left(R_{n}\right)}P_{\mu}\left(R_n\right)\right)/\left(G(\Lambda_{n+1})\times_{G\left(R_{n+1}\right)}P_{\mu}\left(R_{n+1}\right)\right)
    \end{multline*}

      By Lemma \ref{DefFib} this quotient is in bijection with the set \[\left(G(R_{n+1})\times_{G(R_n)}P_{\mu^{-1}}(R_n)\right)/P_{\mu}(R_{n+1})\] as desired.
\end{proof}

\begin{rmk}
    Let $f:S'\onto S$ be a nilpotent thickening of quasi-syntomic rings. Then, any lift of an adjoint-nilpotent prismatic display over $S$ to $S'$ is again adjoint-nilpotent: Choose $X\in G(\prism_{S'})$. Since the kernel of $f$ is nilpotent, it induces an isomorphism on the reduced quotient rings \[\prism_{S'}/\sqrt{(p,[p]_q)}\cong S'/\sqrt{(p)}\xrightarrow{\sim} S/\sqrt{(p)}\cong \prism_S/\sqrt{(p,[p]_q)}\] and so  modulo $\sqrt{(p,[p]_q)}$, elements $X$ and $f(X)$ are equal. Therefore, by Remark \ref{RemAdjNilpModRadical} 1), the display represented by $f(X)$ is adjoint-nilpotent if and only if the display represented by $X$ is adjoint-nilpotent.
\end{rmk}

\begin{pro}\label{BanalGrothendiquemessing}
    Let $R$ be a flat integral perfectoid $R_0$-algebra and $S'\onto S$ a first-order thickening of fgpp $R$-algebras. There is a natural bijection between the isomorphism classes of deformations of an adjoint-nilpotent banal prismatic higher $(G,\mu)$-display along $S'\onto S$ and the set $$\left(G(S')\times_{G(S)}P_{\mu}(S)\right)/P_{\mu}(S')$$ 
\end{pro}

\begin{proof}
    By Proposition \ref{prismaticuniquelifitng} and identities (\ref{EqPrismDisp=Windows}) it is enough to determine the fibers of $${\textup{Win}}_{G,\mu}^{{\textup{ban}}}(\mathscr{B/B'})\to {\textup{Win}}_{G,\mu}^{{\textup{ban}}}(\mathscr{B})$$ 
    
    Similar arguments as before show that these are in bijection with $$\left(G(S')\times_{G(S)}P_{\mu}(S)\right)/P_{\mu}(S')$$
\end{proof}

Now we want to prove this for a general thickening of quasi-syntomic rings that satisfy $(\leftmoon)$. First, we need a lemma:

\begin{lem}
    Any deformation of a banal prismatic higher $(G,\mu)$-displays along nilpotent thickenings is banal.
\end{lem}

\begin{proof}
    Using Proposition \ref{banality}, we need to show that the $P_\mu$-torsor corresponding to the deformation of a banal $(G,\mu)$-display is trivial. This $P_\mu$-torsor deforms the $P_\mu$-torsor of the original prismatic higher $(G,\mu)$-display, which is trivial by Proposition \ref{banality}. But since $P_\mu$ is a smooth group scheme, any deformation of a trivial $P_\mu$-torsor is trivial. Indeed, if $\SQ$ is a $P_\mu$-torsor over a scheme $X$, $\underline{\Isom}(\SQ,P_\mu\times X)$ is a $P_\mu$-torsor, and therefore is smooth (cf. \cite[Ch. III, Prop. 4.2]{7c7cb5a6-e464-3506-b185-42a30a3f8fde}), and so, can lift its $S$-point to a $S'$-point.
\end{proof}

\begin{thm}\label{grothendique-messing}
    Let $S$ be a quasi-syntomic ring and $\CD$ an adjoint-nilpotent prismatic higher $(G,\mu)$-display on $S$. Let $(S',f)$ be an object of $\mathfrak{Nil}_S^1$ and assume that $S$ and $S'$ satisfy condition $(\leftmoon)$. The map $\FH:\Def_{\CD}(S',f)\to\Def_{\FH(\CD)}(S',f)$ is a bijection.
\end{thm}

\begin{proof}
    Both categories $\mathcal{B}^{\textup{qsyn}}_{(G,\mu)}$ and $\SH_{(G,\mu)}$ are stacks, so one can prove this after going to a cover: after going to some \'etale cover, we can assume that the display on $S$ becomes trivial (one can lift the thickening to this cover because of the invariance of \'etale site on nilpotent thickening). Then, using Andr\'e's lemma (\cite[Theorem 7.14]{10.4007/annals.2022.196.3.5}) we can find a suitable quasi-syntomic cover $\tilde{S}'$ of $S'$ such that both $\tilde{S}'$ and $S\otimes_{S'}\tilde{S'}$ are fgpp algebras over $R_0$. Over this cover, we have the bijection between lifts of a banal display and lifts of the associated Hodge filtration. By descent, we obtain the desired bijection over $S'$.
\end{proof}

One can also think about some thickenings where $S$ is quasi-syntomic and $S'$ is not. Let us restrict ourselves to the following special class of thickenings, although our techniques should work in greater generality.

\begin{dfn}\label{DefAdmisThick}
    Let $R$ be a flat integral perfectoid $R_0$-algebra. We say that a thickening $\tilde{S}\onto S$ is \emph{admissible} if $S=R/I$ is an fgpp algebra over $R$ and $\tilde{S}=R/I^2$.
\end{dfn}

\begin{example}\label{universalexample}
    The ``universal'' example of an admissible thickening is $\tilde{S}(n)\onto S(n)$ where $$\tilde{S}(n)=R_0\langle X_1^{1/p^\infty},\cdots,X_n^{1/p^\infty}\rangle/(X_1,\cdots,X_n)^2$$
    
    We can define a natural prism over $\tilde{S}(n)$: Let $\prism_{\tilde{S}(n)}$ be the $(p,[p]_q)$-completion of the sub-$\delta$-ring of $\prism_{S(n)}$, generated by $A_0$, $X_i^{1/p^\infty}$  and $\frac{X_iX_j}{[p]_q}$ for $1\le i,j \le n$. 

    By \cite[Lemma 12.5]{10.4007/annals.2022.196.3.5}, this ring actually contains $\frac{(Y_iY_j)^t}{[\lfloor t\rfloor]_q!}$ for all $t,1\le i,j\le n$. Therefore, $\prism_{\tilde{S}(n)}$ is the $(p,[q]_p)$-adic completion of the ring generated by elements $\frac{(Y_iY_j)^t}{[\lfloor t\rfloor]_q!}$ over $A_0$. By this explicit description and the fact that multiplication by $[p]_q$ preserves the natural grading, one can easily check that $\CN^1_{\tilde{S}}$ is generated by elements $$[p]_{q^{1/p}} \quad \text{and} \quad \frac{(Y_iY_j)^t}{[\lfloor t\rfloor]_q!}\quad (t\in \mathbb{N},1\le i,j\le n)$$ and we have $\prism_{\tilde{S}(n)}/\CN^1_{\tilde{S}(n)}=\tilde{S}(n)$.
\end{example}

\begin{dfn}
    Let $S$ be a $R_0$-algebra. The \emph{flat absolute prismatic site} of $S$ is the site whose underlying category is the opposite category of those prisms over $S$ that are flat over $A_0$. We denote this site by $(S)_{\prism}^{\textup{fl}}$.    
\end{dfn}

\begin{pro}\label{initialflat}
    The prism $(\prism_{\tilde{S}(n)},([p]_q))$ is the final object of the site $(\tilde{S}(n))_{\prism}^{\textup{fl}}$
\end{pro}

\begin{proof}
    The initial prism over $R_n=R_0\langle X_1^{1/p^\infty},\cdots,X_n^{1/p^\infty}\rangle$ is $A_0\langle X_1^{1/p^\infty},\cdots,X_n^{1/p^\infty}\rangle$ and there is a map $R_n\to \tilde{S}(n)$. So, for any prism $(B,J)$ over $\tilde{S}(n)$ we get a map $$A_0\langle X_1^{1/p^\infty},\cdots,X_n^{1/p^\infty}\rangle\to B$$
    
    If we apply \cite[Lemma 12.5]{10.4007/annals.2022.196.3.5} to the image of $Y_iY_j$ under this map, we conclude that $\frac{(Y_iY_j)^t}{[\lfloor t\rfloor]_q!}\in B$ and we get a map $\prism_{\tilde{S}(n)}\to B$.
\end{proof}

\begin{dfn}
    Let $S$ be a $R_0$-algebra. If $(S)_{\prism}^{\textup{fl}}$ has a final object, we denote it by $(\prism_{S},\BI_{S})$.
\end{dfn}

\begin{pro}
    If $\tilde{S}\onto S$ is an admissible thickening, then $(\tilde{S})_{\prism}^{\textup{fl}}$ has a final object.
\end{pro}

\begin{proof}
    Since $S$ is an fgpp $R$-algebra, we have the following identifications: \[S \cong S(n) \widehat{\otimes}_{R_n} R \quad \text{and} \quad \tilde{S} \cong \tilde{S}(n) \widehat{\otimes}_{R_n} R\] 
    Define $$\prism_{\tilde{S}}:=\prism_{\tilde{S}(m)}\overset{\doublewedge}{\otimes}_{A_m} A_{\textup{inf}}(R)$$ 
    
    Since $R$ is flat over $R_0$, similarly to the proof of Proposition \ref{initialflat}, one sees that $\prism_{\tilde{S}}$ is the final object of $(\tilde{S})_{\prism}^{\textup{fl}}$. 
\end{proof}

\begin{pro}\label{nygaardfliterationnice}
    In the setting of the previous proposition $\CN^1_{\tilde{S}}$ is the ideal generated by elements $$[p]_{q^{1/p}} \quad \text{and} \quad \frac{(f_if_j)^t}{[\lfloor t\rfloor]_q!}\quad (t\in \mathbb{N},1\le i,j\le n)$$ 
    
    In particular, we have $\prism_{\tilde{S}}/\CN^1_{\tilde{S}}=\tilde{S}$.
\end{pro}

\begin{proof}
    Since $A_{\textup{inf}}(R)$ is a flat $A_r$-algebra and $\sigma:A_{\textup{inf}}(R)\to A_{\textup{inf}}(R)$ is a bijection, we have $\CN^1_{\tilde{S}}=\CN^1_{\tilde{S}(n)}\otimes_{A_r} A_{\textup{inf}(R)}$ (note that $\CN^1_{\tilde{S}}=\kernel(\sigma:\prism_{\tilde{S}}\to \prism_{\tilde{S}}/[p]_q)$). Thus, we are reduced to the explicit computations in Example \ref{universalexample}.
\end{proof}

\begin{dfn}
    Let $S$ be an $R_0$-algebra. If $(S)_{\prism}^{\textup{fl}}$ has a final object $(\prism_{S},\BI_{S})$, we define the groupoid of \emph{prismatic higher $(G,\mu)$-displays} over $S$ to be the groupoid $\mathcal{B}^{\prism_{/(\prism_{S},\BI_S)}}_{(G,\mu)}(\prism_{S},\BI_S)$ (cf. Definition \ref{DefPrimStackPrimDisp}).
\end{dfn}

\begin{rmk}
    By Corollary \ref{equivofetaleandflat}, if $S$ is a qrsp $R_0$-algebra, this is compatible with the definition of prismatic higher $(G,\mu)$-displays over $S$, i.e., objects of $\mathcal{B}^{\textup{qsyn}}_{(G,\mu)}(S)$.
\end{rmk}

\begin{dfn}\label{DefNilAdCat}
    Let $R$ and $S$ be as above. We denote by $\mathfrak{AdNil}_S$ the category whose objects are \emph{admissible thickenings} of $S$. A morphism from $(\tilde{S}_1,f_1)$ to $(\tilde{S}_2,f_2)$ is a morphism $g:\tilde{S}_1\to \tilde{S}_2$ such that $f_2\circ g=f_1$.
\end{dfn}

\begin{dfn}
    Let $S$ be a quasi-syntomic ring. Let $\CD$ and $\SH$ be, respectively, a prismatic higher $(G,\mu)$-display and a $(G,\mu)$-Hodge filtration on $S$. As in Definitions \ref{DefDeformationFunctor} and \ref{DefDeformationFunctorHodge}, we define the deformation functors $\Def_{\CD}$ and $ \Def_{\SH}$ from $\mathfrak{AdNil}_S$ to $\Set$.
\end{dfn}

\begin{rmk}
    Let $S$ be an $R_0$-algebra and assume that $(S)_{\prism}^{\textup{fl}}$ has a final object $(\prism_{S},\BI_{S})$. As in Remark \ref{RemStackofHodgeFilt}, based on Corollary \ref{assoctiedtorsor}, we can attach a $(G,\mu)$-Hodge filtration to a prismatic higher $(G,\mu)$-display $\CD$ over $S$, and denote it by $\FH(\CD)$.
\end{rmk}

\begin{thm}\label{ThmDefAdmThink}
    Let $S$ be a quasi-syntomic ring and $\CD$ an adjoint-nilpotent prismatic higher $(G,\mu)$-display on $S$. Let $f:\tilde{S}\onto S$ be an admissible thickening. The map $\FH:\Def_{\CD}(\tilde{S},f)\to\Def_{\FH(\CD)}(\tilde{S},f)$ is a bijection.
\end{thm}

\begin{proof}
    The proof is similar to the proof for thickenings of quasi-syntomic rings (Theorem \ref{grothendique-messing}).  First, we prove the analogs of Lemmas \ref{Nygaardofperfectlypresented}, \ref{universalthickening} and \ref{topologicallynilpotent} using Proposition \ref{nygaardfliterationnice}. Now consider the frames $$\tilde{\SB}=(\prism_{\tilde{S}},\SN^1_{\tilde{S}},\sigma,\tilde{\sigma}),\quad \text{and} \quad \SB/\tilde{\SB}=(\prism_{S},J,\sigma,\tilde{\sigma})$$ where $J$ is the kernel of the surjection $\prism_S\onto \tilde{S}$ as in Lemma \ref{universalthickening}. We want to prove that $\tilde{\SB}\to \SB/\tilde{\SB}$ is nil-crystalline. Define the frames $\tilde{\SB}/[w]^n$ and $\left(\SB/\tilde{\SB}\right)/[w]^n$ and repeat the proof of Proposition \ref{prismaticuniquelifitng}. The rest follows by descent. For the descent part, note that every \'etale cover of an admissible thickening is again an admissible thickening by Lemma \ref{etalecoverofqrsp}.
\end{proof}

\begin{pro}\label{prismaticdieudonnetheory}
    Let $\tilde{S}\onto S$ be an admissible thickening. There is a natural functor from the category of $p$-divisible groups over $\tilde{S}$ to the category of prismatic displays over $\prism_{\tilde{S}}$. Moreover, the Hodge filtration of the $p$-divisible group is the same as the Hodge filtration of the associated display.
\end{pro}

\begin{proof}
    Let $\mathcal{G}$ be a $p$-divisible group over $\tilde{S}$. Suppose that $v$ is the morphism of topoi:$$v:\textup{Shv}((\tilde{S})_{\prism})\to \textup{Shv}((\tilde{S})_{\textup{QSYN}})$$
    induced by the cocontinuous map that sends the prism $((A,I),R\to A/I)$  to $(R\to A/I)$ (this is continuous by \cite[proposition 7.11]{10.4007/annals.2022.196.3.5} ).

    Define $$M_{\prism}(\mathcal{G}):=\mathcal{E}xt^1_{(\tilde{S})_\prism}(v^{-1}\mathcal{G},\mathcal{O}_{\prism})(\prism_{\tilde{S}})$$
    
    By \cite[Proposition 4.69]{anschutz2023prismatic}, this is a finitely generated projective $\prism_{\tilde{S}}$-module with a Frobenius semi-linear endomorphism coming from the Frobenius on $\mathcal{O}_{\prism}$. Similarly to \cite[Theorem 4.71]{anschutz2023prismatic}, $M_{\prism}(\mathcal{G})$ is an admissible Breuil-Kisin-Fargues module: We have a surjection $R\onto \tilde{S}$ from an integral perfectoid ring $R$ to $\tilde{S}$ with Henselian kernel, and we can lift $\CG$ to a $p$-divisible group $\mathcal{G}_R$ over $R$. Now, $M_{\prism}(\mathcal{G}_R)$ is admissible over $A_{\textup{inf}}(R)$ and by base change we are done.

    By Proposition \ref{char} an admissible Breuil-Kisin-Fargues module over $\prism_{\tilde{S}}$ is the same as a prismatic display over $\prism_{\tilde{S}}$. The statement about Hodge filtration follows from base change and a similar statement in the perfectoid setting (\cite[Theorem 17.5.2]{scholze2020berkeley}).
\end{proof}

\begin{thm}\label{ClasspDivGrps}
    Let $R$ be a flat integral perfectoid $R_0$-algebra and $\tilde{S}$ an $R$-algebra. Assume that there exists an admissible thickening $\tilde{S}\onto S$. There is an equivalence of categories between the category of formal $p$-divisible groups over $\tilde{S}$ and the category of adjoint-nilpotent prismatic displays over $\prism_{\tilde{S}}$.
\end{thm}

\begin{proof}
    By Proposition \ref{prismaticdieudonnetheory} we have a functor from the category of $p$-divisible groups over $\tilde{S}$ to the category of prismatic displays over $\prism_{\tilde{S}}$. We want to show that this is an equivalence of categories when we restrict to formal p-divisible groups. By a standard argument, it is enough to prove this for the underlying groupoids (i.e., ignoring non-invertible morphisms). Indeed, if $f:\mathcal{G}\to \mathcal{G}'$ is a morphism of $p$-divisible groups, one can encode it in the automorphism
    \begin{equation*}
        \begin{bmatrix}
            1 & 0 \\
            f & 1
        \end{bmatrix}
    \end{equation*}
    of $\mathcal{G}\oplus\mathcal{G'}$. For details, see \cite{ito2023deformation}.

    To prove this for the underlying groupoids, note that the following diagram is commutative, where $\CB\CT$ and, respectively, $\mathcal{D}\textit{isp}$ denote the groupoid of formal $p$-divisible groups and, respectively, adjoint-nilpotent prismatic displays. The bottom arrow is an equivalence by \cite{anschutz2023prismatic} and the fibers of the vertical arrows are in bijection by Grothendieck-Messing:
    \[\begin{tikzcd}
        \CB\CT(\tilde{S}) \arrow{r}\arrow{d}& \mathcal{D}\textit{isp}(\tilde{S})\arrow{d}\\
        \CB\CT(S) \arrow{r}& \mathcal{D}\textit{isp}(S)
    \end{tikzcd}\]
\end{proof}

We saw that $\prism_{\tilde{S}}$ is the initial object in the category of prisms over $\tilde{S}$ that are flat over $A_0$. There is a different way of saying this that does not depend on $A_0$ and therefore can be generalized.

\begin{dfn}
    An oriented prism $(A,(d))$ is called \emph{transversal} if $d,p$ is a regular sequence. 
\end{dfn}

\begin{lem}\cite[Lemma 2.1.7]{anschutz2023prismatic}
    If $(A,(d))$ is a transversal prism, then for all $i, j>0$, $i\neq j$, the pair $\phi^i(d),\phi^j(d)$ is a regular sequence.
\end{lem}

\begin{rmk}
    Let $R$ be an integral perfectoid ring that is flat over $R_0$. Let $\tilde{S}$ be an $R$-algebra and assume that there exists an admissible thickening $\tilde{S}\onto S$. Then $\prism_{\tilde{S}}$ is the initial object in the category $(\tilde{S})_{\prism}^{\textup{t}}$ of transversal prisms over $\tilde{S}$ because every such prism has a map from $A_{\textup{inf}(R)}$ and the proof of \cite[Lemma 12.5]{10.4007/annals.2022.196.3.5} only uses the fact that $\sigma^{i}(d),\sigma^j(d)$ is a regular sequence.    
\end{rmk}

\begin{dfn}
    Let $S$ be a ring. We say that $(S)_\prism$ has \emph{enough transversal prisms} if there is a transversal prism in $(S)_\prism$ that weakly covers the final object.
\end{dfn}

\begin{question}
    Let $\tilde{S}\onto S$ be a thickening with $S$ quasi-syntomic. Assume that $(S)_\prism$ has enough transversal prisms. Is there an equivalence of categories between the category of $p$-divisible groups over $\tilde{S}$ and the category of admissible prismatic crystals over $(\tilde{S})_{\prism}^{\textup{t}}$? 
\end{question}

More generally, it is natural to ask whether there is a relation between crystals in Breuil-Kisin-Fargues modules with $G$-structure over $(\tilde{S})_{\prism}^{\textup{st}}$ and $G$-torsors on the stack of $F$-gauges defined by Bhatt-Lurie \cite{bhatt2022prismatic}, at least when considering those objects that are bounded by a minuscule cocharacter $\mu$ in some suitable sense.

\section{Generalizing to General Thickenings}\label{removingcondition}
Let $\textbf{QRSP}^\bullet$ denote the ($\infty$-)category of cosimplicial qrsp rings, i.e., the category of cosimplicial rings $S^\bullet$ such that each $S^i$ is a qrsp ring.

\begin{dfn}\label{simplicialdisplaygroup}
Let $S^\bullet \in \textbf{QRSP}^\bullet$. We define the cosimplicial groups $L^+G(S^\bullet)$ and $G^\mu(S^\bullet)$ by setting $L^+G(S^\bullet)^i := L^+G(S^i)$ and $G^\mu(S^\bullet)^i := G^\mu(S^i)$ for each $i$.
\end{dfn}

\begin{lem}\label{goodresolution}
If $S$ is a qrsp ring and $S \to S^\bullet$ is a resolution, then the natural map
\[
G^\mu(S) \to G^\mu(S^\bullet)
\]
is a quasi-isomorphism.
\end{lem}

\begin{proof}
From \cite[Sections 12.4 and 12.5]{bhatt2021prismatic}, we have a quasi-isomorphism $\prism_S \to \prism_S^\bullet$ and an identification $\CN_S^i \cong \CN_{S^\bullet}^i$.
\end{proof}

\begin{dfn}\label{simplicialdisplay}
Let $S^\bullet \in \textbf{QRSP}^\bullet$. We define the cosimplicial groupoid of banal $(G,\mu)$-displays over $S^\bullet$ as the totalization of the groupoids $\CB^{\textup{ban}}(G,\mu)(S^i)$.
\end{dfn}

\begin{lem}\label{compatiblewithdirectlimit}
Assume that $S' \to S$ is a first-order thickening of qrsp rings, and furthermore that $S'$ is the filtered direct limit of qrsp rings $S'_i$. Let $D \in \CB^{\textup{ban}}(G,\mu)(S)$. Then there is a natural equivalence:
\[
\Def_{\CD}(S', f) \xrightarrow{\sim} \lim \Def_{\CD}(S'_i, f).
\]
\end{lem}
\begin{proof}
First note that, by the Hodge–Tate comparison, the natural maps $\prism_{S'} \to \lim \prism_{S_i'}$ and $\CN^1_{S'} = \lim \prism_{S'_i}$ are isomorphisms. Since $G$ is of finite type, to specify an object or morphism in $\Def_{\CD}(S',f)$ it suffices to give finitely many elements in $\prism_{S'}$ satisfying certain relations. This is equivalent to giving an object or morphism in the limit, as desired.
\end{proof}

\begin{pro}\label{BanalGrothendiquemessinggeneral}
Let $R$ be an integral perfectoid ring, and let $S' \twoheadrightarrow S$ be a first-order thickening of finitely generated $p$-power torsion-free (fgpp) $R$-algebras. Then there is a natural bijection between the isomorphism classes of deformations of an adjoint-nilpotent banal prismatic higher $(G,\mu)$-display along $S' \twoheadrightarrow S$ and the set
\[
\left(G(S') \times_{G(S)} P_{\mu}(S)\right) / P_{\mu}(S').
\]
\end{pro}

\begin{proof}
By \cite[Lemma 4.28]{bhatt2019topological}, one can construct a cosimplicial resolution of $S$ such that each $S^i$ is a quotient of a flat integral perfectoid $R_0$-algebra. There is a natural map from the isomorphism classes of deformations to the quotient
\[
\left(G(S') \times_{G(S)} P_{\mu}(S)\right) / P_{\mu}(S').
\]
Using the associated spectral sequence, it suffices to show that this map is an isomorphism for each $S^i$. This follows from Theorem \ref{grothendique-messing} and Lemma \ref{compatiblewithdirectlimit}.
\end{proof}

\begin{thm}\label{grothendique-messinggeneral}
Let $S$ be a quasi-syntomic ring and $\CD$ an adjoint-nilpotent prismatic higher $(G,\mu)$-display over $S$. Let $(S',f)$ be an object in $\mathfrak{Nil}_S^1$. Then the natural map
\[
\FH\colon \Def_{\CD}(S',f) \to \Def_{\FH(\CD)}(S',f)
\]
is a bijection.
\end{thm}

\begin{proof}
Repeat the proof of Theorem \ref{grothendique-messing}, using Proposition \ref{BanalGrothendiquemessinggeneral} in place of Proposition \ref{BanalGrothendiquemessing}.
\end{proof}

\section{Relation with Shtukas and \texorpdfstring{$\Bun_{G}$}{Bun\_G}}\label{Secrelationshtukas}

We want to define a functor $\textup{Fib}_G$ from $\mathcal{B}^{\perf}_{(G,\mu)}$ to the stack of shtukas. For a review of the theory of shtukas, see Appendix \ref{shtukas}. Let $\Spa(R,R^{+})$ be an affinoid perfectoid. First, note that by Proposition \ref{fargueskisinmodule}, we have a functor from $\mathcal{B}^{\textup{perf}}_{(G,\mu)}$ to the stack of Breuil-Kisin-Fargues modules with $G$-structure. There is also a functor from the category of Breuil-Kisin-Fargues modules with $G$-structure over $\Spa\left(A_{\textup{inf}}(R^+)\right)$ to the category of $G$-shtukas over $R$, by restricting to $\Spa\left(A_{\textup{inf}}(R^+)\right)\setminus V(p,w)$. Composing these two functors, we obtain a functor $$\textup{Fib}_G:\mathcal{B}^{\perf}_{(G,\mu)}\to\Sht_G$$

\begin{example}
    For every affinoid perfectoid $\Spa(R,R^+)$, there is an equivalence of categories between the category of $p$-divisible groups over $\Spa(R,R^+)$ and the category of shtukas with one leg over $\Spa(R,R^+)$ \cite{scholze2020berkeley}. By Example \ref{prismaticdiudennetheory}, for a minuscule cocharacter $\mu$ we also have an equivalence between $(\GL_n,\mu)$-displays and $p$-divisible groups of height $n$ and Hodge polygon $\mu$. This means that we obtain a functor from the category of $(\GL_n,\mu)$-displays to that of shtukas, and by unwinding the definitions, one can see that this is the same functor as $\textup{Fib}_{\GL_n}$.
\end{example}

Let $S$ be an affinoid perfectoid space of characteristic $p$. We define the functor $\pi:\Sht_G\to \Bun_G$ where $\Bun_G$ is the stack of $G$-bundles on the Fargues-Fontaine curve: Pick a real number $\epsilon>0$ such that $x\not\in \mathcal{Y}_{[0,\epsilon],S}$. The Fargues-Fontaine curve $X_S$ is the quotient of $\mathcal{Y}_{(0,\epsilon],S}$ by Frobenius and $F_{|_{\mathcal{Y}_{(0,\epsilon],S}}}$ is an isomorphism. Choose a $G$-shtuka $(\SP,x,\iota)$ over $S$. This map defines a descent datum for the $G$-torsor $\SP_{|_{\mathcal{Y}_{(0,\epsilon],S}}}$. Therefore, we get a $G$-torsor on $X_S$. By composing with $\textup{Fib}_G$, we obtain a functor $\mathcal{B}^{\textup{perf}}_{(G,\mu)}\to \Bun_G$.

We want to study the topological image of this map. First, we give an alternative definition of it on algebraically closed fields. Let $(C,C^+)$ be an algebraically closed perfectoid field, and let $\kappa$ be the residue field of $C^{\flat+}$. The surjection $A_{\textup{inf}}(C^+)=W(C^{\flat+})\onto W(\kappa)$ defines a functor $$\mathcal{B}^{\perf}_{(G,\mu)}(C,C^{+})\to {\textup{Win}}_{(G,\mu)}(\mathscr{W}({\kappa}))$$ 

We claim that this is surjective. Since $C$ is algebraically closed, by the Almost Purity Theorem \cite[Theorem 7.4.5]{scholze2020berkeley} every display over $(C,C^+)$ is banal, and it is enough to prove that the map $G\left(W(C^{\flat+})\right)\to G(W(\kappa))$ is surjective. Since $G$ is smooth and affine over $\mathbb{Z}_p$ and $W(C^{\flat+})$ is $W(\mathfrak{m})$-adically complete ($\mathfrak{m}$ being the maximal ideal of $C^{\flat+}$), by Lemma \ref{FormalSmoothLift} the map $G\left(W(C^{\flat+})\right)\to G\left(W(C^{\flat+})/W(\Fm)\right)=G(W(\kappa))$ is surjective. 

Consider $B(G)$, the set of isocrystals over an algebraically closed field of characteristic $p$. We have a functor ${\textup{Win}}_{(G,\mu)}(\mathscr{W}({\kappa}))\to B(G)$ that sends a display represented by $U$ to the isocrystal represented by $U\mu(p)$. The image of this map is often denoted by $B(G,\mu)$.

\begin{dfn}
    The set $B(G,\mu)$ consists of $G$-isocrystals $[b]$ such that the set \[\left\{g\in G\left(W(\kappa)[\frac{1}{p}]\right)\,|\,g^{-1}bg^{\sigma}\in G\left(W(\kappa)\right)\mu(p)G(W(\kappa))\right\}\] is nonempty.
\end{dfn}

\begin{rmk}
    Equivalently, this consists of classes such that $\kappa_G(b)=\mu$ and $v_b\le \mu$ where $\kappa_G(b)$ and respectively $v_b$ are the Kottwitz and Newton invariants defined in \cite{rapoport1996classification}.
\end{rmk}

\begin{pro}
    The set of topological points in the image of the morphism $\mathcal{B}_{(G,\mu)}^{\textup{perf}}\to \Bun_G$ is the set $B(G,\mu)$.
\end{pro}

\begin{proof}
     Assume that $(C,C^+)$ is an algebraically closed field of characteristic $p$. There is an isomorphism $B(G)\to \Bun_G(C)$ \cite[Theorem 5.3]{anschutz2022gbundles}, so if we can show that the diagram
     \[\begin{tikzcd}
            \mathcal{B}^{\textup{perf}}_{(G,\mu)}(C^\sharp) \arrow{r}\arrow{d}& \Bun_G(C) \\
            {\textup{Win}}_{(G,\mu)}(\mathscr{W}({\kappa}))\arrow{r}& B(G)\arrow{u}\\
        \end{tikzcd}\]
    is commutative, we conclude that the image of $\mathcal{B}^{\textup{perf}}_{(G,\mu)}(C^\sharp)\to \Bun_G$ is the same as the image of the map ${\textup{Win}}_{(G,\mu)}(\mathscr{W}({\kappa}))\to B(G)$ which is $B(G,\mu)$.

    To see that the diagram is commutative, note that $C$ is algebraically closed, so all objects are banal. By the Tannakian formalism, we can reduce to the case $G=\GL_n$. In that case, we send $X\in \mathcal{B}^{\textup{perf}}_{(G,\mu)}(C^\sharp)$ to $b=X\mu(\phi(d))$. By \cite[Lemma 13.4.1]{scholze2020berkeley} there is a $k\in G(\CO_{Y_{[0,\epsilon]}})$ such that $kb\phi(k)^{-1}$ is given by a standard matrix $N$ in $B(G)$. If we denote the image of $b$ in the residue field of the point $0$ by $\bar{b}$, we have $\bar{b}=X\mu(pv)$ for some unit $v$.

    Now we have $\bar{b}=\bar{k}N\sigma(\bar{k})^{-1}$  which means that its image is $N\in B(G)$ as desired.
\end{proof}

One can also introduce a version with a ``fixed Newton polygon''. Consider an affinoid perfectoid $\Spa(R,R^{+})\to \Spd W(\kappa) $. Let $[b]\in B(G,\mu)$ be an isocrystal defined over $\kappa$, choose a representative $b=U\mu(p)$ for $[b]$ and consider the banal Witt display $\bar{\CD}_b$ represented by the image of $U$ under the map $W(\kappa)\to W(R^{\sharp +}/p)$. 

\begin{dfn}
    A \emph{quasi-isogeny} between displays $\CD$ and $\CD'$ is an isomorphism between the $\SL G$-torsors associated with $\CD$ and $\CD'$ that commutes with Frobenius. 
\end{dfn}

\begin{dfn}\label{DefStackGMuB}
    We define the \'etale stack $\mathcal{B}^{\textup{perf}}_{(G,\mu,b)}$ as follows: for an affinoid perfectoid given by $\Spa (R,R^+)\to \Spd W(\kappa)$, objects of $\mathcal{B}^{\textup{perf}}_{(G,\mu,b)}\left(\Spa (R,R^+)\to \Spd W(\kappa)\right)$ are pairs $(\CD,\eta)$ where $\CD\in \mathcal{B}^{\textup{perf}}_{(G,\mu)}(R,R^+)$ and $\eta$ is a quasi-isogeny $\eta:\bar{\CD}\dashrightarrow \bar{\CD}_b$ where $\bar{\CD}$ is the reduction of $\CD$ mod $p$ (see Proposition \ref{ReductionModpPerfVer}). A morphism from $(\CD,\eta)$ to $(\CD',\eta')$ is a morphism $k:\CD\to \CD'$ in $\mathcal{B}^{\textup{perf}}_{(G,\mu)}$ such that $\eta'\circ k=\eta$. 
\end{dfn}

Let $(G,\mu,b)$ be an integral local Shimura datum and let $\mathcal{M}^{\textup{int}}(G,\mu,b)$ be the associated integral local Shimura variety (see Definitions \ref{IntShimDat} and \ref{DefIntLocShVar}). There is a functor $\widetilde{\textup{Fib}}_G:\mathcal{B}^{\textup{perf}}_{(G,\mu,b)}\to {\mathcal{M}^{\textup{int}}(G,\mu,b)}$ defined in \cite{bartling2022mathcal}.

\begin{lem}\label{translationofquasiisogenies}
    A quasi-isogeny $\eta:\bar{\CD}\dashrightarrow\bar{\CD_b}$ lifts uniquely to an isomorphism $\tilde{\eta}:\textup{Fib}_G(\CD)_{[r,\infty)}\to \SP^b_{[r,\infty)}$ for some $r>0$.
\end{lem}

\begin{proof}
    Given that $B_{\textup{crys}}^+(R)$ is a subring of $H^0(\mathcal{Y}_{[\frac{1}{p^{p-1}},\infty)},\CO_{\mathcal{Y}_{[\frac{1}{p^{p-1}},\infty)}})$, it suffices to construct an isomorphism $\textup{Fib}_G(\CD)_{B_{\textup{crys}}^+}\to \SP^b_{B_{\textup{crys}}^+}$. This is \cite[Lemma 8.1]{bartling2022mathcal}.
\end{proof}

It follows from the lemma that we have a well-defined fully faithful functor $$\widetilde{\textup{Fib}}_G:\mathcal{B}^{\textup{perf}}_{(G,\mu,b)}\to {\mathcal{M}^{\textup{int}}(G,\mu,b)}$$

\begin{thm}\label{localshimuralocalperfectoid}
    The functor $\widetilde{\textup{Fib}}_G:\mathcal{B}^{\textup{perf}}_{(G,\mu,b)}\to {\mathcal{M}^{\textup{int}}(G,\mu,b)}$ is an equivalence of categories.
\end{thm}

\begin{proof}
     It is enough to prove that it is essentially surjective. This is \cite[Remark 25]{bartling2022mathcal}.
\end{proof}

\section{Relation to Integral Local Shimura Varieties}\label{Secrelationlocalshimuravariety}

In this section, we study the relation between the stack of prismatic higher $(G,\mu)$-displays and integral local Shimura varieties (see Appendix \ref{shtukas} for details). Consider the category $\textbf{\textup{FQSyn}}$ of $p$-torsion-free quasi-syntomic rings. This is a site with the quasi-syntomic topology. As we saw, there is a reduction map $\mathcal{B}^{\textup{qsyn}}_{(G,\mu)}(R)\to {\textup{Win}}_{(G,\mu)}(\mathscr{W}(R/p))$ that sends $\CD$ to $\bar{\CD}$ (cf. Proposition \ref{ReductionModp}). Consider an integral local Shimura datum $(G,\mu,b)$.

\begin{dfn}
    Consider the groupoid $\mathcal{B}^{\textup{qsyn}}_{(G,\mu,b)}$ over $\textbf{\textup{FQSyn}}$ defined as follows: objects of $\mathcal{B}^{\textup{qsyn}}_{(G,\mu,b)}(R)$ are pairs $(\CD,\eta)$ where $\CD\in \mathcal{B}^{\textup{qsyn}}_{(G,\mu)}(R)$ and $\eta$ is a quasi-isogeny $\eta:\bar{\CD}\dashrightarrow \bar{\CD}_b$ between the reductions mod $p$. A morphism between $(\CD,\eta)$ and $(\CD',\eta')$ is a morphism $k:\CD\to \CD'$ in $\mathcal{B}^{\textup{qsyn}}_{(G,\mu)}$ such that $\eta'\circ k=\eta$. 
\end{dfn}

\begin{pro}
    There is a natural transformation $$\mathfrak{Fib}_G:\mathcal{B}^{\textup{qsyn}}_{(G,\mu,b)}\to{\mathcal{M}^{\textup{int}}(G,\mu,b)}$$ on $\textbf{\textup{FQSyn}}$.
\end{pro}

\begin{proof}
   By definition, to have an element of ${\mathcal{M}^{\textup{int}}(G,\mu,b)}(R)$, for every affinoid perfectoid $S$ of characteristic $p$, and every untilt of $S$ over $R$ like $R\to S^\sharp$, we should give a $G$-Shtuka $\SP$ on $\mathcal{Y}_S$ with a leg at $S^\sharp$ bounded by $\mu$, and an isomorphism between $\SP_{[r,\infty)}$ and $\SP^b_{[r,\infty)}$. 
   
    By functionality, the map $R\to S^\sharp$ induces a map $$\mathcal{B}^{\textup{qsyn}}_{(G,\mu)}(R)\to \mathcal{B}^{\textup{qsyn}}_{(G,\mu)}(S^\sharp)=\mathcal{B}^{\textup{perf}}_{(G,\mu)}(S^\sharp)$$ 
    
    So, $(\CD,\eta)\in \mathcal{B}^{\textup{qsyn}}_{(G,\mu,b)}(R)$ gives an element of $\mathcal{B}^{\textup{perf}}_{(G,\mu)}(S^\sharp)$ still denoted by $\CD$, and $\SP=\textup{Fib}_G(\CD)$ is a $G$-shtuka on $S$ bounded by $\mu$. By Lemma \ref{translationofquasiisogenies}, one can lift $\tau$ to an isomorphism $\SP_{[r,\infty)}\to \SP^b_{[r,\infty)}$. This defines the natural transformation.
\end{proof}

\begin{pro}\label{perfectoidconjecturerightversion}
    Over integral perfectoid rings, the functor $\mathfrak{Fib}_G$ is an isomorphism.
\end{pro}

\begin{proof}
    The proof is similar to \cite[Remark 25]{bartling2022mathcal}, and we only give a rough sketch of it. The main step is to prove that if $(\SP,\tau)$ is a shtuka over $\Spa(R,R^+)$ coming from a display over $R^+$, then one can extend the vector bundle $\SP$ from $\Spa W(R^+)\setminus V((p,[w]))$ to $W(R^+)$. This extension is unique if it exists because $p,[w]$ is a regular sequence, and $\mathfrak{Fib}_G$ is fully faithful on perfectoid rings. Therefore, by descent, it is enough to extend the vector bundle locally. For valuation rings, this is done in \cite[Theorem 2.7]{kedlaya2020some}. In general, for every point $\Spa(k,k^+)$ on $\Spa(R^+,R^+)$, by the argument given in \cite[Remark 25]{bartling2022mathcal} one can extend $\SP$ to a vector bundle over $W(k^+)$ and then spread this out to an open neighborhood, the key ingredient being the generalization of Bouville-Laszlo Theorem \cite[\S3, Théorème]{zbMATH00810980} given in \cite[Proposition 5.6]{zbMATH06670053}.
\end{proof}

\begin{pro}\label{PropFib_GSurj}
    If $\mathcal{M}^{\textup{int}}(G,\mu,b)$ is representable by the diamond of a formal scheme, then for any $R\in \textbf{\textup{FQSyn}}$, the map $$\mathfrak{Fib}_G(R):\mathcal{B}^{\textup{qsyn}}_{(G,\mu,b)}(R)\to{\mathcal{M}^{\textup{int}}(G,\mu,b)}(R)$$ is surjective.
\end{pro}

\begin{proof}
    It is enough to prove this over a basis, so we can assume that $R$ is qrsp. There is a surjective map from a perfectoid ring $\tilde{R}$ to $R$ with the kernel $J$. By \cite[Corollary 2.10]{anschutz2023prismatic} the completion of the henselisation of $\tilde{R}$ at $J$ is again perfectoid. So, one can assume that this map is Henselian along its kernel. We know that $\mathcal{M}^{\textup{int}}(G,\mu,b)$ is smooth. Therefore, by Lemma \ref{FormalSmoothLift}, we can lift every point of $\mathcal{M}^{\textup{int}}(G,\mu,b)(R)$ to $\tilde{R}$, which means that we can reduce the question to perfectoid rings.
\end{proof}

\begin{lem}\label{existanceofliftsHenselian}
    Let $(A,I)\to (B,IB)$ be a surjective map of prisms and assume that $A$ is Henselian with respect to its kernel. Set $R:=A/\CN_A^1$ and $R'=B/\CN_B^1$. Let $\CD$ be a prismatic higher $(G,\mu)$-display on $(B,IB)$ with $(G,\mu)$-Hodge filtration $\SL$ and let $\tilde{\SL}$ be a lift of $\SL$ to $R$. Then $\CD$ has a lift $\tilde{\CD}$ to $(A,I)$ with $(G,\mu)$-Hodge filtration $\tilde{\SL}$.
\end{lem}

\begin{proof}
    By the Tannakian formalism, we can reduce to the case of $\GL_n$. In that case, we can lift $\tilde{\SL}$ to a decomposition on $A$ because $A$ is Henselian with respect to the kernel of $A\to R$ (Proposition \ref{Henselianpair}). This decomposition is a lift of some normal decomposition of $B$ because it has the right image in $R'$. So, by the definition of displays, we only need to lift the divided Frobenius ${\textup{div}}\,F$ to $B$. But by assumption, $A$ is Henselian with respect to the kernel of the map $A\to B$ and so we can lift every isomorphism between projective modules on $B$ to $A$, including ${\textup{div}} F$.
\end{proof}

\begin{thm}\label{equivalencetolocalshimura}
    If the integral local Shimura variety is representable by the diamond of a formal scheme, then the natural transformation $\mathfrak{Fib}_G$ is an isomorphism.
\end{thm}

\begin{proof}
    By Proposition \ref{PropFib_GSurj}, it suffices to show that $\mathfrak{Fib}_G$ is injective. It is enough to verify this over a basis, and so we can restrict to fgpp rings. Assume that $\CD$ and $\CD'\in \mathcal{B}^{\textup{qsyn}}_{(G,\mu,b)}(S)$ map to the same object $x$ in ${\mathcal{M}^{\textup{int}}(G,\mu,b)}(S)$, in particular they have the same Hodge filtration $\SL$.

    We want to lift these displays to displays on some perfectoid ring with a common Hodge filtration. Our strategy is the same as \cite[Section 4.9]{anschutz2023prismatic}.

    First choose a surjection $f:S(m) \onto S$. By \cite[Lemma 4.86]{anschutz2023prismatic} the pairs $(S(m),\kernel f)$ and $(\prism_{S(m)},\kernel \prism_f)$ are Henselian. By Lemma \ref{existanceofliftsHenselian}, we can lift $\CD$ and $\CD'$ to displays on $S(m)$ with Hodge filtration $\tilde{\SL}$. They still have the same image in $\mathcal{M}^{\textup{int}}(G,b,\mu)$. We know that $S(m)$ is flat over $\mathbb{Z}_p$, therefore it is enough to check that they have the same image on the generic fiber $\mathcal{M}(G,b,\mu)$. The map $\mathcal{M}(G,b,\mu)\to \Gr^{G,\mu}$ is \'etale, and $\CD$ and $\CD'$ have the same image in $\mathcal{M}(G,b,\mu)(S)$ and the same Hodge filtration. By Lemma \ref{formaletaleHenselian}, $\CD$ and $\CD'$ have the same image in $\mathcal{M}(G,b,\mu)\left(S(m)\right)$.

    Therefore, the problem is reduced to the case of $S(m)$. By the proof of \cite[Lemma 4.87]{anschutz2023prismatic} there is a surjection from some perfectoid ring $R\onto S(m)$. The proof also shows that there is some prism $B$, such that $\prism_R\to B$ is Henselian and there is a crystalline map $\prism_S\to B$. Therefore, by Lemma \ref{existanceofliftsHenselian}, one can find lifts $\tilde{\CD}$ and $\tilde{\CD'}$, over $R$, of $\CD$ and $\CD'$, respectively, with Hodge filtration $\tilde{\SL}$. Taking the henselization along the map $R\onto S(m)$, we can assume that this map is Henselian with respect to its kernel, and as before $\tilde{\CD}$ and $\tilde{\CD'}$ have the same image in $\mathcal{M}^{\textup{int}}(G,b,\mu)(R)$. But by assumption, the map $\mathcal{B}^{\textup{perf}}_{(G,\mu,b)}(R)\to \mathcal{M}^{\textup{int}}(G,b,\mu)(R)$ is an equivalence and therefore $\tilde{\CD}=\tilde{\CD'}$ which means that $\CD$ and $\CD'$ are equal.
\end{proof}

\begin{cor}
    Assume that $(G,\mu)$ is of Hodge type, then $\mathfrak{Fib}_G$ is an isomorphism.
\end{cor}

\begin{proof}
    Under the assumptions, the representability of ${\mathcal{M}^{\textup{int}}(G,\mu,b)}$ is known \cite{pappas2022p}.
\end{proof}

We believe that the theorem should be true unconditionally:

\begin{conjecture}\label{ConjFibEquiv}
    The natural transformation $\mathfrak{Fib}_G:\mathcal{B}^{\textup{qsyn}}_{(G,\mu,b)}\to{\mathcal{M}^{\textup{int}}(G,\mu,b)}$ is an isomorphism on $\textbf{\textup{FQSyn}}$.
\end{conjecture}

\begin{rmk}
    Of course, the assumption of the theorem is expected to be true; however, we think it might be possible to prove the conjecture without any recourse to the other conjectures!
\end{rmk}

\begin{appendices}

\section{The Structure of the Display Group}\label{proofofstructure}

Consider a prism $(A,I)$ and a display datum $(G,\mu)$ over $\BZ_q$. We want to understand the geometry of the display group $G^\mu$ (Definition \ref{DefDispGrp}). The following construction is essentially the one in \cite{bueltel2020g} and \cite{lau2021higher}. Since the two definitions differ slightly, we briefly recall the construction. Let us denote by $U_{\mu^{-1}}$ and $P_{\mu}$ the unipotent and, respectively, parabolic subgroups of $G_{\BZ_q}$ given by $\mu^{-1}$, and respectively, $\mu$. There is a scheme-theoretic isomorphism $\exp:\textup{Sym}\left(\Lie U_{\mu^{-1}}\right)\to U_{\mu^{-1}}$ which respects the action of $\BG_m$ on both sides and induces the identity on the Lie algebras. When $\mu$ is minuscule, such a morphism is unique (for a proof, see Lemmas 6.1.1 and 6.3.2 of \cite{lau2021higher}). In this case, it is also a group isomorphism: Set $h(b)=\exp(ab)-\exp(a)$, which is an isomorphism that respects the action of $\mu$, and hence by uniqueness, we should have $h(b)=\exp(b)$. The inverse of $\exp$ is $\log:U_{\mu^{-1}}\to \textup{Sym}\left(\Lie U_{\mu^{-1}}\right)$. The cocharacter $\mu^{-1}$ induces an action of $\BG_m$ on $\Lie\, G_{\BZ_q}$ by adjoint representation. Let $w_1,w_2,\cdots,w_n$ be the weights of this action.

\begin{rmk}
    We always have $P_\mu(A)\subset G^\mu(A,I)$: an element of $P_\mu(A)\subset G(A)$ is by definition a map $g:\CO_{G}\otimes_{\BZ_p}\BZ_q\to A$ whose restriction to $\CF^0$ is zero (using notations of (\ref{EqGradDispGrp})), and therefore trivially respects the filtration on both sides.
\end{rmk}

\begin{rmk}\label{minusculecase}
    Assume that $\mu$ is minuscule. Let $\pi:G(A)\to G(A/\CN_A^1)$ be the reduction mod $\CN_A^1$. One can see that $$G^\mu(A,I)=\pi^{-1}\left(P_\mu(A/\CN_A^1)\right).$$ 
\end{rmk}

\begin{lem}
    Consider the subgroup $U_{\mu^{-1}}^\mu(A,I)=U_{\mu^{-1}}(A)\, \cap\, G^\mu(A,I)$ of $G^\mu(A,I)$. We have a natural bijection $$G^\mu(A,I)\cong P_{\mu}(A)\times U_{\mu^{-1}}^\mu(A,I)$$
\end{lem}

\begin{proof}
    We have a natural map $P_{\mu}(A)\times U_{\mu^{-1}}^\mu(A,I)\to G^\mu(A,I)$. It is injective because, in general, $P_{\mu}\times U_{\mu^{-1}}\to G$ is an open immersion. We only need to show that it is surjective.
 
    Let $\pi:G(A)\to G(A/\CN_A^1)$ be the projection map. Note that $H=P_{\mu}\times U_{\mu^{-1}}\subset G$ is open and $\CN_A^1$ is a subset of $\textup{Rad}(A)$, which means that $\Spec A/\CN_A^1$ contains all closed points of $\Spec A$. Therefore, if we have a map $g:\Spec A\to G$ whose restriction to $\Spec A/\CN_A^1$ factors through $H$, $g$ itself should factor through $H$. In other words $\pi^{-1}\left(H(A/\CN_A^1)\right)\subset H(A)$. 
 
    If $g\in G^\mu(A,I)$, then the image of $g$ in $G(A/\CN^1A)$ lies in $P_{\mu}(A/\CN_A^1)$ (because by definition $g$ sends $\CF^1$ to $\CN^1_A$). We can write $g=xv$ where $x\in P_{\mu}(A)$ and $v\in U_{\mu^{-1}}(A)$. We also have $v=gx^{-1}\in G^\mu(A,I)$, and so $v\in U_{\mu^{-1}}^\mu(A,I)$.
\end{proof}

\begin{pro}\label{structure}
    As a set $G^\mu(A,I)$ is in a canonical bijection with the set $$P_{\mu}(A)\times \log \left(\prod \CN_A^{w_i}\right)$$
\end{pro}

\begin{proof}
    We have to show that $U_{\mu^{-1}}^\mu(A,I)=\log(\prod \CN^{w_i}_A)$. Write $U_{\mu^{-1}}=\Spec S^{<0}$ and consider the filtration on $S^{<0}$ induced by the action of $\BG_m$ on $S^{<0}$. By definition, giving an element of $U_{\mu^{-1}}^\mu(A,I)$ is equivalent to giving a morphism $v:S^{<0}\to A$ that respects the filtration on both sides. This is equivalent to giving a map $v \circ \exp:\textup{Sym}\left(\Lie U_{\mu^{-1}}\right)\to A$, which respects the filtration on both sides.

    Now $\textup{Sym}\left(\Lie U_{\mu^{-1}}\right)$ is a polynomial ring $\BZ_q[x_1,...,x_d]$ and the filtration comes from a grading where $\deg(x_i)=w_i$. Therefore, as sets, we have an equality $$\exp\left(U_{\mu^{-1}}^\mu(A,I)\right)=\left(\prod \CN_A^{w_i}\right)$$
\end{proof}

\begin{rmk}\label{grouphomorphism2}
    If $\mathscr{F}=(A,I,\sigma,\dot{\sigma})$ is a frame and $\mu$ is minuscule, then one obtains a bijection $$G^\mu(\SF)\cong P_\mu(A)\times \prod I^{w_i}$$ where $I^0=A$ and $I^1=I$. 
\end{rmk}

The following definition and remark are basically given in \cite{lau2021higher} and \cite{bultel2020displays} with slightly different notations.

Recall that $P_\mu$ is defined as the subgroup of $G_{\BZ_q}$ consisting of those $g\in G$ with the property that the map $\BG_m\to G_{\BZ_q}$ defined by $t\mapsto\mu(t)g\mu(t^{-1})$ extends to a (necessarily unique) map $\BG_a\to G_{\BZ_q}$.
So, the action of $\BG_m$ on $P_\mu$ defined by $t\cdot g=\mu(t)g\mu(t^{-1})$ extends to an action $\textup{Int}_\mu$ of $\BG_a$ on $P_\mu$.

\begin{dfn}\label{frameaction}
    Let $\mathscr{F}=(A,I,\sigma,\dot{\sigma})$ be a frame, with frame constant $\theta$. Assume that $\mu$ is minuscule. By the above remark, to give a map $G^\mu(\mathscr{F}) \to G(A)$ it is enough to give three maps $P_\mu(A)\to G(A)$, $A\to A$ and $I\to A$. We define the map $\Phi_\theta^\mu:G^\mu(\mathscr{F}) \to G(A)$ to be the one defined as $\sigma\circ\textup{Int}_{\mu}(\theta)$ on $P_\mu$, $\sigma$ on $A$ and $\dot{\sigma}$ on $I$. Assume that this map is a group homomorphism. Then there is a natural action of the group $G^\mu(\SF)$ on $G(A)$, called \emph{$\Phi^\mu$-conjugation}, given as follows: for $k\in G^\mu(\SF)$ and $X\in G(A)$, we set $k\cdot X:=k^{-1}X\Phi_\theta^\mu(k)$. This action is used to define the groupoid of banal windows (cf. Definition \ref{DefBanWind}).
\end{dfn}

\begin{rmk}\label{grouphomorphism1}
    Assume that $\theta$ is a non-zero-divisor in $A$. Then for any $k\in G^\mu(\SF)$, elements $\Phi_\theta^\mu(k)$ and $\mu(\theta)\sigma(k)\mu(\theta)^{-1}$ coincide, and therefore $\Phi_\theta^\mu$ is a group isomorphism. In general, assume that for each pair of elements $k_1$ and $k_2$ of $G^\mu(\SF)$ there is a frame $\SG=(B,J,\tau,\dot{\tau})$, with a non-zero-divisor frame constant $\eta$, a frame morphism $f:\SG\to \SF$ and elements $\tilde{k}_1,\tilde{k}_2\in B$ such that $f(\tilde{k_i})=k_i$ (for $i=1,2$), then $\Phi_\theta^\mu$ is a group isomorphism. This is the case, for example, for prisms and Witt frames. 
\end{rmk}

\section{Elements from Perfectoid Geometry}\label{shtukas}

We want to review some definitions from the theory of perfectoid spaces. The main reference is \cite{scholze2020berkeley}.

\begin{dfn}
    Let $\Perf$ be the site of perfectoid spaces in characteristic $p$ with the \'etale topology.
\end{dfn}

Let us recall the definition of the v-topology.

\begin{dfn}
    We say that a collection of maps $\{f_i:X_i\to Y\}_{i\in I}$ of perfectoid spaces is a \emph{v-cover} if for every quasi-compact open subset $U\subset Y$ there exists a finite subset $J\subset I$ and quasi-compact open subsets $V_j\subset X_j\, ({j\in J})$ such that $U= \bigcup_{j\in J}f_j(V_j)$. 
\end{dfn}

\begin{dfn}
    Let $\Perfd$ be the site of perfectoid spaces (of any characteristic) with the v-topology. Sheaves on this site are called \emph{v-sheaves}. There is a functor from the category of formal schemes over $\mathbb{Z}_p$ to the category of v-sheaves, which is fully faithful when we restrict it to the normal and flat formal schemes of finite type over $\mathcal{O}_E$ where $\mathcal{O}_E$ is the ring of integers of a finite extension $E$ of $\mathbb{Q}_p$.
\end{dfn}

\begin{dfn}
    Let $X$ be a formal scheme over $\mathbb{Z}_p$. The \emph{diamond} of $X$, denoted by $X^\diamondsuit$, is the v-sheaf on $\Perf$ that sends an object $S$ to the set of isomorphism classes of pairs $(S^\sharp, S^\sharp\to X)$, where $S^\sharp$ is an untilt of $S$ and $S^\sharp\to X$ is a morphism of adic spaces (to see that this is actually a v-sheaf, see \cite[Lemma 15.1]{scholze2017etale}.)
\end{dfn}

\begin{dfn}
    Let $I$ be a Cartier divisor of a uniform adic space $X$ with support $Z$ and $U=X\setminus Z$. Let $\mathcal{F}$ be an analytic sheaf on $X$, then a section of $\mathcal{F}(U)$ is called \emph{meromorphic} if it extends to an element of $H^0(X,\mathcal{F}\otimes \varinjlim_n I^{\otimes -n})$, where for a positive $n$, we set $I^{\otimes -n}=\underline{\Hom}(I^{\otimes n},\mathcal{O}_X)$. A \emph{meromorphic} isomorphism between sheaves $\mathcal{F}$ and $\mathcal{F}'$ away from $Z$ is a meromorphic section of the sheaf $\underline{\Isom}(\mathcal{F},\mathcal{F}')$.
\end{dfn}

Let $\Spa(R,R^+)$ be an affinoid perfectoid space of characteristic $p$. Fix a pseudo-uniformizer $\omega$. Set $$\mathcal{Y}_R:=\Spa A_{\textup{inf}}(R^+)\setminus\{p,[\omega]\}$$ 

The action of Frobenius on $\mathcal{Y}_R$ is free and the quotient adic space, $\mathcal{X}_R$, is the Fargues-Fontaine curve. Recall that there is a radius function $\kappa:\mathcal{Y}_R\to [0,\infty]$ that sends a valuation $|.|_v$ to $\frac{|w|_v}{|p|_v}$. Let $\mathcal{Y}_{[a,b],R}$ be the inverse image of the interval $[a,b]$ under $\kappa$. If $\mathcal{F}$ is an analytic sheaf on $\mathcal{Y}_R$, we set $\mathcal{F}_{[a,b]}=\mathcal{F}_{\vert_{\mathcal{Y}_{[a,b],R}}}$. If $S$ is a perfectoid space of characteristic $p$, one can define $\mathcal{Y}_S$ by taking an affinoid cover $S=\bigcup \Spa (R_i,R_i^+)$ for $S$ and gluing $\mathcal{Y}_{R_i}$. For more details, see \cite[Section 2.1.1]{fargues2021geometrization}.

Let $S=\Spa(C,C^+)$ be a geometric point. Then every point $x\in \mathcal{Y}_S$ defines an untilt $\Spa(C(x),C(x)^+)$ of $S$. We denote by $B_{dR}^+(C(x))$ the completion of $A_{\inf}(C)[1/p]$ with respect to the adic topology defined by $\kernel\theta_{C(x)}$. It is a complete discrete valuation ring with residue field $C(x)$. We denote its fraction field by $B_{dR}(C(x))$.

\begin{dfn}\label{dfnshtuka}
    Let $S$ be an affinoid perfectoid space of characteristic $p$, and let $G$ be a reductive group over $\mathbb{Z}_p$. A \emph{$G$-shtuka} over $S$ with \emph{leg $x$} is a triple $(\SP,\{x\},\iota)$, where
    \begin{itemize}
        \item $\SP$ is a $G$-torsor over $\mathcal{Y}_{(0,\infty],S}$
        \item $x$ is an untilts of $S$, called \emph{the leg of the shtuka} and $\Gamma_{x}$ is the Cartier divisor of $\mathcal{Y}_{(0,\infty],S}$ associated with $x$.
        \item $\iota:\sigma^*\SP\to \SP$ is a meromorphic isomorphism away from $\Gamma_{x}$.
    \end{itemize}

    The moduli space of $G$-shtukas is denoted by $\Sht_G$.
\end{dfn}

\begin{rmk}
    By \cite[Proposition 19.5.3.]{scholze2020berkeley} $\Sht_G$ is actually a v-sheaf.
\end{rmk}

Let $(\SP,\{x\},\iota)$ be a $G$-shtuka over $S$ and $s$ be a geometric point of $S$. Let $\Spa(C(x),C(x)^+)$ be the untilt of $s$ associated with $x$. For every trivialization of $\SP_s$ on $B_{dR}(C(x))$ one can see $\iota$ as an element of $G\left(B_{dR}(C(x))\right)$.

\begin{dfn}\label{relativposition}
    We say that a $G$-shtuka $(\SP,\{x\},\iota)$ is \emph{bounded by $\mu$} if at every geometric point $s$ and every trivialization of $\SP_s$ at $B_{dR}(x)$ we have $$\iota\in G\left(B_{dR}^+(C(x)\right)\mu(d)G\left(B_{dR}^+(C(x)\right)$$
    (where $d$ is a generator of $\kernel\theta_{R^+}$). The moduli space of $G$-shtukas bounded by $\mu$ is denoted by $\Sht_{(G,\mu)}$.
\end{dfn} 

We also need the definition of the integral local Shimura variety.

\begin{dfn}
\label{IntShimDat}
    An \emph{integral local Shimura datum} is a triple $({G},\mu,b)$, where ${G}$ is a reductive group over $\mathbb{Z}_p$, $\mu$ is a cocharacter of ${G}$ defined over an unramified extension $\BZ_q$ of $\mathbb{Z}_p$, and $b$ is an element of $B({G},\mu)$.
\end{dfn}

\begin{dfn}
    Let $(\SP^b,\emptyset,\iota^b)$ be the shtuka where $\SP^b$ is the trivial $G$-torsor on $\CY_S$, and $\iota^b$ is given by $b^\sigma$.
\end{dfn}

\begin{dfn}\label{DefIntLocShVar}
    Let $(G,\mu,b)$ be an integral local Shimura datum. The associated \emph{integral local Shimura variety}, denoted by $\mathcal{M}^{\textup{int}}(G,\mu,b)$ and defined over $\Spd \mathcal{O}_E$, is the v-sheaf on $\Perfd$ that sends an untilt of $S=\Spa(R,R^+)$ like $S^{\sharp}=\Spa(R^{\sharp},R^{\sharp+})$ to the set of isomorphism classes of triples $(\SP,\tau,i)$ where $(\SP,\tau)$ is a $G$-shtuka on $\Spa(R,R^+)$ with a leg at $\Gamma_{S^\sharp}$ bounded by $\mu$, and $i:\SP_{[r,\infty)}\to \SP^b_{[r,\infty)}$ is a Frobenius-equivariant isomorphism. 
\end{dfn}

\begin{rmk}
    By \cite[Theorem 11]{gleason2020specialization}, $\Gr_{\CO_E}^{(G,\mu)}$ is the local model for the Shimura datum $(G,\mu,b)$, and it has the same moduli interpretation as our Definition \ref{relativposition}. This means that our definition of integral local Shimura variety is consistent with the definition in \cite{scholze2020berkeley}.
\end{rmk}

\begin{rmk}
    The completed local ring of $\mathcal{M}^{\textup{int}}(G,\mu,b)$ at any point $x$ is isomorphic to the completed local ring of $\Gr_{\CO_E}^{(G,\mu)}$ \cite[Theorem 5.3.3]{ito2023deformation}. In particular, the integral local Shimura variety is formally smooth.
\end{rmk}

\section{Henselian Pairs}\label{SecHenselianpair}

In this section, we review the definition and some lifting properties of Henselian pairs.

\begin{dfn}
    A pair $(A,I)$ where $A$ is a ring and $I\subset A$ is an ideal is called \emph{Henselian} if $I$ is contained in the Jacobson radical of $A$ and for any monic polynomial $f\in A[T]$, and factorization $\bar{f}=g_0h_0$ in $A/I[T]$ to monic polynomials $g_0$ and $h_0$ that generate the unit ideal, there exists a factorization $f=gh$ in $A[T]$ with $g$ and $h$ monic and $g_0=\bar{g},h_0=\bar{h}$.
\end{dfn}

\begin{lem}{\cite[\href{https://stacks.math.columbia.edu/tag/09XK}{Lemma 09XK}]{stacks-project}}\label{LemHenselQuot}
    Let $(A,I)$ be a Henselian pair and $A\onto B$ be a surjective ring homomorphism. Then $(B,IB)$ is a Henselian pair.
\end{lem}

\begin{lem}{\cite[\href{https://stacks.math.columbia.edu/tag/0DYD}{Lemma 0DYD}]{stacks-project}}\label{LemHenselSubId}
    Let $(A,J)$ be a Henselian pair and let $I\subset J$ be an ideal. Then $(A,I)$ is also a Henselian pair.
\end{lem}

\begin{lem}\label{FormalSmoothLift}
    Let $R$ be a ring, $A$ an $R$-algebra, and $I\subset A$ an ideal. Let $X$ be a smooth $R$-scheme.
    \begin{itemize}
        \item[(1)] If $(A,I)$ is a Henselian pair, then the canonical map \[X(A)\to X(A/I)\] is surjective.
        \item[(2)] If $A$ is $I$-adically complete and $X$ is affine, then the canonical map \[X(A)\to \varprojlim X(A/I^n)\] is bijective.
    \end{itemize}
\end{lem}

\begin{proof}
    When $X$ is affine, statement (I) is \cite[Th\'eor\`eme I.8]{Gruson1972}, and it follows from  \cite[\href{https://stacks.math.columbia.edu/tag/00T4}{Lemma  00T4}]{stacks-project}, \cite[\href{https://stacks.math.columbia.edu/tag/09XI}{Lemma  09XI}]{stacks-project} and \cite[\href{https://stacks.math.columbia.edu/tag/07M7}{Lemma 07M7}]{stacks-project}. For general $X$, this is \cite[Proposition 6.1.1.(a)]{Cesnavicius} or for quasi-separated schemes, shown in the proof of \cite[Proposition A.0.4]{elmanto2020modules}. Injectivity in (ii) follows from the fact that $X$ is affine and $A$ is separated. Surjectivity follows from the fact that the system $X(A/I^n)$ is Mittag-Leffler by (i) and \cite[\href{https://stacks.math.columbia.edu/tag/0598}{Lemma 0598}]{stacks-project}. 
\end{proof}

\begin{lem}\label{DefFib}
    Let $R$ be a ring, $A$ an $R$-algebra, and $I\subset J\subset J$ be ideals. Assume that $(A,I)$ is a Henselian pair. Let $H$ be an affine smooth group scheme over $R$ and $P\subset H$ be a closed subgroup. The projection $H(A)\to H(A/I)$ induces a bijection on the quotients
\begin{multline*}
\left(H(A)\times_{H(A/J)}P(A/J)\right)/\left(H(A)\times_{H(A/I)}P(A/I)\right)\overset{\cong}{\longrightarrow}\\
    \left(H(A/I)\times_{H(A/J)}P(A/J)\right)/P(A/I)
\end{multline*}

\end{lem}

\begin{proof}
    Injectivity is straightforward. Surjectivity follows from surjectivity of $H(A)\to H(A/I)$ given by the previous lemma.
\end{proof}

\begin{lem}\cite[\href{https://stacks.math.columbia.edu/tag/0D4A}{Lemma  0D4A}]{stacks-project}\label{liftofprojectivemodules}
    If $(A,I)$ is a Henselian pair, the functor $P\to P/IP$ defines a bijection between the set of isomorphism classes of finitely generated projective $A$-modules and finitely generated projective $A/I$-modules.
\end{lem}

\begin{lem}\label{formaletaleHenselian}
    Let $(A,I)$ be a Henselian pair and $f:Y\to X$ be an \'etale and separated morphism. Then, for any $\bar{y}\in Y(A/I)$, $f:Y(A)\to X(A)$ induces a bijection between lifts of $\bar{y}$ to $Y(A)$ and lifts of $f(\bar{y})\in X(A/I)$ to $X(A)$.
\end{lem}

\begin{proof}
    Choose a map $\Spec A\to X$ that lifts $f(\bar{y})$. By changing the base along this map, we can assume that $X=\Spec A$. So, we are reduced to showing that the natural map \[\Hom_A(\Spec A,Y)\to\Hom_A(\Spec (A/I),Y)\] is a bijection. Surjectivity follows from Lemma \ref{FormalSmoothLift}. Injectivity follows from \cite[Corollary 3.13]{7c7cb5a6-e464-3506-b185-42a30a3f8fde}, noting that every connected component of $\Spec A$ has a closed point, and so $\Spec (A/I)$ meets all connected components of $\Spec A$.
\end{proof}

\end{appendices}

\subsection*{Acknowledgements}

We are grateful to Oliver B\"ultel, Simon H\"aberli, and Aliakbar Hosseini for many helpful discussions. We thank Kiran Kedlaya for his interest and for answering some of our questions, and Torsten Wedhorn for suggesting a reference.


\bibliography{myref}{}
\bibliographystyle{alphaurl}


\end{document}